\documentclass{amsart}

\usepackage{amssymb}
\usepackage{eepic}
\usepackage{color}

\allowdisplaybreaks

\numberwithin{equation}{section}

\newtheorem{theorem}{Theorem}[section]
\newtheorem{lemma}[theorem]{Lemma}

\newtheorem{corollary}[theorem]{Corollary}
\newtheorem{remark}[theorem]{Remark}

\newtheorem{Atheorem}{Theorem}

\newtheorem{Alemma}[Atheorem]{Lemma}
\newtheorem{Aremark}[Atheorem]{Remark}

\newtheorem{Btheorem}{Theorem}
\newtheorem{Bproposition}[Btheorem]{Proposition}

\newtheorem{Blemma}[Btheorem]{Lemma}
\newtheorem{Bremark}[Btheorem]{Remark}
\newtheorem{Bproblem}[Btheorem]{Problem}

\newtheorem{Ctheorem}{Theorem}

\newtheorem{Cproposition}[Ctheorem]{Proposition}
\newtheorem{Clemma}[Ctheorem]{Lemma}

\newtheorem{TheoA}{Theorem A1}
\newtheorem{TheoAA}{Theorem A2}
\newtheorem{TheoB}{Theorem B1}
\newtheorem{TheoBB}{Theorem B2}
\newtheorem{TheoC}{Theorem C}

\newcommand{\N}{\mathbb{N}}
\newcommand{\Z}{\mathbb{Z}}
\newcommand{\R}{\mathbb{R}}
\newcommand{\C}{\mathbb{C}}
\newcommand{\F}{\mathbb{F}}

\newcommand{\summ}{\sum\nolimits}

\def\G{\mathrm{G}}
\def\1{\mathbf{1}}
\def\H{\mathcal{H}}

\def\M{\mathcal{M}}
\def\A{\mathcal{A}}

\def\P{\mathcal{P}}

\def\S{\mathcal{S}}

\def\V{\mathrm{\mathcal{L}(G)}}

\newcommand{\dem}{\noindent {\bf Proof. }}
\newcommand{\demA}{\noindent {\bf Proof of Theorem A1. }}
\newcommand{\demAA}{\noindent {\bf Proof of Theorem A2. }}
\newcommand{\demB}{\noindent {\bf Proof of Theorem B1. }}
\newcommand{\demBB}{\noindent {\bf Proof of Theorem B2. }}
\newcommand{\demC}{\noindent {\bf Proof of Theorem C i). }}
\newcommand{\demCC}{\noindent {\bf Proof of Theorem C ii). }}

\newcommand{\fin}{\hspace*{\fill} $\square$ \vskip0.2cm}

\def\esssup{\mathop{\mathrm{ess \, sup}}}

\addtocounter{tocdepth}{0}

\begin{document}

\null

\vskip-30pt

\null

\title[Riesz transforms and Fourier multipliers]{Noncommutative Riesz transforms \\ ---dimension free bounds and Fourier multipliers---}

\author[Junge, Mei, Parcet]
{Marius Junge, Tao Mei, Javier Parcet}

\maketitle

\null

\vskip-30pt

\null

\begin{abstract}
We obtain dimension free estimates for noncommutative Riesz transforms associated to conditionally negative length functions on group von Neumann algebras. This includes Poisson semigroups, beyond Bakry's results in the commutative setting. Our proof is inspired by Pisier's method and a new Khintchine inequality for crossed products. New estimates include Riesz transforms associated to fractional laplacians in $\R^n$ (where Meyer's conjecture fails) or to the word length of free groups. Lust-Piquard's work for discrete laplacians on LCA groups is also generalized in several ways. In the context of Fourier multipliers, we will prove that H\"ormander-Mihlin multipliers are Littlewood-Paley averages of our Riesz transforms. This is highly surprising in the Euclidean and (most notably) noncommutative settings. As application we provide new Sobolev/Besov type smoothness conditions. The Sobolev type condition we give refines the classical one and yields dimension free constants. Our results hold for arbitrary unimodular groups.
\end{abstract}

\tableofcontents

\addtolength{\parskip}{+1ex}

\renewcommand{\theequation}{MP}
\addtocounter{equation}{-1}

\null

\vskip-60pt

\null

\section*{{\bf Introduction}}

Classical Riesz transforms $R_jf \! = \! \partial_j (-\Delta)^{-\frac12} \! f$ are higher dimensional forms of the Hilbert transform in $\R$. Dimension free estimates for the associated square functions $\mathcal{R}f = |\nabla (-\Delta)^{-\frac12}f|$ were first proved by Gundy/Varopoulos \cite{GV} and shortly after by Stein \cite{St3}, who pointed out the significance of a dimension free formulation of Euclidean harmonic analysis. The aim of this paper is to provide dimension free estimates for a much broader class of Riesz transforms and apply them for further insight in Fourier multiplier $L_p$--theory. Our approach is surprisingly simple and it is valid in the general context of group von Neumann algebras. 

A relevant generalization appeared in the work of P.A. Meyer \cite{M}, continued by Bakry, Gundy and Pisier \cite{B1,B2,G, PRiesz} among others. The probabilistic approach consists in replacing $-\Delta$ by the infinitesimal generator $A$ of a nice semigroup acting on a probability space $(\Omega,\mu)$. The gradient form $\langle \nabla f_1, \nabla f_2 \rangle$ is also replaced by the so-called \lq\lq carr\'e du champs\rq\rq${}$ $\Gamma(f_1,f_2) = \frac12 ( \overline{A(f_1)} f_2 + \overline{f_1} A(f_2) - A(\overline{f_1}f_2))$ and Meyer's problem for $(\Omega,\mu,A)$ consists in determining whether 
\begin{equation} \label{MeyerProblem}
\big\| \Gamma(f,f)^\frac12 \big\|_p \sim_{c(p)} \big\| A^\frac12 f \big\|_p \qquad (1 < p < \infty)
\end{equation} 
holds on a dense subspace of $\mathrm{dom}A$. Meyer proved this for the Ornstein-Uhlenbeck semigroup, while Bakry considered other diffusion semigroups assuming the $\Gamma^2 \ge 0$ condition. This condition yields in turn a lower bound for the Ricci curvature in Riemannian manifolds \cite{B2,Li}. Clifford algebras were considered by Lust-Piquard \cite{LuP1,LuP2} and other topics concerning Hermite or Laguerre processes and optimal linear estimates can be found in \cite{CD,DV,HTV,HdRST} and the references therein. Further dimension free estimates for maximal functions appear in \cite{B,C,NT,St2}.

\numberwithin{equation}{section}

Contrary to what might be expected, \eqref{MeyerProblem} \emph{fails} for the Poisson semigroup in $\R^n$ when $1 < p \le \frac{2n}{n+1}$ even allowing constants depending on $n$, see Appendix C for details. Bakry's argument heavily uses commutative diffusion properties and hence the failure of (MP) for subordinated processes and $p < 2$ does not contradict his work. Moreover, besides the heat semigroup, convolution processes have not been studied systematically. Lust-Piquard's theorem on discrete laplacians for LCA groups \cite{LuP} seems to be the only exception. Our first goal is to fill this gap and study Meyer's problem for Markov convolution semigroups in the Euclidean case and other group algebras. In this paper, we introduce a new form of \eqref{MeyerProblem}  which holds in much larger generality. As we shall see, this requires to follow here the tradition of noncommutative Khintchine inequalities \cite{LuP0,PX2} which imposes to consider an infimum over decompositions in two terms when $p<2$.  In the terminology from noncommutative geometry, our decomposition takes  place in the space of differential forms of order 1, see appendix Lemma B1. Indeed, the deeper understanding of derivations in noncommutative analysis provides a better understanding of Riesz transforms, even for semigroups of classical Fourier multipliers.

Let us first consider a simple model. Given a discrete abelian group $\G$, let $(\Omega,\mu)$ be the compact dual group with its normalized Haar measure. Construct the group characters $\chi_g: \Omega \to \mathbb{T}$ and consider $f \in L_p(\Omega)$. According to Schoenberg's theorem \cite{Sc} a convolution semigroup $\S_{\psi,t}: \chi_g \mapsto e^{-t \psi(g)} \chi_g$ is Markovian in $\Omega$ if and only if $\psi(e)=0$ for the identity $e$, $\psi(g) = \psi(g^{-1})$ for all $g \in \G$, and $\summ_g a_g = 0 \Rightarrow \summ_{g,h} \overline{a_g} a_h \psi(g^{-1}h) \le 0$. Any such function $\psi$ is called a conditionally negative length. $A_\psi(\chi_g) = \psi(g) \chi_g$ is the generator, which determines the gradient form $\Gamma_\psi$. Does \eqref{MeyerProblem} or a generalization of it holds for arbitrary pairs $(\G,\psi)$? To answer this question we first widen the scope of the problem and consider its formulation for non-abelian discrete groups $\G$. The former r\^ole of $L_\infty(\Omega,\mu)$ for $\G$ abelian is now played by the group von Neumann algebra $\V$, widely studied in noncommutative geometry and operator algebras \cite{BO,Con,Dav}. 

Let $\G$ be a discrete group with left regular representation $\lambda: \G \to \mathcal{B}(\ell_2(\G))$ given by $\lambda(g) \delta_h = \delta_{gh}$, where the $\delta_g$'s form the unit vector basis of $\ell_2(\G)$. Write $\mathcal{L}(\G)$ for its group von Neumann algebra, the weak operator closure of the linear span of $\lambda(\G)$ in $\mathcal{B}(\ell_2(\G))$. Consider the standard trace $\tau(\lambda(g)) = \delta_{g=e}$ where $e$ denotes the identity of $\G$. Any $f \in \mathcal{L}(\G)$ has a Fourier series expansion of the form $\sum_g \widehat{f}(g) \lambda(g)$ with $\tau(f) = \widehat{f}(e)$. If we consider $(\V,\tau)$ as a noncommutative measure space, we may define the $L_p$ space associated to it as $$L_p(\widehat{\mathbf{G}}) = L_p(\mathcal{L}(\mathrm{G}), \tau) \equiv \mbox{Closure of $\mathcal{L}(\G)$ wrt $\|f\|_{L_p(\widehat{\mathbf{G}})} = \big( \tau |f|^p \big)^\frac1p$}.$$ It turns out that $L_p(\widehat{\mathbf{G}}) = L_p(\widehat{\G}) = L_p(\Omega,\mu)$ for abelian $\G$, after identifying $\lambda(g)$ with the character $\chi_g$. In general, the (unbounded) operator $|f|^p$ is obtained from functional calculus on the Hilbert space $\ell_2(\G)$, see \cite{PX2} or Appendix A for further details. Harmonic analysis on $\V$ places the group on the frequency side. This approach is inspired by the remarkable results of Cowling/Haagerup \cite{CH,H} on the approximation property and Fourier multipliers on group algebras. This paper is part of an effort \cite{JMP0,JPPP,PRo} to extend modern harmonic analysis to the unexplored context of group von Neumann algebras. Markovian semigroups acting on $\V$ are composed of self-adjoint, completely positive and unital maps. Schoenberg's theorem is still valid and \vskip-10pt $$\S_{\psi,t}f \, = \, \sum_{g \in \G}^{\null} e^{-t\psi(g)} \widehat{f}(g) \lambda(g)$$ will be Markovian if and only if $\psi: \G \to \R_+$ is a conditionally negative length.

Riesz transforms should look like $R_{\psi,j}f = \partial_{\psi,j} A_\psi^{-\frac12}f$ where the former laplacian is now replaced by $A_\psi(\lambda(g)) = \psi(g) \lambda(g)$ and $\partial_{\psi,j}$ is certain differential operator playing the role of a directional derivative. Unlike for $\R^n$, there is no standard differential structure for an arbitrary discrete $\G$. The additional structure comes from the length $\psi$, which allows a broader interpretation of tangent space in terms of the associated cocycle. Namely, conditionally negative lengths are in one-to-one correspondence with affine representations $(\H_\psi, \alpha_\psi, b_\psi)$, where $\alpha_\psi: \G \to O(\H_\psi)$ is an orthogonal representation over a real Hilbert space $\H_\psi$ and $b_\psi: \G \to \H_\psi$ is a mapping satisfying the cocycle law \vskip-10pt $$b_\psi(gh) = \alpha_{\psi,g}(b_\psi(h)) + b_\psi(g).$$ 

\noindent Since $\partial_j (\exp(2 \pi i \langle x, \cdot \rangle)) = 2 \pi i x_j \exp(2 \pi i \langle x, \cdot \rangle)$, it is natural to define \vskip-5pt $$R_{\psi,j} f = \partial_{\psi,j} A_\psi^{- \frac12} f = 2 \pi i \sum_{g \in \G} \frac{\langle b_\psi(g), e_j \rangle_{\H_\psi}^{\null}}{\sqrt{\psi(g)}} \widehat{f}(g) \lambda(g)$$ for some orthonormal basis $(e_j)_{j \ge 1}$ of $\H_\psi$. Recalling that $b_\psi(g) / \sqrt{\psi(g)}$ is always a normalized vector, we recover the usual symbol of $R_j$ as a Fourier multiplier. Note also that classical Riesz transforms can be seen from this viewpoint. Namely, de Leeuw's theorem \cite{dL} allows us to replace $\R^n$ by its Bohr compactification, whose $L_p$ spaces come from the group von Neumann algebra $\mathcal{L}(\R^n_{\mathrm{disc}})$ of $\R^n$ equipped with the discrete topology. Then, classical Riesz transforms arise from the standard cocycle where $\psi(\xi) = |\xi|^2$ (generating the heat semigroup) and $\H_\psi = \R^n$ with the trivial action $\alpha_\psi$ and the identity map $b_\psi$ on $\R^n$. Moreover, the classical Riesz transforms vanish on the $L_p$-functions fixed by the heat semigroup: the constant functions in the $n$-torus and the zero function in the Euclidean space. This is also the case here and $R_{\psi,j}$ will be properly defined on \vskip-5pt 
$$L_p^\circ(\widehat{\mathbf{G}}) = \Big\{ f \in L_p(\widehat{\mathbf{G}}) \, \big| \, \widehat{f}(g)=0 \ \mbox{whenever} \ b_\psi(g)=0 \Big\}.$$ 

\noindent An elementary calculation shows that $$\Gamma_\psi \big( A_\psi^{-\frac12} f, A_\psi^{-\frac12} f \big) \, = \, \sum_{j \ge 1} \big| R_{\psi,j} f \big|^2.$$ By Khintchine inequality, \eqref{MeyerProblem} in the commutative setting is equivalent to $$\Big\|  \sum_{j \ge 1} \gamma_j R_{\psi,j}f \Big\|_{L_p(\Omega \times \widehat{\G})} \, \sim_{c(p)} \, \|f\|_{L_p(\widehat{\G})}$$ for any family $(\gamma_j)_{j \ge 1}$ of centered independent gaussians. We have already pointed that this fails when $A = \Delta^\frac12$ is the generator of the Poisson semigroup. According to the standard gaussian measure space construction (see below), we may construct a canonical action $\beta_\psi: \G \curvearrowright L_\infty(\Omega)$ determined by $\psi$. As we shall justify in this paper, a natural revision of Meyer's problem \eqref{MeyerProblem} is to ask whether $$\Big\| \sum_{j \ge 1} \gamma_j \rtimes R_{\psi,j}f \Big\|_{L_p(L_\infty(\Omega) \times \widehat{\G})} \, \sim_{c(p)} \, \|f\|_{L_p(\widehat{\G})}$$ Our first result is that this form of \eqref{MeyerProblem} holds for all Markov semigroups on group von Neumann algebras, including the classical Poisson semigroups. The Riesz transforms above were introduced in \cite{JMP0} under the additional assumption that $\dim \H_\psi < \infty$. Our dimension free estimates below ---in cocycle form, see Theorem A2 for a Meyer's type formulation--- allow to consider Riesz transforms associated to infinite-dimensional cocycles. 

\begin{TheoA} 
Let $f \in L_p^\circ(\widehat{\mathbf{G}})$ and $1 < p < \infty \! :$
\begin{itemize}
\item[i)] If $1 < p \le 2$, we get $$\|f\|_p \, \sim_{c(p)} \, \inf_{R_{\psi,j}f=a_j+b_j} \Big\| \Big( \sum_{j \ge 1} a_j^*a_j \Big)^{\frac12} \Big\|_p + \Big\| \Big( \sum_{j \ge 1} \widetilde{b}_j \widetilde{b}_j^* \Big)^{\frac12} \Big\|_p.$$ 

\vskip3pt

\item[ii)] If $2 \le p < \infty$, we get $$\hskip18pt \|f\|_p \sim_{c(p)} \max \left\{ \Big\| \Big( \sum_{j \ge 1} (R_{\psi,j}f)^*(R_{\psi,j}f) \Big)^{\frac12} \Big\|_p, \Big\| \Big( \sum_{j \ge 1} (R_{\psi,j}f^*)^*(R_{\psi,j}f^*) \Big)^{\frac12} \Big\|_p \right\}.$$
\end{itemize}
\end{TheoA}

The infimum in Theorem A1 i) runs over all decompositions $R_{\psi,j} f = a_j + b_j$ in the tangent module, the noncommutative analogue of the module of differential forms of order one. A more precise description will be possible after the statement of Theorem A2. A crucial aspect comes from the $\widetilde{b}_j$'s, twisted forms of $b_j$'s that will be rigorously defined in the paper. The failure of Theorem A1 for $b_j$'s instead of their twisted forms was discovered by Lust-Piquard \cite[Proposition 2.9]{LuP} in her work with discrete laplacians on LCA groups. It shows certain \lq intrinsic noncommutativity\rq${}$ in the problem, since the statement for $p < 2$ does not simplify for $\G$ abelian unless the action $\alpha_\psi$ is trivial! We refer to Remark \ref{HilbertModule} or Appendix B for more on this and to Remark \ref{RemRowColumnRiesz} for a more transparent row/column behavior when $p \ge 2$.  

A great variety of new and known estimates for Riesz transforms and other Fourier multipliers arise from Theorem A1, by considering different lengths. All conditionally negative length functions appear as deformations of the canonical inner cocycle for the left regular representation. Namely, if we consider the space $\Pi_0$ of trigonometric polynomials in $\V$ whose Fourier coefficients have vanishing sum ---finite sums $\sum_g a_g \lambda(g)$ with $\sum_g a_g = 0$--- then $\psi: \G \to \R_+$ is conditionally negative iff $\psi(g) = \tau_\psi ( 2 \lambda(e) - \lambda(g) - \lambda(g^{-1}) )$ for some positive linear functional $\tau_\psi: \Pi_0 \to \C$. This characterization will be useful along the paper and we will prove it in Appendix A. 
Once we have identified the exact form of conditionally negative lengths, let us now illustrate Theorem A1 with a few examples which will be analyzed in the body of the paper:
\begin{itemize}
\item[a)] \textbf{Fractional laplacians in $\R^n$.} Our results also hold for group algebras over arbitrary unimodular groups. In the particular case $\G = \R^n$ we may consider conditionally negative lengths of the form $$\hskip32pt \psi (\xi) \, = \, 2 \int_{\R^n} \big( 1 - \cos (2 \pi \langle x, \xi \rangle ) \big) \, d\mu_\psi(x)$$ for a positive Borel measure $\mu_\psi$ satisfying $\psi(\xi) < \infty$ for all $\xi \in \R^n$. If $d\mu_\psi(x) \sim dx/|x|^{n+2\beta}$ for any $0 < \beta < 1$, we get $\psi(\xi) = |\xi|^{2\beta}$. This will  provide us dimension free estimates for Riesz transforms associated to fractional laplacians, which are new. The estimates predicted by Meyer fail, see Appendix C. Contrary to the case $\beta = 1$, we find highly nontrivial cocycles. The vast family of measures $\mu_\psi$ are explored in further generality in the second part of this paper.

\vskip3pt

\item[b)] \textbf{Discrete laplacians in LCA groups.} Let $\Gamma_0$ be a locally compact abelian group and $s_0 \in \Gamma_0$ be torsion free. If $\partial_j f(\gamma) = f(\gamma) - f (\gamma_1, \ldots, s_0 \gamma_j, \ldots, \gamma_n)$ stand for discrete directional derivatives in $\Gamma = \Gamma_0 \times \Gamma_0 \times \cdots \times \Gamma_0$, we may consider the laplacian $\mathcal{L} = \sum_j \partial_j^* \partial_j$ and $R_j = \partial_j \mathcal{L}^{-1/2}$. Lust-Piquard provided dimension free estimates for these Riesz transforms in her paper \cite{LuP}. If we set $\sigma_j = (0, \ldots, 0, s_0, 0, \ldots, 0)$ with $s_0$ in the $j$-th entry, consider the sum of point-masses $\mu_\psi = \sum_j \delta_{\sigma_j}$. Then we shall recover Lust-Piquard's theorem via Theorem A1 taking $$\hskip34pt \psi(g) = \int_{\Gamma} \big( 2 - \chi_g - \chi_{g^{-1}} \big) (\gamma) \, d\mu_\psi(\gamma) \quad \mbox{for} \quad g \in \G = \widehat{\Gamma}.$$ The advantage is that we do not need to impose $s_0$ to be torsion free and may take $\Gamma_0 = \Z_k$ for any $k \ge 2$, which was left open in \cite{LuP}. Moreover, our formulation holds for any finite sum of point-masses, so that we may allow the shift $s_0$ to depend on the entry $j$ or even the group $\Gamma$ not to be given in a direct product form... This solves the problem of discrete laplacians in a very general form, continuous analogues can also be given.

\vskip3pt

\item[c)] \textbf{Word-length laplacians.} Consider a finitely generated group $\G$ and write $|g|$ to denote the word length of $g$, its distance to $e$ in the Cayley graph. If it is conditionally negative ---like for free, cyclic, Coxeter groups--- a natural laplacian is $A_{|\ |}(\lambda(g)) = |g| \lambda(g)$ and the Riesz transforms $$\hskip20pt R_{| \ |, j}f \, = \, \partial_{| \ |,j} A_{| \ |}^{-\frac12}f \, = \, 2 \pi i \summ_g \frac{\langle b_{| \ |}(g), e_j \rangle_{\H_{| \ |}}}{\sqrt{|g|}} \widehat{f}(g) \lambda(g)$$ satisfy \hskip-1pt Theorem A1. \hskip-2pt Many other \hskip-1pt transforms arise from \hskip-1pt other \hskip-1pt conditionally negative lengths. The natural example given above is out of the scope of the method in \cite{JMP0}. It yields new interesting inequalities, here are two examples in the (simpler) case $p \ge 2$. When $\G = \Z_{2m}$ $$\hskip36pt \Big\| \sum_{j \in \Z_{2m}} \widehat{f}(j) e^{2\pi i \frac{j}{2m} \cdot} \Big\|_p \sim_{c(p)} \Big\| \Big( \sum_{k \in \Z_{2m}} \Big| \sum_{j \in \Lambda_k} \frac{\widehat{f}(j)}{\sqrt{j \wedge (2m-j)}} e^{2 \pi i \frac{j}{2m} \cdot} \Big|^2 \Big)^\frac12 \Big\|_p$$ for $\Lambda_k = \big\{ j \in \Z_{2m} : j-k \equiv s \ (2m) \ \mbox{with} \ 0 \le s \le m-1 \big\}$. When $\G = \mathbb{F}_n$ $$\hskip32pt \|f\|_p \, \sim_{c(p)} \, \Big\| \Big( \sum_{h \neq e} \Big| \sum_{g \ge h} \frac{1}{\sqrt{|g|}} \widehat{f}(g) \lambda(g) \Big|^2 + \Big| \sum_{g \ge h} \frac{1}{\sqrt{|g|}} \overline{\widehat{f}(g^{-1})} \lambda(g) \Big|^2 \Big)^\frac12 \Big\|_p.$$
\end{itemize}
 
Let us now go back to Meyer's problem \eqref{MeyerProblem} for convolution Markov semigroups. Integration by parts gives $- \Delta = \nabla^* \circ \nabla$. According to Sauvageot's theorem \cite{Sa}, we know that a similar factorization takes place for Markovian semigroups. Namely,  there exists a Hilbert $\V$-bimodule $\mathsf{M}_\psi$ and a densely defined closable symmetric derivation $$\delta_\psi: L_2(\widehat{\mathbf{G}}) \to \mathsf{M}_\psi \quad \mbox{such that} \quad A_\psi = \delta_\psi^* \delta_\psi^{\null}.$$ If $B: e_j \in \H_\psi \mapsto \gamma_j \in L_2(\Omega, \Sigma, \mu)$ denotes the standard gaussian measure space construction, we will find in our case that $\mathsf{M}_\psi = L_\infty(\Omega, \Sigma,\mu) \rtimes_{\beta_{\psi}}Ê\G$ where $\G$ acts via the cocycle action as $\beta_{\psi,g}(B(h))= B(\alpha_{\psi,g}(h))$. The derivation is $$\delta_\psi: \lambda(g) \mapsto B(b_\psi(g)) \rtimes \lambda(g),$$ $$\delta_\psi^*: \rho \rtimes \lambda(g) \mapsto \big\langle \rho, B(b_\psi(g)) \big\rangle \, \lambda(g).$$ If we consider the conditional expectation onto $\V$ $$\mathsf{E}_\V: \summ_g \rho_g \rtimes \lambda(g) \in \mathsf{M}_\psi \mapsto \summ_g \int_\Omega \rho_g \hskip1pt d\mu \, \lambda(g) \in \V,$$ we may obtain the following solution to Meyer's problem for the pair $(\G,\psi)$.

\begin{TheoAA}
The following norm equivalences hold$\, :$
\begin{itemize}
\item[i)] If $1 < p \le 2$, we have $$\big\| A_\psi^\frac12 f \big\|_p \, \sim_{c(p)} \, \inf_{\delta_\psi f = \phi_1 + \phi_2} \Big\| \mathsf{E}_\V \big( \phi_1^* \phi_1 \big)^{\frac12} \Big\|_p + \Big\| \mathsf{E}_\V \big( \phi_2 \phi_2^* \big)^{\frac12} \Big\|_p.$$

\item[ii)] If $2 \le p < \infty$, we have $$\hskip12pt \big\| A_\psi^\frac12 f \big\|_p \, \sim_{c(p)} \, \max \Big\{ \Big\| \mathsf{E}_\V \big( \delta_\psi f^* \delta_\psi f \big)^{\frac12} \Big\|_p, \Big\| \mathsf{E}_\V \big( \delta_\psi f \delta_\psi f^* \big)^{\frac12} \Big\|_p\Big\}.$$
\end{itemize}
\end{TheoAA}

We are now ready to describe the families of operators along which we allow our decompositions in Theorems A1 i) and A2 i) to run over. Let us consider the subspace $G_p(\C) \rtimes \G$ of $L_p(L_\infty(\Omega) \rtimes \G)$ formed by operators of the form $$\phi = \sum_{g \in \G} \underbrace{\sum_{j \ge 1} \phi_{g,j} B(e_j)}_{\phi_g} \rtimes \lambda(g)$$ with $\phi_{g,j} \in \C$ and $\phi_g \in L_p(\Omega)$. The infimum in Theorem A2 is taken over all possible decompositions $\delta_\psi f = \phi_1 + \phi_2$ with $\phi_1, \phi_2 \in G_p(\C) \rtimes \G$. On the other hand, to describe the infimum in Theorem A1 we introduce the maps $$u_j: G_p(\C) \rtimes \G \ni \phi \mapsto \sum_{g \in \G} \big\langle \phi_g, B(e_j) \big\rangle_{L_2(\Omega)} \lambda(g) \in L_p(\widehat{\mathbf{G}}).$$ Then $R_{\psi,j}f = a_j + b_j$ runs over $(a_j, b_j) = (u_j(\phi_a), u_j(\phi_b))$ for $\phi_a, \phi_b \in G_p(\C) \rtimes \G$.

The relation with Meyer's formulation comes from $$\Gamma_\psi(f_1,f_2) \, = \, \mathsf{E}_\V \big( \delta_\psi f_1^* \delta_\psi f_2 \big).$$ As in Theorem A1, we recover Meyer's inequalities \eqref{MeyerProblem} when $\G$ is abelian and the cocycle action is trivial, the general case is necessarily more involved. The infimum can not be reduced to decompositions $f=f_1+f_2$, see Remark \ref{NoEasyInf}. The main result in \cite{JM} provides lower estimates for $p \ge 2$ and regular Markov semigroups satisfying $\Gamma^2 \ge 0$. In the context of group algebras, Theorem A2 goes much further. 

\renewcommand{\theequation}{RI}
\addtocounter{equation}{-1} 

Theorems A1 and A2 are valid for arbitrary unimodular groups and generalize to matrix-valued functions. We refer to Remark \ref{BestConstant} for a brief analysis on optimal constants. Theorem A1 follows by standard manipulations from Theorem A2. The proof of the latter is inspired by a crossed product extension of Pisier's method \cite{PRiesz} which ultimately relies on a Khintchine type inequality of independent interest. The key point in Pisier's argument is to identify the Riesz transform as a combination of the transferred Hilbert transform $$Hf(x,y) \, = \, \mbox{p.v.} \, \frac{1}{\pi} \int_\R \beta_t f(x,y) \frac{dt}{t} \quad \mbox{where} \quad \beta_tf(x,y) = f(x+ty)$$ and the gaussian projection $Q: L_p(\R^n, \gamma) \to \overline{L_p-\mbox{span}} \big\{ B(\xi) \, | \, \xi \in \R^n \big\}$. Here the gaussian variables are given by $B(\xi)(y) = \langle \xi,y \rangle$, homogenous polynomials of degree $1$. The following identity can be found in \cite{PRiesz} for smooth $f$ 
\begin{equation} \label{PisierFormula}
\sqrt{\frac{2}{\pi}} \, \delta (-\Delta)^{- \frac12} f = (id_{\R^n} \otimes Q) \Big( \mbox{p.v.} \, \frac{1}{\pi} \int_\R \beta_t f \frac{dt}{t} \Big),\end{equation} 
where $\delta: C^\infty(\R^n) \to C^\infty(\R^n \times \R^n)$ is the following derivation $$\delta(f)(x,y) = \sum_{k=1}^n \frac{\partial f}{\partial x_k} y_k = \big\langle \nabla f(x), y \big\rangle.$$ Our Khintchine inequality allows to generalize this formula for pairs $(\G,\psi)$. It seems fair to say that for the Euclidean case, this kind of formula has it roots in the work of Duoandikoetxea and Rubio de Francia \cite{DR} through the use of Calder\'on's method of rotations. Pisier's main motivation was to establish similar identities including Riesz transforms for the Ornstein-Ulhenbeck semigroup.

\numberwithin{equation}{section}

Our class of $\psi$-Riesz transforms becomes very large when we move $\psi$. This yields a fresh perspective in Fourier multiplier theory, mainly around H\"ormander-Mihlin smoothness conditions in terms of Sobolev and (limiting) Besov norms. We refer to \cite{JMP0} for a more in depth discussion on smoothness conditions for Fourier multipliers defined on discrete groups. The main idea is that this smoothness may and should be measured through the use of cocyles, via lifting multipliers $\widetilde{m}$ living in the cocycle Hilbert space, so that $m = \widetilde{m} \circ b_\psi$. Let $\mathsf{M}_p(\widehat{\mathbf{G}})$ be the space of multipliers $m: \G \to \C$ equipped with the $p \to p$ norm of the map $\lambda(g) \mapsto m(g) \lambda(g)$. Let us consider the classical differential operators in $\R^n$ $$\widehat{\mathsf{D}_\alpha f} (\xi) = |\xi|^\alpha \widehat{f}(\xi) \quad \mbox{for} \quad \xi \in \R^n \!\! = \! \H_\psi.$$ We shall also need the fractional laplacian lengths $$\psi_\beta(\xi) = 2 \int_{\R^n} \big( 1 - \cos (2 \pi \langle \xi, x \rangle) \big) \, \frac{dx}{|x|^{n+2\beta}} = \mathrm{k}_n(\beta) |\xi|^{2\beta}.$$ Our next result provides new Sobolev conditions for the lifting multiplier.

\begin{TheoB} 
Let $(\G,\psi)$ be a discrete group equipped with a conditionally negative length giving rise to an $n$-dimensional cocycle $(\H_\psi, \alpha_\psi, b_\psi)$. Let $(\varphi_j)_{j \in \Z}$ denote a standard radial Littlewood-Paley partition of unity in $\R^n$. Then, if $1 < p < \infty$ and $\varepsilon > 0$, the following estimate holds $$\|m\|_{\mathsf{M}_p(\widehat{\mathbf{G}})} \, \lesssim_{c(p,n)} \, |m(e)| + \inf_{m = {\widetilde{m}}^{\null} \circ b_\psi} \left\{ \sup_{j \in \Z} \Big\| \mathsf{D}_{\frac{n}{2} + \varepsilon} \Big( \sqrt{\psi_\varepsilon} \, \varphi_j \, \widetilde{m} \Big) \Big\|_{L_2(\R^n)} \right\}.$$ The infimum runs over all lifting multipliers $\widetilde{m}: \H_\psi \to \C$ such that $m = \widetilde{m} \circ b_\psi$.
\end{TheoB}

Our Sobolev type condition in Theorem B1 is formally less demanding than the standard one (see below) and our argument is also completely different from the classical approach used in \cite{JMP0}. As a crucial novelty, we will show that every H\"ormander-Mihlin type multiplier (those for which the term on the right hand side is finite, in particular the classical ones) is in fact a Littlewood-Paley average of Riesz transforms associated to a single infinite-dimensional cocycle! The magic formula comes from an isometric isomorphism between the Sobolev type norm in Theorem B1 and mean-zero elements of $L_2(\R^n, \mu_\varepsilon)$ with $d\mu_\varepsilon(x) = |x|^{-(n+2\varepsilon)} dx$. In other words, $m: \mathbb{R}^n \to \mathbb{C}$ satisfies $$\Big\| \mathsf{D}_{\frac{n}{2} + \varepsilon} \big( \sqrt{\psi_\varepsilon} m \big) \Big\|_{L_2(\mathbb{R}^n)} < \infty$$ iff there exists a mean-zero $h \in L_2(\mathbb{R}^n, \mu_\varepsilon)$ such that $$m(\xi) \, = \, \frac{\big\langle h, b_\varepsilon(\xi) \big\rangle_{\mu_\varepsilon}}{\sqrt{\psi_\varepsilon(\xi)}} \quad \mbox{and} \quad \|h\|_{L_2(\R^n, \mu_\varepsilon)} = \Big\| \mathsf{D}_{\frac{n}{2} + \varepsilon} \big( \sqrt{\psi_\varepsilon} m \big) \Big\|_{L_2(\mathbb{R}^n)}.$$ A few remarks are in order:
\begin{itemize}
\item Theorem B1 also holds for any unimodular group.

\vskip3pt

\item Our condition is bounded above by the classical one $$\hskip20pt \sup_{j \in \Z} \Big\| \big( 1 + | \ |^2 \big)^{\frac{n}{4} + \varepsilon} \big( \varphi_0 \, \widetilde{m}(2^j \cdot) \big)^{\wedge} \Big\|_{L_2(\R^n)}.$$ A crucial fact is that our Sobolev norm is dilation invariant, more details in Corollary \ref{SobVsClassical}. Moreover, our condition is more appropriate in terms of dimensional behavior of the constants, see Remark \ref{ConstantsSobolev}.

\vskip3pt

\item Our result is stronger than the main result in \cite{JMP0} in two respects. First we obtain Sobolev type conditions, which are way more flexible than the Mihlin assumptions $$\sup_{\xi \in \R^n \setminus \{0\}} \sup_{|\beta| \le [\frac{n}{2}] + 1} |\xi|^{|\beta|} \big|\partial_\beta \widetilde{m}(\xi) \big| < \infty.$$ Second, we avoid the modularity restriction in \cite{JMP0}. Namely, there we needed a simultaneous control of left/right cocycles for non-abelian discrete groups, when there is no spectral gap. In the lack of that, we could also work only with the left cocycle at the price of extra decay in the smoothness condition. In Theorem B1, it suffices to satisfy our Sobolev-type conditions for the left cocycle. On the other hand the approach in \cite{JMP0} is still necessary. First, it explains the connection between H\"ormander-Mihlin multipliers and Calder\'on-Zygmund theory for group von Neumann algebras. Second, the Littlewood-Paley estimates obtained in \cite{JMP0} are crucial for this paper and \cite{PRo}. Third, our approach here does not give $L_\infty \to \mathrm{BMO}$ estimates.  
\end{itemize}

The dimension dependence in the constants of Theorem B1 has its roots in the use of certain Littlewood-Paley inequalities on $\V$, but not on the Sobolev type norm itself. This yields a form of the H\"ormander-Mihlin condition with dimension free constants, replacing the compactly supported smooth functions $\varphi_j$ by certain class $\mathcal{J}$ of analytic functions which arises from Cowling/McIntosh holomorphic functional calculus, the simplest of which is $x \mapsto x e^{-x}$ that already appears in the work of Stein \cite{St}. Our result is the following.

\begin{TheoBB} 
Let $\G$ be a discrete group and $\Lambda_\G$ the set of conditionally negative lengths $\psi: \G \to \R_+$ giving rise to a finite-dimensional cocycle $(\H_\psi, \alpha_\psi, b_\psi)$. Let $\varphi: \R_+ \to \C$ be an analytic function in the class $\mathcal{J}$. Then, if $1 < p < \infty$ and $\varepsilon > 0$ the following estimate holds $$\|m\|_{\mathsf{M}_p(\widehat{\mathbf{G}})} \, \lesssim_{c(p)} \, |m(e)| + \inf_{\begin{subarray}{c} \psi \in \Lambda_\G \\ m = {\widetilde{m}} \circ b_\psi \end{subarray}} \left\{ \esssup_{s > 0} \Big\| \mathsf{D}_{\frac{\dim \H_\psi}{2} + \varepsilon} \Big( \sqrt{\psi_\varepsilon} \, \varphi (s \, | \cdot |^2) \, \widetilde{m} \Big) \Big\|_2 \right\}.$$ The infimum runs over all $\psi \in \Lambda_\G$ and all  $\widetilde{m}: \H_\psi \to \C$ such that $m = \widetilde{m} \circ b_\psi$.
\end{TheoBB}

Taking the trivial cocycle in $\R^n$ whose associated length function is $|\xi|^2$, we find a Sobolev condition which works up to dimension free constants, we do not know whether this statement is known in the Euclidean setting. The versatility of Theorems B1 and B2 for general groups is an illustration of what can be done using other conditionally negative lengths to start with. Replacing for instance the fractional laplacian lengths by some others associated to limiting measures when $\varepsilon \to 0$, we may improve the Besov type conditions \`a la Baernstein/Sawyer \cite{BS}, see also the related work of Seeger \cite{Se3,Se2} and  \cite{C2,Lz,Se1}. The main idea is to replace the former measures $$d\mu_\varepsilon(x) = \frac{dx}{|x|^{n+2\varepsilon}}$$ used to prove Theorem B1, by the limiting measure $d\nu(x) = u(x) dx$ with $$u(x) = \frac{1}{|x|^n} \Big( \chi_{\mathrm{B}_1(0)}(x) + \frac{1}{1 + \log^2 |x|} \chi_{\R^n \setminus \mathrm{B}_1(0)}(x) \Big).$$ Let us also consider the associated length $$\gamma(\xi) \, = \, 2 \int_{\R^n} \big( 1 - \cos(2 \pi \langle \xi, x \rangle) \big) \, u(x) \, dx.$$ Then, if $1 < p < \infty$ and $\dim \H_\psi=n$, we prove in Theorem \ref{BesovThm} that
\begin{eqnarray*}
\|m\|_{\mathsf{M}_p(\widehat{\mathbf{G}})} & \lesssim_{c(p,n)} & |m(e)| + \inf_{m = {\widetilde{m}}^{\null} \circ b_\psi} \left\{ \sup_{j \in \Z} \Big( \sum_{k \in \Z} 2^{nk} \mathrm{w}_k \big\| \widehat{\varphi}_k * ( \sqrt{\gamma} \, \varphi_j \, \widetilde{m} ) \big\|_2^2 \Big)^{\frac12} \right\},
\end{eqnarray*}
where $(\varphi_j)_{j \in \Z}$ is a standard radial Littlewood-Paley partition of unity in $\R^n$ and the weights $\mathrm{w}_k$ are of the form $\delta_{k \le 0} + k^2 \delta_{k > 0}$ for $k \in \Z$. A more detailed analysis of this result will be given in Paragraph \ref{Besov}. The huge variety of infinite-dimensional cocycles will be further explored within the context of Euclidean harmonic analysis in a forthcoming paper. In fact, an even more general construction is possible which relates Riesz transforms with \lq\lq Sobolev type norms\rq\rq${}$ directly constructed in group von Neumann algebras. We will not explore this direction here, further details in Remarks \ref{Abstract1} and \ref{Abstract2}.

Let us now consider a given branch in the Cayley graph of $\mathbb{F}_\infty$, the free group with infinitely many generators. Of particular interest are two applications we have found for operators (frequency) supported by such a branch. If we fix a branch $\mathrm{B}$ of $\mathbb{F}_\infty$ let us set $$L_p(\widehat{\mathbf{B}}) = \Big\{ f \in L_p(\mathcal{L}(\mathbb{F}_\infty)) \ \big| \ \widehat{f}(g) = 0 \ \mbox{for all} \ g \notin \mathrm{B} \Big\}.$$ As usual, we shall write $| \ |$ to denote the word length of the free group $\mathbb{F}_\infty$.

\begin{TheoC}
Given any branch $\mathrm{B}$ of $\mathbb{F}_\infty \! :$
\begin{itemize}
\item[i)] \textsf{\emph{H\"ormander-Mihlin multipliers.}} If $m: \mathbb{Z}_+ \to \C$ \vskip-8pt $$\hskip20pt \big\| \widetilde{m} \circ | \ | \ \big\|_{\mathsf{M}_p(\widehat{\mathbf{B}})} \, \lesssim_{c(p)} \, \sup_{j \ge 1} \, |\widetilde{m}(j)| + j |\widetilde{m}(j) - \widetilde{m}(j-1)|,$$ where $\mathsf{M}_p(\widehat{\mathbf{B}})$ denotes the space of $L_p(\widehat{\mathbf{B}})$-bounded Fourier multipliers. 

\vskip3pt

\item[ii)] \textsf{\emph{Twisted Littlewood-Paley estimates.}} Consider a standard Littlewood-Paley partition of unity $(\varphi_j)_{j \ge 1}$ in $\R_+$, generated by dilations of a function $\phi$ with $\sqrt{\phi}$ Lipschitz. Let $\Lambda_j: \lambda(g) \mapsto \sqrt{\varphi_j(|g|)} \lambda(g)$ denote the corresponding radial multipliers in $\mathcal{L}(\mathbb{F}_\infty)^{\null}$. Then, the following estimates hold for any $f \in L_p(\widehat{\mathbf{B}})$ and $1 < p < 2$ $$\hskip40pt \inf_{\Lambda_jf=a_j+b_j} \Big\| \Big( \sum_{j \ge 1} a_j^* a_j + \widetilde{b}_j \widetilde{b}_j^* \Big)^\frac12 \Big\|_{L_p(\widehat{\mathbf{B}})} \, \lesssim_{c(p)} \, \|f\|_{L_p(\widehat{\mathbf{B}})},$$ $$\hskip40pt \|f\|_{L_p(\widehat{\mathbf{B}})} \, \lesssim_{c(p)} \, \inf_{\Lambda_jf=a_j+b_j} \Big\| \Big( \sum_{j \ge 1} a_j^* a_j + b_j b_j^* \Big)^\frac12 \Big\|_{L_p(\widehat{\mathbf{B}})}.$$
\end{itemize}
\end{TheoC}

In analogy with Theorem A1, the first infimum runs over all decompositions with $(a_j, b_j) = (v_j(\phi_a), v_j(\phi_b))$ with $\phi_a, \phi_b \in G_p(\C) \rtimes \G$ and $v_j(\phi) = \sum_g \langle \phi_g, B(h_j) \rangle \lambda(g)$ for certain $h_j \in \H_{| \, |}$ to be defined in the paper. The second infimum runs over $(a_j, b_j) \in C_p(L_p) \times R_p(L_p)$, the largest space where it is meaningful. Theorem C shows that H\"ormander-Mihlin multipliers on branches of $\mathbb{F}_\infty$ behave like in the $1$-dimensional groups $\Z$ or $\R$. However, general branches have no group structure and $L_p$-norms admit less elementary combinatorics ($p \in 2\Z$) than the trivial ones $g^{\null}, g^2, g^3, \ldots$ with $g$ a generator. The key idea for our Littlewood-Paley inequalities is to realize the Littlewood-Paley partition of unity as a family of Riesz transforms. The crucial difference with Theorems A1 and A2 is that this family does not arise from an orthonormal basis, but from a quasi-orthonormal incomplete system. It is hence very likely that norm equivalences do not hold for nontrivial branches. On the contrary, our result shows that the untwisted square function is greater than the twisted one, and both coincide when the cocycle action is trivial and the product commutes. This is the case for trivial branches (associated to subgroups isomorphic to $\Z$) since we may replace the word length by the one coming from the heat semigroup on $\mathbb{T}$, which yields a trivial cocycle action. We may also obtain lower estimates for $p > 2$, see Corollary \ref{UpperBranch}. At the time of this writing, we do not know an appropriate upper bound for $\|f\|_p$ ($p>2$) since standard duality fails due to the twisted nature of square functions. Bo\.zejko-Fendler theorem \cite{BF} indicates that sharp truncations might not work for all values of $1 < p< \infty$.

Our approach requires some background on noncommutative $L_p$ spaces, group von Neumann algebras and geometric group theory. A brief survey of the main notions and results used in this paper is included in Appendix A for the non-expert reader. Appendix B contains a geometric analysis of our results in terms of the tangent module naturally associated to the infinitesimal generator $A_\psi$.

\section{{\bf Riesz transforms}}
\label{DFEst}

In this section we shall focus on our dimension free estimates for noncommutative Riesz transforms. More specifically, we will prove Theorems A1 and A2. We shall also illustrate our results with a few examples which provide new estimates both in the commutative and in the noncommutative settings.

\subsection{Khintchine inequalities}

Our results rely on Pisier's method \cite{PRiesz} and a modified version of Lust-Piquard/Pisier's noncommutative Khintchine inequalities \cite{Lu, LuP0}. Given a noncommutative measure space $(\M,\varphi)$, we set $RC_p(\M)$ as the closure of finite sequences in $L_p(\M)$ equipped with the norm $$\big\| (f_k) \big\|_{RC_p(\M)} = \begin{cases} \displaystyle \inf_{f_k = g_k + h_k} \mbox{$ \Big\| \big( \summ_k g_k^* g_k \big)^\frac12 \hskip1pt \Big\|_p + \Big\| \big( \summ_k h_k h_k^* \big)^\frac12 \Big\|_p$} & \mbox{if $1 \le p \le 2$,} \\ \hskip10pt \max \hskip2.5pt \Big\{ \Big\| \big( \summ_k f_k^* f_k \big)^\frac12 \Big\|_p \hskip2pt , \hskip2pt \Big\| \big( \summ_k f_k f_k^* \big)^\frac12 \Big\|_p \Big\} & \mbox{if $2 \le p < \infty$.} \end{cases}$$ The noncommutative Khintchine inequality reads as $G_p(\M)=RC_p(\M)$, where $G_p(\M)$ denotes the closed span in $L_p(\Omega, \mu; L_p(\M))$ of a family $(\gamma_k)$ of centered independent gaussian variables on $(\Omega,\mu)$. The specific statement for $1 \le p < \infty$ is $$\Big( \int_\Omega \Big\| \summ_k \gamma_k(w) f_k \Big\|_{L_p(\M)}^p \, d \mu(w) \Big)^{\frac1p} \, \sim_{c(p)} \, \big\| (f_k) \big\|_{RC_p(\M)}.$$ Our goal is to prove a similar result adding a group action to the picture. Let $\H$ be a separable real Hilbert space. Choosing a orthonormal basis $(e_j)_{j \ge 1}$, we consider the linear map $B: \H \to L_2(\Omega,\mu)$ given by $B(e_j) = \gamma_j$. Let $\Sigma$ stand for smallest $\sigma$-algebra making all the $B(h)$'s measurable. Then, the well known gaussian measure space construction \cite{CCJJV} tells us that, for every real unitary $\alpha$ in $O(\H)$, we can construct a measure preserving automorphism $\beta$ on $L_2(\Omega, \Sigma, \mu)$ such that $\beta(B(h)) = B(\alpha(h))$. Now, assume that a discrete group $\G$ acts by real unitaries on $\H$ and isometrically on some finite von Neumann algebra $\M$. In particular, $\G$ acts isometrically on $L_\infty(\Omega,\Sigma,\mu) \bar\otimes \M$ and we may consider the space $G_p(\M) \rtimes \G$ of operators of the form $$\sum_{g \in \G} \underbrace{\sum_{j \ge 1} \big( B(e_j) \otimes f_{g,j} \big)}_{f_g} \rtimes \lambda(g) \in L_p(\mathcal{A})$$ with $\mathcal{A} = (L_\infty(\Omega,\Sigma,\mu) \bar\otimes \M) \rtimes \G$ and $f_g \in G_p(\M)$. We will also need the conditional expectation $\mathsf{E}_{\M \rtimes \G} ( f_g \rtimes \lambda(g) ) = ( \int_{\Omega} f_g d\mu ) \rtimes \lambda(g)$, which takes $L_p (\mathcal{A})$ contractively to $L_p(\M \rtimes \G)$. The conditional $L_p$ norms $$L_p^{rc}(\mathsf{E}_{\M \rtimes \G}) \, = \, \begin{cases} L_p^r(\mathsf{E}_{\M \rtimes \G}) + L_p^c(\mathsf{E}_{\M \rtimes \G}) & \mbox{if} \ 1 \le p \le 2, \\ L_p^r(\mathsf{E}_{\M \rtimes \G}) \hskip1pt \cap L_p^c(\mathsf{E}_{\M \rtimes \G}) & \mbox{if} \ 2 \le p < \infty, \end{cases}$$ are determined by $$\|f\|_{L_p^r(\mathsf{E}_{\M \rtimes \G})} = \big\| \mathsf{E}_{\M \rtimes \G}(ff^*)^{\frac12} \big\|_p \quad \mbox{and} \quad \|f\|_{L_p^c(\mathsf{E}_{\M \rtimes \G})} = \| \mathsf{E}_{\M \rtimes \G}(f^*f)^{\frac12} \|_p,$$ \cite{JP1} for more precise definitions of these spaces. Define $RC_p(\M) \rtimes \G$ as the gaussian space $G_p(\M) \rtimes \G$, with the norm inherited from $L_p^{rc}(\mathsf{E}_{\M \rtimes \G})$. Then we may generalize the noncommutative Khintchine inequality as follows.

\begin{theorem} \label{qqq}
If $1 < p < \infty$, we have $$C_1 \sqrt{\frac{p-1}{p}} \, \|f\|_{RC_p(\M) \rtimes \G} \, \le \, \|f\|_{G_p(\M) \rtimes \G} \, \le \, C_2 \sqrt{p} \, \|f\|_{RC_p(\M) \rtimes \G}.$$ $G_p(\M) \rtimes \G$ is complemented in $L_p\big( L_{\infty}(\Omega, \Sigma, \mu) \bar\otimes \M \rtimes \G \big)$ with $c(p) \sim \sqrt{p^2 / p-1}$.
\end{theorem}

\dem Let us first assume $p > 2$, the case $p=2$ is trivial. Then, the lower estimate holds with constant $1$ from the continuity of the conditional expectation on $L_{p/2}$. The upper estimate relies on a suitable application of the central limit theorem. Indeed, assume first that $f$ is a finite sum $\sum_{g,h} (B(h) \otimes f_{g,h}) \rtimes \lambda(g)$. Fix $m \ge 1$, use the diagonal action (copying the original action on $\H$ entrywise) on $\ell_2^m(\H)$ and repeat the gaussian measure space construction on the larger Hilbert space resulting in a map $B_m: \ell_2^m(\H) \to L_2(\Omega_m, \Sigma_m, \mu_m)$. Let $\phi_m: \H \to \ell_2^m(\H)$ denote the isometric diagonal embedding $h \mapsto \frac{1}{\sqrt{m}} \sum_j h \otimes e_j$ and $F_1, F_2, \ldots, F_k$ be bounded functions on the complex plane. Then $$\pi \big( F_1(B(h_1)) \cdots F_k(B(h_k)) \big) \, = \, F_1(B_m(\phi_m(h_1)) \cdots F_k(B_m(\phi_m(h_k))$$ extends to a measure preserving $*$-homomorphism $L_{\infty}(\Omega,\Sigma,\mu) \to L_{\infty}(\Omega_m, \Sigma_m, \mu_m)$ which is in addition $\G$-equivariant, i.e. $\pi(\beta_g(f))=\beta_g(\pi(f))$. Thus, we obtain a trace preserving isomorphism $\pi_\G = (\pi \otimes id_\M) \rtimes id_\G$ from $(L_{\infty}(\Omega, \Sigma, \mu) \bar\otimes \M) \rtimes \G$ to the larger space $(L_{\infty}(\Omega_m, \Sigma_m, \mu_m) \bar\otimes \M) \rtimes \G$. This implies $$\|f\|_{G_p(\M) \rtimes \G} \ = \ \Big\| \frac{1}{\sqrt{m}} \sum_{j=1}^m \sum_{g,h} \big( B_m(h \otimes e_j) \otimes f_{g,h} \big) \rtimes \lambda(g) \Big\|_p.$$ The random variables $f_j = \sum_{g,h} (B_m(h \otimes e_j) \otimes f_{g,h}) \rtimes \lambda(g)$ are clearly mean-zero and independent over $\mathsf{E}_{\M \rtimes \G}$, see Appendix A for precise definitions. Hence, the noncommutative Rosenthal inequality ---\eqref{Rosenthal} in the Appendix--- yields 
\begin{eqnarray*}
\lefteqn{\|f\|_{G_p(\M) \rtimes \G}} \\ & \le & \frac{C p}{\sqrt{m}} \Big[ \Big( \sum_{j=1}^m \|f_j\|_p^p \Big)^{\frac{1}{p}} + \Big\| \Big( \sum_{j=1}^m \mathsf{E}_{\M \rtimes \G}(f_j^* f_j) \Big)^{\frac12} \Big\|_p + \Big\| \Big( \sum_{j=1}^m \mathsf{E}_{\M \rtimes \G}(f_j f_j^*) \Big)^{\frac12} \Big\|_p \Big].
\end{eqnarray*}
Note that $\mathsf{E}_{\M \rtimes \G}(f_j f_j^*) = \mathsf{E}_{\M \rtimes \G}(f f^*)$ and $\mathsf{E}_{\M \rtimes \G}(f_j^* f_j) = \mathsf{E}_{\M \rtimes \G}(f^* f)$ for all $j$. Moreover, we also have $\|f_j\|_p = \|f\|_p$. Therefore, the second inequality follows sending $m \to \infty$. An improved Rosenthal's inequality \cite{JZ} actually yields $$\|f\|_{G_p(\M) \rtimes \G} \le C \sqrt{p} \, \|f\|_{RC_p(\M) \rtimes \G},$$ which provides the correct order of the constant in our Khintchine inequality. 

Let us now consider the case $1 < p < 2$. We will proceed by duality as follows. Define the gaussian projection by $$Q(f) = \summ_k \big( \int_\Omega f \gamma_k \, d\mu \big) \gamma_k,$$ which is independent of the choice of the basis. Let $\widehat{Q} = (Q \otimes id_\M) \rtimes id_\G$ be the amplified gaussian projection on $L_p\big( L_{\infty}(\Omega, \Sigma, \mu) \bar\otimes \M \rtimes \G \big)$. It is clear that $G_p(\M) \rtimes \G$ is the image of this $L_p$ space under the gaussian projection. Similarly $RC_p(\M) \rtimes \G$ is the image of $L_p^{rc}(\mathsf{E}_{\M \rtimes \G})$. Note that $$\widehat{Q}: L_p^{rc}(\mathsf{E}_{\M \rtimes \G}) \to RC_p(\M) \rtimes \G$$ is a contraction. Indeed, we have $$\mathsf{E}_{\M \rtimes \G}(ff^*) = \mathsf{E}_{\M \rtimes \G}(\widehat{Q}f\widehat{Q}f^*) + \mathsf{E}_{\M \rtimes \G}(\widehat{Q}^\perp f \widehat{Q}^\perp f^*) \ge \mathsf{E}_{\M \rtimes \G}(\widehat{Q}f\widehat{Q}f^*)$$ by orthogonality and the same holds for the column case. This immediately gives $$\|f\|_{RC_p(\M) \rtimes \G} = \sup_{\|g\|_{RC_q} \le 1} \big| \mathrm{tr} (fg) \big| \le \Big( \sup_{\|g\|_{RC_q} \le 1} \|g\|_{G_q(\M) \rtimes \G} \Big) \, \|f\|_{G_p(\M) \rtimes \G}$$ with $1/p + 1/q = 1$. In conjunction with our estimates for $q \ge 2$, this proves the lower estimate for $p \le 2$. The upper estimate is a consequence of the continuous inclusion $L_p^{rc}(\mathsf{E}_{\M \rtimes \G}) \to L_p$ for $p \le 2$, see \cite[Theorem 7.1]{JX1}. 

It remains to prove the complementation result. Since the gaussian projection is self-adjoint we may assume $p \le 2$. Moreover, the upper estimate in the first assertion together with the contractivity of the gaussian projection on $L_p^{rc}(\mathsf{E}_{\M \rtimes \G})$ give rise to
\begin{eqnarray*}
\big\| \widehat{Q}f \big\|_{G_p(\M) \rtimes \G} & \lesssim & \big\| \widehat{Q}f \big\|_{RC_p(\M) \rtimes \G} \\ [4pt] & = & \sup_{\|g\|_{RC_q} \le 1} \big| \mathrm{tr} (\widehat{Q}f g) \big| \ = \ \sup_{\|g\|_{RC_q} \le 1} \big| \mathrm{tr} (f \widehat{Q} g) \big| \\ & \le &  \sup_{\|g\|_{RC_q} \le 1} \big\| \widehat{Q} g \big\|_{G_q(\M) \rtimes \G} \|f\|_{L_p(\mathcal{A})} \ \le \ C_2 \sqrt{q} \, \|f\|_{L_p(\mathcal{A})}.
\end{eqnarray*}
The last inequality follows from the first assertion for $q$, the proof is complete. \fin

\subsection{Riesz transforms in Meyer form} \label{DFRiesz}

Denote by $\lambda$ the Lebesgue measure and let us write $\gamma$ for the normalized gaussian measure in $\R^n$. With this choice, the maps $\beta_t f (x,y) = f(x+ty)$ are measure preserving $*$-homomorphisms from $L_\infty(\R^n,\lambda)$ to $L_{\infty}(\R^n \times \R^n, \lambda \times \gamma)$. We may also replace $\lambda$ by the Haar measure $\nu$ on the Bohr compactification $\R_{\mathrm{bohr}}^n$ of $\R^n$, this latter case including $n=\infty$. Moreover, if $\G$ acts on $\R^n$ then $\beta_t$ commutes with the diagonal action. As we already recalled in the Introduction, \eqref{PisierFormula} takes the form $$\sqrt{\frac{2}{\pi}} \, \delta (-\Delta)^{- \frac12} f = (id_{\R^n} \otimes Q) \Big( \mbox{p.v.} \, \frac{1}{\pi} \int_\R \beta_t f \frac{dt}{t} \Big),$$ with $Q$ the gaussian projection and $\delta: C^\infty(\R^n) \to C^\infty(\R^n \times \R^n)$ the derivation $$\delta(f)(x,y) = \sum_{j=1}^n \frac{\partial f}{\partial x_j} y_j = \big\langle \nabla f(x), y \big\rangle = B(\nabla f(x)) (y).$$

\begin{lemma} \label{kr}
If $1 < p< \infty$, we have $$\delta (-\Delta)^{- \frac12} \rtimes id_\G: L_p \big( L_\infty(\R_{\mathrm{bohr}}^n, \nu) \rtimes \G \big) \to L_p \big( L_\infty(\R_{\mathrm{bohr}}^n \times \R^n, \nu \times \gamma) \rtimes \G \big)$$ with norm bounded by $C p^3 / (p-1)^{3/2}$. Moreover, the same holds when $n = \infty$.
\end{lemma}

\dem The cross product extension of \eqref{PisierFormula} reads as $$\sqrt{\frac{2}{\pi}} \, \big( \delta (-\Delta)^{- \frac12} \rtimes id_\G \big) f \, = \, \big( (id_{\R^n} \otimes Q) \rtimes id_\G \big) \Big( \mathrm{p.v.} \, \frac{1}{\pi} \int_\R ( \beta_t \rtimes id_\G) f \frac{dt}{t} \Big).$$ This gives $\delta (-\Delta)^{- \frac12} \rtimes id_\G = \sqrt{\pi/2} \, \widehat{Q} (H \rtimes id_\G)$ where $H$ is the transferred Hilbert transform $$Hf(x,y) \, = \, \mbox{p.v.} \, \frac{1}{\pi} \int_\R \beta_t f(x,y) \frac{dt}{t}.$$ By de Leeuw's theorem \cite{dL}, the Hilbert transform is bounded on $L_p(\R_{\mathrm{bohr}},\nu)$ with the same constants as in $L_p(\R,\lambda)$. The operator above can be seen as a directional Hilbert transform at $x$ in the direction of $y$, which also preserves the same constants for $y$ fixed. In particular, a Fubini argument combined with a gaussian average easily gives that $H: L_p(\R_{\mathrm{bohr}}^n,\nu) \to L_p(\R_{\mathrm{bohr}}^n \times \R^n, \nu \times \gamma)$ again with the classical constants, even for $n = \infty$. To analyze the crossed product $H \rtimes id_\G$ we note that $\beta_t \rtimes id_\G$ is a trace preserving $*$-automorphism on the algebra $L_\infty(\R_{\mathrm{bohr}}^n \times \R^n, \nu \times \gamma) \rtimes \G$. According to the Coifman-Weiss transference principle \cite{CW} and the fact that $H$ is $\G$-equivariant, we see that $H \rtimes id_\G$ extends to a bounded map on $L_p$ with constant $c(p) \sim p^2 / p-1$. Indeed, it is straightforward that the proof of the transference principle for one-parameter automorphisms translates verbatim to the present setting. Then, the assertion follows from the complementation result in Theorem \ref{qqq} for the gaussian projection. The proof is complete. \fin 

\renewcommand{\theequation}{$\psi$-Riesz}
\addtocounter{equation}{-1}

Given a length $\psi: \G \to \R_+$, consider the derivation $\delta_\psi: \V \to L_\infty(\R^n, \gamma) \rtimes \G$ determined by $\delta_\psi(\lambda(g)) = B(b_\psi(g)) \rtimes \lambda(g)$ where $b_\psi: \G \to \R^n$ is the cocycle map associated to $\psi$ and the crossed product is defined via the cocycle action $\alpha_\psi: \G \to \mathrm{Aut}(\R^n)$. We include the case $n=\infty$, so that any length function/cocycle is admissible. The Leibniz rule $\delta_\psi(\lambda(gh)) = \delta_\psi(\lambda(g)) \lambda(h) + \lambda(g) \delta_\psi(\lambda(h))$ follows from the cocycle law. Given the infinitesimal generator $A_\psi(\lambda(g)) = \psi(g) \lambda(g)$, our estimates for Riesz transforms read as follows
\begin{equation} \label{OneLineDimFree}
\big\| A_\psi^{\frac12}f \big\|_{L_p(\widehat{\mathbf{G}})} \, \sim_{c_1(p)} \, \|\delta_\psi f\|_{G_p(\C) \rtimes \G} \, \sim_{c_2(p)} \, \|\delta_\psi f\|_{RC_p(\C) \rtimes \G}.
\end{equation}

\numberwithin{equation}{section}

\demAA According to the definition of the norm in $RC_p(\C) \rtimes \G$ it suffices to prove our one-line formulation \eqref{OneLineDimFree} above. The second norm equivalence follows from the Khintchine inequality in Theorem \ref{qqq} with constant $c_2(p) \sim \sqrt{p^2/p-1}$, let us prove the first one. The map $$\pi: \lambda(g) \in \V \mapsto \exp \big( 2 \pi i \langle b_\psi(g), \cdot \rangle \big) \rtimes \lambda(g) \in L_\infty(\R^n_{\mathrm{bohr}}, \nu) \rtimes \G$$ is a trace preserving $*$-homomorphism which satisfies $$\big( \delta (-\Delta)^{-\frac12} \rtimes id_\G \big) \circ \pi \, = \, i \big( id_{\R^n} \rtimes \pi \big) \circ \delta_\psi A_\psi^{-\frac12}.$$ Indeed, if we let the left hand side act on $\lambda(g)$ we obtain
\begin{eqnarray*}
\big( \delta (-\Delta)^{-\frac12} \rtimes id_\G \big) \circ \pi (\lambda(g)) & = & \delta (-\Delta)^{-\frac12} \Big( \exp \big( 2 \pi i \langle b_\psi(g), \cdot \rangle \big) \Big) \rtimes \lambda(g) \\ [2pt]& = & \frac{1}{2 \pi \, \|b_\psi(g)\|_{\H_\psi}} \delta \Big( \exp \big( 2 \pi i \langle b_\psi(g), \cdot \rangle \big) \Big) \rtimes \lambda(g) \\ & = & \frac{i}{\sqrt{\psi(g)}} \exp \big( 2 \pi i \langle b_\psi(g), \cdot \rangle \big) \otimes B(b_\psi(g)) \rtimes \lambda(g),
\end{eqnarray*}
which is $i \big( id_{\R^n} \rtimes \pi \big) \circ \delta_\psi A_\psi^{-\frac12} (\lambda(g))$. By Lemma \ref{kr}, both sides in the intertwining identity above define bounded maps $L_p(\V) \to L_p(L_\infty(\R_{\mathrm{bohr}}^n \times \R^n, \nu \times \gamma) \rtimes \G)$ for $1 < p < \infty$. In particular, we obtain the inequality $$\|\delta_\psi f\|_p = \big\| (id_{\R^n} \rtimes \pi) \delta_\psi A_\psi^{-\frac12} (A_\psi^{\frac12} f ) \big\|_p \lesssim \frac{p^3}{(p-1)^{\frac32}} \big\| A_\psi^{\frac12}f \big\|_p.$$ The constant given above also follows from Lemma \ref{kr}. The reverse estimate follows with the same constant from a duality argument. Indeed, if we fix $f$ to be a trigonometric polynomial, there exists another trigonometric polynomial $f'$ with $\|f'\|_{p'} = 1$ and such that $$(1 - \varepsilon) \big\| A_\psi^\frac12 f \big\|_p \le \tau \big( f' A_\psi^\frac12 f \big).$$ Note that $A_\psi^{-1/2}$ is only well-defined on $f'' = \sum_{\psi(g) \neq 0} \widehat{f}'(g) \lambda(g)$. However, since $\G_0 = \{ g \in \G \, | \, \psi(g) = 0 \}$ is a subgroup, we may consider the associated conditional expectation $E_{\G_0}$ on $\V$ and obtain $f'' = f' - E_{\G_0}f'$ so that $\|f''\|_{p'} \le 2$. On the other hand, we note the crucial identity
\begin{eqnarray*}
\hskip-15pt \lefteqn{\mathrm{tr}_{L_\infty(\Omega) \rtimes \G} \big( \delta_\psi f_1^* \delta_\psi f_2 \big)} \\ [5pt] & = & \sum_{g_1, g_2} \overline{\widehat{f}_1(g_1)} \widehat{f}_2(g_2) \Big( \int_\Omega \alpha_{g_1^{-1}} \big( \overline{B(b_\psi(g_1))} B(b_\psi(g_2)) \big) \, d\mu \Big) \tau(\lambda(g_1^{-1} g_2)) \\ & = & \sum_{g \in \G} \overline{\widehat{f}_1(g)} \widehat{f}_2(g) \big\langle b_\psi(g), b_\psi(g) \big\rangle_{\psi} \, = \, \sum_{g \in \G} \overline{\widehat{f}_1(g)} \widehat{f}_2(g) \psi(g) \, = \, \tau(A_\psi^\frac12 f_1^* A_\psi^\frac12 f_2).
\end{eqnarray*}
Combining both results we get
\begin{eqnarray*}
\big\| A_\psi^\frac12 f \big\|_p & \le & \frac{1}{1-\varepsilon} \tau \big( f' A_\psi^\frac12 f \big) = \frac{1}{1-\varepsilon} \tau \big( f'' A_\psi^\frac12 f \big) \\ & = & \frac{1}{1-\varepsilon} \ \mathrm{tr}_{L_\infty(\Omega) \rtimes \G} \ \big( \delta_\psi (A_\psi^{-\frac12} f'') \delta_\psi f \big) \\ & \le & \frac{1}{1 - \varepsilon} \big\| \delta_\psi A_\psi^{-\frac12} f'' \big\|_{p'} \, \|\delta_\psi f\|_p \ \lesssim \ \frac{p^3}{(p-1)^{\frac32}} \, \|\delta_\psi f\|_p.
\end{eqnarray*}
The last estimate follows from Lemma \ref{kr} and yields $c_1(p) \lesssim p^3/(p-1)^{3/2}$. \fin

\begin{remark}
\emph{Let us write $\Gamma_\psi(f_1,f_2)$ for the gradient form $$\Gamma_\psi(f_1,f_2) = \frac12 \Big( A_\psi(f_1^*)f_2 + f_1^*A_\psi(f_2) - A_\psi(f_1^*f_2) \Big).$$ In the Introduction we related Theorem A2 with Meyer's formulation in terms of $\Gamma_\psi$ via the identity $\mathsf{E}_\V(\delta_\psi f_1^* \delta_\psi f_2) = \Gamma_\psi(f_1,f_2)$. The proof is not difficult. Arguing as above, we find $$\mathsf{E}_\V \big( \delta_\psi f_1^* \delta_\psi f_2 \big) \, = \, \summ_{g,h} \overline{\widehat{f}_1(g)} \widehat{f}_2(h) \, \big\langle b_\psi(g), b_\psi(h) \big\rangle_{\psi} \, \lambda(g^{-1}h) \, = \, \Gamma_\psi(f_1,f_2)$$ since $\langle b_\psi(g), b_\psi(h) \rangle_{\psi} = \frac12 (\psi(g) + \psi(h) - \psi(g^{-1}h))$, as explained in Appendix A.}
\end{remark}

\begin{remark} \label{NoEasyInf}
\emph{When $1 < p < 2$, we may consider decompositions of $f = f_1 + f_2$ so that $\delta_\psi f = \phi_1 + \phi_2$ with $\phi_j = \delta_\psi f_j$ in our result. These particular decompositions give rise to $$\big\| A_\psi^\frac12 f \big\|_p \, \le \, c(p) \, \inf_{f = f_1 + f_2} \big\| \Gamma_\psi(f_1,f_1) \big\|_p + \big\|\Gamma_\psi(f_2^*,f_2^*) \big\|_p.$$ Somehow surprisingly, the reverse inequality does not hold. Indeed, using the arguments in the next section this would imply that Theorem A1 holds with the untwisted operators $(b_j)$, but this was already disproved in \cite{LuP} by F. Lust-Piquard}.
\end{remark}

\begin{remark} \label{BestConstant}
\emph{Our constants grow like $p^{3/2}$ as $p \to \infty$ with a dual behavior as $p \to 1$. According to the results in the literature, one might expect a linear growth of the constant $\sim p$. It is however not clear to us whether this is true in our context since we admit semigroups which are not diffusion in the sense of Bakry, like the Poisson semigroup. It is an interesting problem to determine the optimal behavior of dimension free constants for (say) the Riesz transform associated to the Poisson semigroup $(e^{-t \sqrt{-\Delta}})$ in $\R^n$. On the other hand, it seems that $p^{3/2}$ is best known for dimension free estimates of Riesz transforms acting on matrix-valued functions in $\R^n$, even for the classical ones associated with the heat semigroup $(e^{t \Delta
})$.}
\end{remark}

\subsection{Riesz transforms in cocycle form}

We are now ready to prove Theorem A1. The main ingredient comes from a factorization of the conditional expectation $\mathsf{E}_\V: L_\infty(\Omega,\mu) \rtimes \G \to \V$ in terms of certain right $\V$-module map. As predicted by Hilbert module theory, this factorizations is always possible. In our case, when $\phi_1, \phi_2 \in G_p(\C) \rtimes \G$ it takes the form $$\mathsf{E}_\V(\phi_1^* \phi_2) = u(\phi_1)^* u(\phi_2),$$ where $u: G_p(\C) \rtimes \G \to C_p(L_p(\widehat{\mathbf{G}}))$ is defined as follows $$u(\phi) \, = \, u \Big( \sum_{g \in \G} \phi_g \rtimes \lambda(g) \Big) \, = \, \sum_{j \ge 1} \Big( \underbrace{\sum_{g \in \G} \big\langle \phi_g, B(e_j) \big\rangle_{L_2(\Omega,\mu)} \lambda(g)}_{u_j(\phi)} \Big) \otimes e_{j1}.$$ Before proving Theorem A1, it is convenient to explain where the infimum is taken and how do we define the twisted form of $(b_j)$. The infimum $R_{\psi,j}f = a_j + b_j$ runs over all possible families of the form $a_j = u_j(\phi_1)$ and $b_j = u_j(\phi_2)$ for some $\phi_1, \phi_2 \in G_p(\C) \rtimes \G$. This is equivalent to require that $\sum_j B(e_j) \rtimes a_j \in G_p(\C) \rtimes \G$ and the same for $(b_j)$. Once this is settled, if we note that $$b_j \, = \, \sum_{g \in \G} \Big\langle \sum_{k \ge 1} \widehat{b}_k(g) e_k , e_j \Big\rangle_{\H_\psi} \lambda(g),$$ then the twisted form of the family $(b_j)_{j \ge 1}$ is determined by the formula $$\widetilde{b}_j \, = \, \sum_{g \in \G} \Big\langle \sum_{k \ge 1} \widehat{b}_k(g) e_k , \alpha_{\psi,g}(e_j) \Big\rangle_{\H_\psi} \lambda(g).$$

\demA If $f \in L_p^\circ(\widehat{\mathbf{G}})$ with $1 < p \le 2$
\begin{eqnarray*}
\|f\|_p & = & \big\| A_\psi^\frac12 A_\psi^{-\frac12} f \big\|_p \\ & \sim & \inf_{\delta_\psi A_\psi^{-\frac12} f = \phi_1 + \phi_2} \big\| \mathsf{E}_\V(\phi_1^* \phi_1)^\frac12 \big\|_p + \big\| \mathsf{E}_\V(\phi_2 \phi_2^*)^\frac12 \big\|_p
\end{eqnarray*} 
by Theorem A2. By the factorization of $\mathsf{E}_\V$, we observe that 
\begin{eqnarray*}
\mathsf{E}_\V(\phi_1^* \phi_1)^\frac12 & = & |u(\phi_1)| \ = \ \Big( \sum_{j \ge 1} u_j(\phi_1)^* u_j(\phi_1) \Big)^\frac12, \\ \mathsf{E}_\V(\phi_2 \phi_2^*)^\frac12 & = & |u(\phi_2^*)| \ = \ \Big( \sum_{j \ge 1} u_j(\phi_2^*)^* u_j(\phi_2^*) \Big)^\frac12.
\end{eqnarray*}
On the other hand, we may write $u(\delta_\psi A_\psi^{-\frac12} f)$ in the familiar way $$u(\delta_\psi A_\psi^{-\frac12} f) \, = \, \sum_{j \ge 1} \Big( \sum_{g \in \G} \frac{\langle b_\psi(g), e_j \rangle_{\H_\psi}}{\sqrt{\psi(g)}} \widehat{f}(g) \lambda(g) \Big) \otimes e_{j1} \, = \, \sum_{j \ge 1} R_{\psi,j}f \otimes e_{j1}.$$ This and the injectivity of $u$ gives that $\|f\|_p$ is comparable to the norm
$$\inf_{R_{\psi,j}f = u_j(\phi_1) + u_j(\phi_2)} \Big\| \Big( \sum_{j \ge 1} u_j(\phi_1)^* u_j(\phi_1) \Big)^\frac12 \Big\|_p + \Big\| \Big( \sum_{j \ge 1} u_j(\phi_2^*)^* u_j(\phi_2^*) \Big)^\frac12 \Big\|_p.$$ Therefore, if $a_j = u_j(\phi_1)$ and $b_j = u_j(\phi_2)$ it suffices to see that $\widetilde{b}_j = u_j(\phi_2^*)^*$ to settle the case $1 < p \le 2$. Since $u$ is injective and $\phi_2 = u^{-1}(\sum_k b_k \otimes e_{k1})$, we may prove such identity as follows 
\begin{eqnarray*}
u_j(\phi_2^*)^* & = & u_j \Big[ \Big( \sum_{g \in \G} \big( \sum_{k \ge 1} \widehat{b}_k(g) B(e_k) \big) \rtimes \lambda(g) \Big)^* \Big]^* \\ & = & u_j \Big[ \sum_{g \in \G} \big( \sum_{k \ge 1} \overline{\widehat{b}_k(g)} B(\alpha_{\psi,g^{-1}}(e_k)) \big) \rtimes \lambda(g^{-1}) \Big]^* \\ & = & \Big[ \sum_{g \in \G} \Big \langle \sum_{k \ge 1} \widehat{b}_k(g) B(\alpha_{\psi,g^{-1}}(e_k)), B(e_j) \Big\rangle \lambda(g^{-1}) \Big]^* \\ & = & \sum_{g \in \G} \Big \langle \alpha_{\psi,g^{-1}} \Big( \sum_{k \ge 1} \widehat{b}_k(g) e_k \Big), e_j \Big\rangle_{\H_\psi} \lambda(g) \ = \ \widetilde{b}_j.
\end{eqnarray*}
Recall that we have implicitly used that $\H_\psi$ and $L_2(\Omega,\mu)$ are real Hilbert spaces. It remains to consider the case $2 \le p < \infty$. However, this case is simpler since arguing as above we obtain $$\mathsf{E}_\V \big( \delta_\psi A_\psi^{-\frac12} f^* \delta_\psi A_\psi^{-\frac12} f \big) \, = \, \big| u (\delta_\psi A_\psi^{-\frac12} f) \big| \, = \, \Big( \sum_{j \ge 1} (R_{\psi,j}f)^* (R_{\psi,j}f) \Big)^\frac12.$$ The row term in Theorem A1 arises replacing $f$ by $f^*$. The proof is complete. \fin

\begin{remark} \label{HilbertModule}
\emph{Let $A$ be the infinitesimal generator of a given Markovian semigroup $\S = (\S_t)_{t \ge 0}$ acting on $(\M,\tau)$. Sauvageot's theorem \cite{Sa} provides a factorization $A = \delta^* \delta$ in terms of certain symmetric derivation $\delta: L_2(\M) \to \mathsf{M}$ taking values in some Hilbert $\M$-bimodule $\mathsf{M}$. As a consequence of our results, it is the nature of the tangent module $\mathsf{M}$ and not of $\M$ itself what dictates the behavior of Riesz transforms on $L_p(\M,\tau)$, in the sense that we find noncommutative phenomena as long as $\mathsf{M}$ is not abelian regardless the nature of $\M$.}
\end{remark}

\begin{remark} \label{RemRowColumnRiesz}
\emph{Theorems A1 and A2 are formulated for left cocycles, although an alternative form is possible for right cocycles or both together, see the precise definitions in Appendix A. The only (cosmetic) change appears in the statement of Theorem A1 for $p \ge 2$. The row version of $R_{\psi,j}$ is $$R_{\psi,j}'f \, = \, 2 \pi i \sum_{g \in \G} \frac{\langle b_\psi(g^{-1}), e_j \rangle_{\H_\psi}}{\sqrt{\psi(g)}} \widehat{f}(g) \lambda(g).$$ Note that $R_{\psi,j}(f^*) = - R_{\psi,j}'(f)^*$ so that $\sum_j (R_{\psi,j}f^*)^* (R_{\psi,j}f^*)$ can be written as a row square function in terms of the row Riesz transforms $R_{\psi,j}'$. Although these two formulations collapse for $\G$ abelian (left and right cocycles coincide) the statement does not simplifies for non-trivial cocycle actions and $1 < p < 2$.}
\end{remark}

\begin{remark} \label{LCA}
\emph{Our discreteness assumption on $\G$ is not essential in Theorems A1 and A2. In fact, it will be instrumental to observe that our arguments can be easily modified to make them work for $\G = \R^n$ or any other LCA group. Moreover, we may extend them to unimodular locally compact groups. Let us sketch how to modify the argument. Let $\G$ be a unimodular locally compact group and $\psi: \G \to \R_+$ a conditionally negative continuous length function, so that we may associate to it a cocycle $(\H_\psi, \alpha_\psi, b_\psi)$ were both the action $\alpha_\psi$ and the cocycle map $b_\psi$ are continuous. Then, we may consider the semidirect product $\Gamma = \H_{\psi} \rtimes_\alpha \G$. Here $\H_{\psi}$ carries the discrete topology and hence $\Gamma$ is again locally compact. We have a $*$-homomorphism $\pi: g \in \G \mapsto b(g) \rtimes g \in \Gamma$ which extends to a trace preserving map on $\V$. Indeed, recall that the trace on $\V$ is given by $$\tau_\G \big( \int_\G \widehat{f}(g) \lambda(g)d\mu(g) \big) = \widehat{f}(e)$$ for $\widehat{f} \in \mathcal{C}_c(\G)$ and we obtain $$\tau_\Gamma \big( \int_\G \widehat{f}(g) \lambda(b_\psi(g) \rtimes g) d\mu(g) \big) = \widehat{f}(e) = \tau_\G \big( \int_\G \widehat{f}(g) \lambda(g) d\mu(g) \big).$$  The only delicate point is to justify that $$\widehat{Q} \ {\rm p.v.} \int_{\G \times \R} \widehat{f}(g) \beta_t(\lambda(b_\psi(g) \rtimes g)) d\mu(g) \frac{dt}{t} = c \int_\G \widehat{f}(g) \frac{B(b_\psi(g))}{\sqrt{\psi(g)}} \otimes \lambda(b_\psi(g) \rtimes g)) d\mu(g).$$ Clearly, we have 
\begin{eqnarray*}  
\lefteqn{\hskip-30pt \widehat{Q} \ {\rm p.v.} \int_{\R} \beta_t(\lambda(b_\psi(g) \rtimes g)) \frac{dt}{t}} \\ & = & \widehat{Q} \ {\rm p.v.} \int_{\R} e^{it \langle B(b_\psi(g)), \cdot \rangle} \otimes \lambda(b_\psi(g) \rtimes g) \frac{dt}{t} \\ & = & \widehat{Q} \ \mathrm{p.v.} \int_{\R} i \sin(t \langle B(b_\psi(g)), \cdot \rangle) \frac{dt}{t} \otimes \lambda(b_\psi(g) \rtimes g) \\ & = & \sqrt{\frac{2}{\pi}} i \frac{B(b_\psi(g))}{\sqrt{\psi(g)}} \otimes \lambda(b_\psi(g) \rtimes g). 
\end{eqnarray*} 
As before, we only have to evaluate the principal value for a single gaussian random variable $B(b_\psi(g))$. Moreover, for $f \in L_1(\widehat{\mathbf{G}}) \cap L_2(\widehat{\mathbf{G}})$ we may justify interchanging the integral and then deduce that our formula is valid in $L_2$. Indeed, the functions $F_{\delta}(x) = \int_{\delta}^{1/\delta} \sin(xt)\frac{dt}{t}$ are uniformly bounded in $\delta$, and hence we have almost everywhere convergence. Using $\|\xi\|_p = \sup_{\|\eta\|_{p'} \le 1} |\mathrm{tr}(\xi \eta)|$ where the supremum is taken over $\eta$ in $L_{p'} \cap L_2$, we deduce first that $$\int_\G \widehat{f}(g) \, \frac{B(b_\psi(g))}{\sqrt{\psi(g)}} \otimes \lambda(b_\psi(g) \rtimes g) d\mu(g) \in L_p(\widehat{\mathbf{G}}).$$ Then, we may take the gaussian projection and conclude as before.}
\end{remark}

\begin{remark} \label{OS}
\emph{On the other hand, it is interesting to note that Theorems A1 and A2 hold in the category of operator spaces \cite{P3}. In other words, the same inequalities are valid with matrix Fourier coefficients when the involved operators act trivially on the matrix amplification.}
\end{remark}

\renewcommand{\theequation}{CN}
\addtocounter{equation}{-1}

\subsection{Examples, commutative or not} \label{ExSect}

In order to illustrate Theorems A1 and A2, it will be instrumental to present conditionally negative lengths in a more analytic way. As  it will be justified in Appendix A (Theorem \ref{CharacterizationCN}), these lengths are all deformations of the standard inner cocycle which acts by left multiplication. Namely, $\psi: \G \to \R_+$ is conditionally negative iff it can be written as 
\begin{equation} \label{CNLength}
\psi(g) \, = \, \tau_\psi \big(2 \lambda(e) - \lambda(g) - \lambda(g^{-1}) \big)
\end{equation}
for a positive linear functional $\tau_\psi: \Pi_0 \to \C$ defined on the space $\Pi_0$ of trigonometric polynomials in $\V$ whose Fourier coefficients have vanishing sum. Having this in mind, let us consider some examples illustrating Theorems A1 and A2.  

\numberwithin{equation}{section}

\vskip3pt

\noindent \textbf{A. Fractional laplacians in $\R^n$.} The Riesz potentials $f \mapsto (-\Delta)^{-\beta/2}f$ are classical operators in Euclidean harmonic analysis \cite{St4}. It seems however that dimension free estimates for associated Riesz transforms in $\R^n$ are unexplored. If we let our infinitesimal generator to be $A_{\beta} = (- \Delta)^\beta$ and $p \ge 2$ (for simplicity), the problem in $\R^n$ consists in showing that $$\Big\| \Big( \sum_{j \ge 1} \big| \partial_{\beta,j} A_\beta^{-\frac12} f \big|^2 \Big)^\frac12 \Big\|_p \, \sim_{c(p)} \, \|f\|_p$$ for some differential operators $\partial_{\beta,j}$ and constants $c(p)$ independent of the dimension $n$. $A_\beta$ generates a Markov semigroup for $0 \le \beta \le 1$. Indeed, the non-elementary cases $0 < \beta < 1$ require to know that the length $\psi_\beta(\xi) = |\xi|^{2\beta}$ is conditionally negative. Since \eqref{CNLength} holds for $\G = \R^n$, the claim follows from the simple identity $$\psi_\beta(\xi) \, = \, |\xi|^{2\beta} \, = \, \frac{1}{\mathrm{k}_n(\beta)} \int_{\R^n} \big( 2 - e^{2\pi i \langle x, \xi \rangle} - e^{-2\pi i \langle x, \xi \rangle} \big) \, d\mu_\beta(x),$$ with $d\mu_\beta(x) = dx / |x|^{n+2\beta}$ and $$\mathrm{k}_n(\beta) = 2 \int_{\R^n} \Big( 1 - \cos \big( 2 \pi \big\langle \frac{\xi}{|\xi|}, s \big\rangle \big) \Big) \, \frac{ds}{|s|^{n+2\beta}} \sim \frac{\pi^{n/2}}{\Gamma(n/2)} \max \Big\{ \frac{1}{\beta}, \frac{1}{1-\beta} \Big\}.$$ The constant $\mathrm{k}_n(\beta)$ only makes sense for $0 < \beta < 1$ and is independent of $\xi$. The associated cocycle $(\H_\beta, \alpha_\beta, b_\beta)$ is given by the action $\alpha_{\beta,\xi}(f) = \exp(2 \pi i \langle \cdot, \xi \rangle) f$ and the cocycle map $\xi \mapsto 1 - \exp(2 \pi i \langle \cdot, \xi \rangle) \in \H_\beta = L_2(\R^n, \mu_\beta/\mathrm{k}_n(\beta))$. Contrary to the standard Riesz transforms for $\beta = 1$, for which $\H_\beta = \R^n$ and the cocycle is trivial, we need to represent $\R^n$ in infinitely many dimensions to obtain the right differential operators $$\widehat{R_{\beta,j}f}(\xi) \, = \, \widehat{\partial_{\beta,j} A_\beta^{-\frac12} f}(\xi) \, = \, \frac{\langle b_\beta(\xi), e_j \rangle_{\H_\beta}}{|\xi|^{\beta}}\widehat{f}(\xi).$$ In particular, Theorem A1 gives norm equivalences for all $1 < p < \infty$ which differ from the classical statement when $1 < p < 2$, since the cocycle action is not trivial anymore for $0 < \beta < 1$. It is also interesting to look at Theorem A2. Taking into account that $\widehat{A_\beta f}(\xi) = |\xi|^{2\beta} \widehat{f}(\xi)$ and the definition of the associated gradient form $\Gamma_\beta$, we obtain $$\big\| (- \Delta)^{\beta/2} f \big\|_p \sim \Big\| \Big( \int_{\R^n} M_{f,\beta}(\cdot,\xi) e^{2\pi i \langle \cdot, \xi \rangle} d\xi \Big)^\frac12 \Big\|_p$$ for $f$ smooth enough and $$\displaystyle M_{f,\beta}(x,\xi) = \frac12 \Big( f(x) \overline{\widehat{f}(-\xi)} + \overline{f(x)} \widehat{f}(\xi) - \widehat{|f|^2}(\xi) \Big) |\xi|^{2\beta}.$$ More applications on Euclidean $L_p$ multipliers will be analyzed in the next section.

\vskip3pt

\noindent \textbf{B. Discrete laplacians in LCA groups.} Let $\Gamma_0$ be a locally compact abelian group and $s_0 \in \Gamma_0$ be torsion free. If $\partial_j f(\gamma) = f(\gamma) - f (\gamma_1, \ldots, s_0 \gamma_j, \ldots, \gamma_n)$ stand for discrete directional derivatives in $\Gamma = \Gamma_0 \times \Gamma_0 \times \cdots \times \Gamma_0$, we may construct  the laplacian $\mathcal{L} = \sum \partial_j^* \partial_j$. Lust-Piquard's main result in \cite{LuP} establishes dimension free estimates in this context for the family of discrete Riesz transforms given by $R_j = \partial_j \mathcal{L}^{-1/2}$ and $R_j^* = \mathcal{L}^{-1/2} \partial_j^*$. If $p \ge 2$, we obtain $$\Big\| \Big( \sum_{j=1}^n |R_j f|^2 + |R_j^* f|^2 \Big)^\frac12 \Big\|_{L_p(\Gamma)} \, \sim_{c(p)} \, \|f\|_{L_p(\Gamma)}.$$ It is not difficult to recover and generalize Lust-Piquard's theorem from Theorem A1. Indeed, let $\Gamma$ be any LCA group equipped with a positive measure $\mu_\psi$. If $\G$ denotes the dual group of $\Gamma$, let us write $\chi_g: \Gamma \to \mathbb{T}$ for the associated characters and $\nu$ for the Haar measure on $\G$. Define
\begin{eqnarray*}
A_\psi f & = & \int_\G \widehat{f}(g) A_\psi (\chi_g) d\nu(g) \\ & = & \int_\G \widehat{f}(g) \Big[ \underbrace{\int_\Gamma \big( 2 \chi_e - \chi_{g} - \chi_{g^{-1}} \big)(\gamma) d\mu_\psi(\gamma)}_{\psi(g)} \Big] \chi_g \, d\nu(g).
\end{eqnarray*} 
If $\psi: \G \to \R_+$, then it is a c.n. length which may be represented by the cocycle 
$$\Big( \H_\psi, \alpha_{\psi,g}, b_{\psi}(g) \Big) \, = \, \Big( L_2(\Gamma,\mu_\psi), \chi_g \cdot [ \ ], 1 -\chi_g \Big).$$ In other words, we have $\psi(g) = \langle b_{\psi}(g), b_{\psi}(g) \rangle_{\H_\psi}$. Of course, we may regard $\H_\psi$ as a real Hilbert space by identifying $\chi_g(\gamma) \in \C \mapsto (\mathrm{Re}(\chi_g(\gamma)), \mathrm{Im}(\chi_g(\gamma))) \in \R^2$ and the product $\chi_g(\gamma) \cdot [ \ ]$ with a rotation by $\mathrm{arg}(\chi_g(\gamma))$. Let us set
\begin{eqnarray*}
R_\gamma^1f & = & \int_\G \frac{\langle b_\psi(g), e_\gamma^1 \rangle_{\H_\psi}}{\sqrt{\psi(g)}} \widehat{f}(g) \chi_g \, d\nu(g) \ = \ \int_\G \frac{\mathrm{Re} \, b_\psi(g,\gamma)}{\sqrt{\psi(g)}} \widehat{f}(g) \chi_g \, d\nu(g), \\ R_\gamma^2f & = & \int_\G \frac{\langle b_\psi(g), e_\gamma^2 \rangle_{\H_\psi}}{\sqrt{\psi(g)}} \widehat{f}(g) \chi_g \, d\nu(g) \ = \ \int_\G \frac{\mathrm{Im} \, b_\psi(g,\gamma)}{\sqrt{\psi(g)}} \widehat{f}(g) \chi_g \, d\nu(g).
\end{eqnarray*}
Then, Theorem A1 takes the following form for LCA groups and $p \ge 2$ $$\Big\| \Big( \int_\Gamma |R_\gamma^1f|^2 + |R_\gamma^2f|^2 \, d\mu_\psi(\gamma) \Big)^\frac12 \Big\|_{L_p(\Gamma)} \, \sim_{c(p)} \, \|f\|_{L_p(\Gamma)}.$$ Now, let us go back to Lust-Piquard's setting $\Gamma = \Gamma_0 \times \Gamma_0 \times \cdots \times \Gamma_0$ and let us write $\sigma_j = (0, \ldots, 0, s_0, 0, \ldots, 0)$ with $s_0$ in the $j$-th entry. If we pick our measure $\mu_\psi$ to be the sum of point-masses $\mu_\psi = \sum_j \delta_{\sigma_j}$ and recall the following simple identities  
\begin{eqnarray*}
R_j & = & \partial_j \mathcal{L}^{-\frac12} \ = \ R_{\sigma_j}^1 - i R_{\sigma_j}^2, \\ R_j^* & = & \mathcal{L}^{-\frac12} \partial_j^* \ = \ R_{\sigma_j}^1 + i R_{\sigma_j}^2,
\end{eqnarray*}
then we deduce Lust-Piquard's theorem as a particular case of Theorem A1 for LCA groups. We may also recover her inequalities  for $1 < p < 2$ ---which has a more intricate statement--- from our result, we leave the details to the reader. The advantage of our formulation is that it is much more flexible. For instance, we do not need to impose $s_0$ to be torsion free and take $\Gamma_0 = \Z_k$ for any $k \ge 2$, which was left open in \cite{LuP}. We may also consider other point-masses, giving rise to other forms of discrete laplacians. In particular, we may allow the shift $s_0$ to depend on $j$ or even the group $\Gamma$ not to be given in a direct product form. Many other (not necessarily finitely supported) measures can also be considered. We will not make an exhaustive analysis of this here.

\vskip3pt

\noindent \textbf{C. Word-length laplacians.} We may work with many other discrete groups equipped with more or less standard lengths. Let us consider one of the most canonical examples, the word length. Given a finitely generated discrete group $\G$, the word length $|g|$ is the distance from $g$ to $e$ in the Cayley graph of $\G$. It is not always conditionally negative, but when this is the case we may represent it via the cocycle $(\H_{| \ |}, \alpha_{| \ |}, b_{| \ |})$, where $\H_{| \ |}$ is the closure of the pre-Hilbert space defined in $\Pi_0$ with $$\langle f_1, f_2 \rangle_{| \ |} = - \frac12 \sum_{g,h \in \G} \overline{\widehat{f}_1(g)} \widehat{f}_2(h) |g^{-1}h|$$ and $(\alpha_{| \ |,g}(f), b_{| \ |}(g)) = (\lambda(g)f, \lambda(e) - \lambda(g))$, as usual. The Riesz transforms $$R_{| \ |, j}f \, = \, \partial_{| \ |,j} A_{| \ |}^{-\frac12}f \, = \, 2 \pi i \sum_{g \in \G} \frac{\langle b_{| \ |}(g), e_j \rangle_{\H_{| \ |}}}{\sqrt{|g|}} \widehat{f}(g) \lambda(g)$$ satisfy Theorem A1 for any orthonormal basis $(e_j)_{j \ge 1}$ of $\H_{| \ |}$. Theorem A2 also applies with $A_{| \ |}(\lambda(g)) = |g| \lambda(g)$ and the associated gradient form $\Gamma_{| \ |}$. Further less standard Riesz transforms arise from other positive linear functionals $\tau_\psi$ on $\Pi_0$. Let us recall a few well-known groups for which the associated word length is conditionally negative:

\begin{itemize}
\item[i)] \textbf{Free groups.} The conditional negativity of the length function for the free group $\mathbb{F}_n$ with $n$ generators was established by Haagerup \cite{H}. A concrete description of the associated cocycle yields an interesting application of Theorem A1. The Gromov form is $$\hskip30pt K(g,h) = \frac{|g| + |h| - |g^{-1}h|}{2} = |\min(g,h)|$$ where $\min(g,h)$ is the longest word inside the common branch of $g$ and $h$ in the Cayley graph. If $g$ and $h$ do not share a branch, we let $\min(g,h)$ be the empty word. Let $\R[\mathbb{F}_n]$ stand for the algebra of finitely supported real valued functions on $\mathbb{F}_n$ and consider the bracket $\langle \delta_g, \delta_h \rangle = K(g,h)$. If $N$ denotes its null space, let $\H$ be the completion of $\R[\mathbb{F}_n]/N$. Then, writing $g^{-}$ for the word which results after deleting the last generator on the right of $g$, the system $\xi_g = \delta_g - \delta_{g^{-}} + N$ for all $g \in \mathbb{F}_n \setminus \{e\}$ is orthonormal and generates $\H$. Indeed, it is obvious that $\delta_e$ belongs to $N$ and we may write $$\hskip32pt \phi = \sum_{g \in \mathbb{F}_n} a_g \delta_g = \sum_{g \in \mathbb{F}_n} \Big( \sum_{h \ge g} a_h \Big) \xi_g$$ for any $\phi \in \R[\mathbb{F}_n]$. Here $h \ge g$ iff $g$ belongs to the (unique) path from $e$ to $h$ in the Cayley graph and we set $\xi_e = \delta_e$. This shows that $N = \R \delta_e$ and $\dim \H = \infty$. It yields a cocycle with $\alpha = \lambda$ and $$\hskip32pt b: g \in \mathbb{F}_n \mapsto \delta_g + \R \delta_e \in \R[\mathbb{F}_n] / \R \delta_e.$$ Considering the ONB $(\xi_h)_{h \neq e}$ we see that $\langle b(g), \xi_h \rangle = \delta_{g \ge h}$ and $$\hskip32pt R_{| \ |,h}f \, = \, \sum_{g \ge h} \frac{1}{\sqrt{|g|}} \widehat{f}(g) \lambda(g)$$ form a natural family of Riesz transforms in $\mathcal{L}(\mathbb{F}_n)$ with respect to the word length. If $p \ge 2$, Theorem A1 gives rise to the inequalities below in the free group algebra with constants $c(p)$ only depending on $p$ $$\hskip32pt \|f\|_p \, \sim_{c(p)} \, \Big\| \Big( \sum_{h \neq e} \Big| \sum_{g \ge h} \frac{1}{\sqrt{|g|}} \widehat{f}(g) \lambda(g) \Big|^2 + \Big| \sum_{g \ge h} \frac{1}{\sqrt{|g|}} \overline{\widehat{f}(g^{-1})} \lambda(g) \Big|^2 \Big)^\frac12 \Big\|_p.$$ The case $1 < p < 2$ can be formulated similarly according to Theorem A1. Note that the $L_p$-boundedness of the Riesz transforms $R_{| \ |,h}$ is out of the scope of Haagerup's inequality \cite{H}, since they are unbounded on $L_\infty$.

\vskip5pt

\item[ii)] \textbf{Finite cyclic groups.} The group $\Z_n$ of $n$-th roots of unity in the unit circle provides a central example to study. Despite its simplicity it may be difficult to provide precise bounds for Fourier multipliers, see for instance \cite{JPPP} for a recent discussion on optimal hypercontractivity bounds. The conditional negativity for the word length $|k| = \min(k,n-k)$ in $\Z_n$ was justified in that paper. For simplicity we will assume that $n = 2m$ is an even integer. Its Gromov form $$\hskip36pt K(j,k) \, = \, \frac{(j \wedge (2m-j)) + (k \wedge (2m-k)) - ((k-j) \wedge (2m-k+j))}{2}$$ can be written as follows for $0 \le j \le k \le 2m-1$ $$\hskip34pt K(j,k) \, = \, \langle \delta_j, \delta_k \rangle \, = \, j \wedge (2m-k) \wedge (m-k+j)_+$$ with $s_+ = 0 \vee s$, details are left to the reader. Using this formula, we may consider the associated bracket in $\Pi_0 \simeq \R[\Z_{2m}] \ominus \R\delta_0$ and deduce the expression below after rearrangement $$\hskip36pt \Big\langle \sum_{j=1}^{2m-1} a_j \delta_j, \sum_{j=1}^{2m-1} b_j \delta_j \Big\rangle \, = \, \sum_{k=1}^m \Big( \sum_{j=k}^{k+m-1} a_j \Big) \Big( \sum_{j=k}^{k+m-1} b_j \Big) \, = \, \sum_{k=1}^m \alpha_k \beta_k.$$ This shows that the null space $N$ of the bracket above is the space of elements $\sum_j a_j \delta_j$ in $\R[\Z_{2m}] \ominus \R\delta_0$ with vanishing coordinates $\alpha_k$. If we quotient out this subspace we end up with a Hilbert space $\H$ of dimension $(2m-1) - (m-1) = m$. Our discussion establishes an isometric isomorphism $$\hskip32pt \Phi: \sum_{j=1}^{2m-1} a_j \delta_j \in \H \mapsto \sum_{k=1}^m \alpha_k e_k \in \R^m.$$ Consider the ONB of $\H$ given by $\phi_k = \Phi^{-1}(e_k)$ for $1 \le k \le m$. We are ready to construct explicit Riesz transforms associated to the word length in $\Z_{2m}$. Namely, since $\langle \delta_j, \phi_k \rangle = \Phi(\delta_j)_k = 1$ iff $$\hskip36pt j \, \in \, \Lambda_k = \Big\{ j \in \Z_{2m} \, \big| \, j-k \equiv s \ \mathrm{mod} \ 2m \ \mbox{ for some } \ 0 \le s \le m-1 \Big\},$$ and vanishes otherwise. Then we may consider $$\hskip32pt R_{| \ |, k}f = \sum_{j \in \Lambda_k} \frac{1}{\sqrt{j \wedge (2m-j)}} \widehat{f}(j) e^{2\pi i \frac{j}{2m} \cdot}.$$ If $p \ge 2$, Theorem A1 establishes the following equivalence in $L_p(\Z_{2m})$ $$\hskip36pt \Big\| \sum_{j \in \Z_{2m}} \widehat{f}(j) e^{2\pi i \frac{j}{2m} \cdot} \Big\|_p \, \sim_{c(p)} \, \Big\| \Big( \sum_{k \in \Z_{2m}} \Big| \sum_{j \in \Lambda_k} \frac{\widehat{f}(j)}{\sqrt{j \wedge (2m-j)}} e^{2 \pi i \frac{j}{2m} \cdot} \Big|^2 \Big)^\frac12 \Big\|_p$$ with constants independent of $m$. Similar computations can be performed for $n \hskip-1pt = \hskip-1pt 2m+1$ odd and Theorem A1 can be easily reformulated in this case.

\vskip5pt

\item[iii)] \textbf{Infinite Coxeter groups.} Any group presented by $$\hskip32pt \G \, = \, \Big\langle g_1, g_2, \ldots, g_m \, \big| \, (g_jg_k)^{s_{jk}}= e \Big\rangle$$ with $s_{jj}=1$ and $s_{jk} \ge 2$ for $j \neq k$ is called a Coxeter group. Bo\.zejko proved in \cite{Bz} that the word length is conditionally negative for any infinite Coxeter group. The Cayley graph of these groups is more involved and we will not construct here a natural ONB for the associated cocycle, we invite the reader to do it. Other interesting examples include conditionally negative lengths in discrete Heisenberg groups or symmetric groups.    
\end{itemize} 

\section{{\bf H\"ormander-Mihlin multipliers}}

In this section we shall prove Theorems B1 and B2. Before that we obtain some preliminary Littlewood-Paley type inequalities. Afterwards we shall also study a few Besov limiting conditions in the spirit of Seeger for group von Neumann algebras that follow from Riesz transforms estimates.

\subsection{Littlewood-Paley estimates} \label{LPEstimates}

Let $(\G,\psi)$ be a discrete group equipped with a conditionally negative length and associated cocycle $(\H_\psi, \alpha_\psi, b_\psi)$. If $h \in \H_\psi$ we shall write $R_{\psi,h}$ for the $h$-directional Riesz transform $$R_{\psi,h}f \, = \, 2\pi i \sum_{g \in \G} \frac{\langle b_\psi(g), h \rangle_{\H_\psi}}{\sqrt{\psi(g)}} \widehat{f}(g) \lambda(g).$$ Let us recollect some useful square function inequalities which will be used below.  

\begin{lemma} \label{Littlewood-Paley}
Given $(h_j)_{j \ge 1}$ in $\H_\psi$ and $1 < p < \infty$ 
$$\hskip10pt \Big\| \Big( \sum_{j \ge 1} |R_{\psi,h_j} f_j|^2 \Big)^\frac12 \Big\|_{L_p(\widehat{\mathbf{G}})} \, \lesssim_{c(p)} \, \Big( \sup_{j \ge 1} \|h_j\|_{\H_\psi} \Big) \Big\| \Big( \sum_{j \ge 1} |f_j|^2 \Big)^\frac12 \Big\|_{L_p(\widehat{\mathbf{G}})}.$$
\end{lemma}

\dem Given an ONB $(e_k)_{k \ge 1}$ of $\H_\psi$, we have
\begin{eqnarray*}
\lefteqn{\Big\| \Big( \sum_{j \ge 1} |R_{\psi,h_j} f_j|^2 \Big)^\frac12 \Big\|_{L_p(\widehat{\mathbf{G}})}} \\ & = & \Big\| \sum_{j \ge 1} \Big( \sum_{g \in \G} \sum_{k \ge 1} \frac{\langle h_j, e_k \rangle_{\H_\psi} \langle b_\psi(g), e_k \rangle_{\H_\psi}}{\sqrt{\psi(g)}}\widehat{f}_j(g) \lambda(g) \Big) \otimes e_{j,1} \Big\|_{S_p(L_p)} \\ & = & \Big\| \sum_{j,k \ge 1} R_{\psi,k}f_j \otimes \langle h_j, e_k \rangle_{\H_\psi} e_{j,1} \Big\|_{S_p(L_p)} = \ \Big\| \sum_{j,k \ge 1} R_{\psi,k}f_j \otimes \Lambda(e_{(j,k),1}) \Big\|_{S_p(L_p)}, 
\end{eqnarray*}
where $\Lambda: \delta_{(j,k)} \in \ell_2(\N \times \N) \mapsto \langle h_j, e_k \rangle_{\H_\psi} \delta_j \in \ell_2(\N)$. Since $\ell_2^c$ is an homogeneous Hilbertian operator space \cite{P3}, the cb-norm of $\Lambda$ coincides with its norm which in turn equals $\sup_j \|h_j\|_{\H_\psi}$. Altogether, we deduce 
\begin{eqnarray*}
\Big\| \Big( \sum_{j \ge 1} |R_{\psi,h_j} f_j|^2 \Big)^\frac12 \Big\|_{L_p(\widehat{\mathbf{G}})} & \le & \Big( \sup_{j \ge 1} \|h_j\|_{\H_\psi} \Big) \Big\| \sum_{j,k \ge 1} R_{\psi,k}f_j \otimes e_{(j,k),1} \Big\|_{S_p(L_p)} \\ & = & \Big( \sup_{j \ge 1} \|h_j\|_{\H_\psi} \Big) \Big\| \Big( \sum_{j,k \ge 1} |R_{\psi,k}f_j|^2 \Big)^\frac12 \Big\|_{L_p(\widehat{\mathbf{G}})} \\ & = & \Big( \sup_{j \ge 1} \|h_j\|_{\H_\psi} \Big) \Big\| \Big( \sum_{k \ge 1} |\widetilde{R}_{\psi,k}f|^2 \Big)^\frac12 \Big\|_{S_p(L_p)} 
\end{eqnarray*}
for $\widetilde{R}_{\psi,k} = R_{\psi,k} \otimes id_{\mathcal{B}(\ell_2)}$ and $f = \sum_j f_j \otimes e_{j,1}$. Now, since Theorem A1 also holds in the category of operator spaces, the last term on the right hand side is dominated by $c(p) \|f\|_{S_p(L_p)}$, which yields the inequality we are looking for. \fin

We need to fix some standard terminology for our next result. Let us consider a sequence of functions $\varphi_j: \R_+ \to \C$ in $\mathcal{C}^{\mathrm{k}_n}(\R_+ \setminus \{0\})$ for $\mathrm{k}_n = [\frac{n}{2}] +1$ such that the following inequalities hold $$\Big( \summ_j \Big| \frac{d}{d\xi^k} \varphi_j(\xi) \Big|^2 \Big)^\frac12 \, \lesssim \, |\xi|^{-k} \quad \mbox{for all} \quad 0 \le k \le \mathrm{k}_n.$$ Let $\mathsf{M}_p(\widehat{\mathbf{G}})$ stand for the space of symbols $m: \G \to \C$ associated to $L_p$-bounded Fourier multipliers in the group von Neumann algebra $\V$, equip any such symbol $m$ with the $p \to p$ norm of the multipler $\lambda(g) \mapsto m(g) \lambda(g)$. Now, given any conditionally negative length $\psi: \G \to \R_+$ and $h$ in the associated cocycle Hilbert space $\H_\psi$, let us consider the symbols $$m_{\psi,h}(g) = \frac{\langle b_\psi(g), h \rangle_\psi}{\sqrt{\psi(g)}}.$$ Then, we may combine families of these symbols in a single Fourier multiplier patching them via the Littlewood-Paley decompositions provided by the families $(\varphi_j)$ and finite-dimensional cocycles on $\G$. The result is the following. 

\begin{lemma} \label{Littlewood-Paley2}
Let $\G$ be a discrete group equipped with two $n$-dimensional cocycles with associated length functions $\psi_1, \psi_2$ and an arbitrary cocycle with associated length function $\psi_3$. Let $(\varphi_{1j})$, $(\varphi_{2j})$ be Littlewood-Paley decompositions satisfying the assumptions above and $(h_j)_{j \ge 1}$ in $\H_{\psi_3}$. Then, the following inequality holds for any $1 < p < \infty$ $$\Big\| \sum_{j \ge 1} \varphi_{1j}(\psi_1(\cdot)) \varphi_{2j}(\psi_2(\cdot)) m_{\psi_3,h_j} \Big\|_{\mathsf{M}_p(\widehat{\mathbf{G}})} \, \lesssim_{c(p,n)} \, \sup_{j \ge 1} \|h_j\|_{\H_{\psi_3}}.$$
\end{lemma}

\dem If $\Lambda_{\psi, \varphi}: \lambda(g) \mapsto \varphi(\psi(g)) \lambda(g)$, we have
\begin{eqnarray*}
\lefteqn{\Big| \tau \Big( \widetilde{f}^* \, \summ_j \Lambda_{\psi_1,\varphi_{1j}} \Lambda_{\psi_2,\varphi_{2j}} R_{\psi_3,h_j} (f) \Big) \Big|} \\ [3pt] & = & \Big| \summ_j \tau \Big( \big( \Lambda_{\psi_1, \varphi_{1j}} \widetilde{f} \, \big)^* \, R_{\psi_3,h_j} \big( \Lambda_{\psi_2, \varphi_{2j}} f \big) \Big) \Big| \\ [3pt] & \le & \Big\| \summ_j \Lambda_{\psi_1, \varphi_{1j}} \widetilde{f} \otimes \delta_j \Big\|_{RC_{p'}(\V)} \Big\| \summ_j R_{\psi_3,h_j} \big( \Lambda_{\psi_2, \varphi_{2j}} f \big) \otimes \delta_j \Big\|_{RC_p(\V)} 
\end{eqnarray*}
by anti-linear duality. On the other hand, observe that Lemma \ref{Littlewood-Paley} trivially extends to $RC_p$-spaces. Indeed, the row inequality follows from the column one applied to $f_j^*$s together with Remark \ref{RemRowColumnRiesz} and the fact that the maps $R_{\psi,h_j}'$ satisfy the same estimates. \hskip-1pt Combining \hskip-1pt that \hskip-1pt with the Littlewood-Paley inequality in \cite[Theorem C]{JMP0} we obtain $$\displaystyle \Big| \tau \Big( \widetilde{f}^* \, \sum_{j \ge 1}^{\null} \Lambda_{\psi_1, \varphi_{1j}} \Lambda_{\psi_2, \varphi_{2j}} R_{\psi_3,h_j}(f) \Big) \Big|  \lesssim_{c(p,n)} \Big( \sup_{j \ge 1} \|h_j\|_{\H_{\psi_3}} \Big) \|f\|_{L_p(\widehat{\mathbf{G}})} \|\widetilde{f}\|_{L_{p'}(\widehat{\mathbf{G}})}.$$ The result follows taking supremums of $\widetilde{f}$ running over the unit ball of $L_{p'}(\widetilde{\mathbf{G}})$. \fin

The only drawback of Lemma \ref{Littlewood-Paley2} is that we need finite-dimensional cocycles to apply our Littlewood-Paley estimates from \cite{JMP0}. We may ignore that requirement at the price of using other square function inequalities from \cite{JLX},  where the former $\varphi_j$ are now dilations of a fixed function. Namely, given $1 < p < \infty$ and $\varphi: \R_+ \to \C$ belonging to a certain class $\mathcal{J}_p$ of analytic functions, it turns our that the column Hardy norm $$\|f\|_{H_p^c(\varphi,\psi)} \, = \, \Big\| \Big( \int_{\R_+} \big| \varphi(s A_\psi) f\big|^2 \frac{ds}{s} \Big)^\frac12 \Big\|_{L_p(\widehat{\mathbf{G}})}$$ does not depend on the chosen function $\varphi \in \mathcal{J}_p$. We know from \cite{JLX} that $$L_p(\widehat{\mathbf{G}}) \, \subset\, \begin{cases} H_p^r(\varphi,\psi) + H_p^c(\varphi,\psi) & \mbox{if } 1 < p \le 2, \\ H_p^r(\varphi,\psi) \cap H_p^c(\varphi,\psi) & \mbox{if } 2 \le p < \infty, \end{cases}$$ for any $\varphi \in \displaystyle  \mathcal{J} = \bigcap_{1 < p < \infty} \mathcal{J}_p \neq \emptyset$ and any conditionally negative length $\psi: \G \to \R_+$.

\begin{lemma} \label{Littlewood-Paley3}
Let $\G$ be a discrete group equipped with three arbitrary cocycles with associated length functions $\psi_1, \psi_2, \psi_3$. Let $\varphi_1, \varphi_2 \in \mathcal{J}$ and $(h_s)_{s > 0}$ in $\H_{\psi_3}$. Then the following inequality holds for any $1 < p < \infty$ $$\Big\| \int_{\R_+} \varphi_1(s \psi_1(\cdot)) \varphi_2(s \psi_2(\cdot)) m_{\psi_3,h_s} \frac{ds}{s} \Big\|_{\mathsf{M}_p(\widehat{\mathbf{G}})} \, \lesssim_{c(p)} \, \esssup_{s > 0} \|h_s\|_{\H_{\psi_3}}.$$
\end{lemma}

\dem Let us write $L_p(\widehat{\mathbf{G}}; L_{rc}^2(\R_+))$ for the space $RC_p(\V)$ in which we replace discrete sums over $\Z_+$ by integrals on $\R_+$ with the harmonic measure. Arguing as in Lemma \ref{Littlewood-Paley}, it is not difficult to show that $$\Big\| \int_{\R_+}^\oplus R_{\psi_3,h_s} \frac{ds}{s} \, : L_p(\widehat{\mathbf{G}}; L_{rc}^2(\R_+)) \to L_p(\widehat{\mathbf{G}}; L_{rc}^2(\R_+)) \Big\| \, \lesssim_{c(p)} \, \esssup_{s > 0} \|h_s\|_{\H_{\psi_3}}.$$ According to the results from \cite{JLX}, the maps 
\begin{eqnarray*}
\Phi_1: \hskip2.5pt f \in L_p(\widehat{\mathbf{G}}) & \mapsto & \big( \varphi_1( \cdot A_{\psi_1}) f \big)_{s > 0} \in L_p(\widehat{\mathbf{G}}; L_{rc}^2(\R_+)), \\ [5pt] \Phi_2: f \in L_{p'}(\widehat{\mathbf{G}}) & \mapsto & \big( \varphi_2( \cdot A_{\psi_2}) f \big)_{s > 0} \in L_{p'}(\widehat{\mathbf{G}}; L_{rc}^2(\R_+)),
\end{eqnarray*} 
are bounded. Therefore, the result follows by noticing that $$\displaystyle \int_{\R_+} \varphi_1(s A_{\psi_1}) \varphi_2(s A_{\psi_2}) R_{\psi_3,h_j} \frac{ds}{s} \, = \, \Phi_2^* \circ \Big( \int_{\R_+}^\oplus R_{\psi_3,h_s} \frac{ds}{s} \Big) \circ \Phi_1,$$ and we get the expected inequality from the estimates given above. \fin 

\vskip8pt

\begin{remark} \label{MoreLPRem}
\emph{Slight modifications also give:}
\begin{itemize}
\item[i)] \emph{If $p \ge 2$} $$\hskip10pt \Big\| \Big( \sum_{j \ge 1} |R_{\psi,h_j} f|^2 \Big)^\frac12 \Big\|_{L_p(\widehat{\mathbf{G}})} \ \lesssim_{c(p)} \ \Big\| \Big( \langle h_j, h_k \rangle \Big)^\frac12 \Big\|_{\mathcal{B}(\ell_2)} \|f\|_{L_p(\widehat{\mathbf{G}})}.$$

\vskip3pt

\item[ii)] \emph{If $1 < p < \infty$} $$\hskip10pt \Big\| \sum_{j \ge 1} \Lambda_{\psi_1,\varphi_j} R_{\psi_3,h_j} f \otimes \delta_j \Big\|_{RC_p(\V)} \ \lesssim_{c(p,n)} \ \Big( \sup_{j \ge 1} \|h_j\|_{\H_{\psi_3}} \Big) \|f\|_{L_p(\widehat{\mathbf{G}})}.$$ \emph{It becomes an equivalence when $\sum_j |\varphi_j|^2 \equiv 1$ and $(h_j)$ is an ONB of $\H_{\psi_2}$.}

\vskip5pt

\item[iii)] \emph{If $1 < p < \infty$ and $0 < \theta < \frac{2}{p} \wedge \frac{2}{p'}$} $$\hskip12pt \Big\| \int_{\R_+} \varphi_1(s \psi_1(\cdot)) \varphi_2(s \psi_2(\cdot)) m_s \frac{ds}{s} \Big\|_{\mathsf{M}_p(\widehat{\mathbf{G}})} \lesssim_{c(p)} \esssup_{s > 0} \|m_s\|_{[\mathrm{X}_{\psi_3}, \ell_\infty(\G)]^\theta},$$ \emph{where $\sqrt{\psi_3(g)} m_s(g) = \langle h_s, b_{\psi_3}(g) \rangle_{\H_{\psi_3}}$ for all $g$ and $\|m_s\|_{\mathrm{X}_{\psi_3}} = \|h_s\|_{\H_{\psi_3}}$.}
\end{itemize}
\emph{Assertion i) follows as in Lemma \ref{Littlewood-Paley} with $\Lambda': \delta_k \in \ell_2(\N) \mapsto \sum_j \langle h_j, e_k \rangle_{\H_\psi} \delta_j \in \ell_2(\N)$ instead of the map $\Lambda$ used there. Moreover, by duality we also find the following inequality for $1 < p \le 2$ 
\begin{itemize}
\item[iv)] $\displaystyle \Big\| \sum_{j \ge 1} R_{\psi,h_j} f_j \Big\|_{L_p(\widehat{\mathbf{G}})} \, \lesssim_{c(p)} \, \Big\| \Big( \langle h_j, h_k \rangle \Big)^\frac12 \Big\|_{\mathcal{B}(\ell_2)} \Big\| \Big( \sum_{j \ge 1} |f_j|^2 \Big)^\frac12 \Big\|_{L_p(\widehat{\mathbf{G}})}.$ 
\end{itemize}
Estimate ii) follows clearly from the argument in Lemma \ref{Littlewood-Paley2}. According to Theorem A1 and \cite[Theorem C]{JMP0}, the equivalence holds when the sequence $|\varphi_j|^2$ forms a Littlewood-Paley partition of unity and the $h_j$ form an ONB of $\H_{\psi_3}$. Finally, iii) is just an improvement of Lemma \ref{Littlewood-Paley3} by interpolation. Namely, assuming (wlog) that $b_{\psi_3}(g)$ span $\H_{\psi_3}$ when $g$ runs over $\G$, it is clear that $h_s$ is univocally determined for $m_s$ and the norm in $\mathrm{X}_{\psi_3}$ is well-defined. Once this is settled, we know that iii) holds for $(p,\theta) = (q,0)$ with any $q < \infty$ and also for $(2,1)$. Interpolation of the maps $L_\infty(\R_+; \mathrm{X}_{\psi_3}) \to \mathsf{M}_q(\widehat{\mathbf{G}})$ and $L_\infty(\R_+; \ell_\infty(\G)) \to \mathsf{M}_2(\widehat{\mathbf{G}})$ with parameters $\frac{1}{p} = \frac{1-\theta}{q} + \frac{\theta}{2}$ yields the expected inequality. Note that $$[L_\infty(\R_+; \mathrm{X}_{\psi_3}), L_\infty(\R_+; \ell_\infty(\G))]^\theta = L_\infty \big( \R_+; [\mathrm{X}_{\psi_3},\ell_\infty(\G)]^\theta \big).$$}
\end{remark}

\subsection{A refined Sobolev condition}

In order to prove Theorem B1, we start with a basic inequality for Euclidean Fourier multipliers which is apparently new. Recall the definition of the differential operators $\mathsf{D}_\alpha$ and the fractional laplacian lengths $\psi_\beta$ from the Introduction. The proof shows how large is the family of Riesz transforms in $\R^n$ when we allow infinite-dimensional cocycles.

\begin{lemma} \label{HMLemma}
If $1 < p < \infty$ and $\varepsilon > 0$, then $$\|m\|_{\mathsf{M}_p(\R^n)} \, \lesssim_{c(p)} \, \Big\| \mathsf{D}_{\frac{n}{2} + \varepsilon} \big( \sqrt{\psi_\varepsilon} m \big) \Big\|_{L_2(\R^n)}.$$
\end{lemma}

\dem If we consider the Sobolev-type space $$\mathsf{W}_{(\frac{n}{2},\varepsilon)}^2(\R^n) = \Big\{ m: \R^n \to \C \ \big| \ m \ \mbox{measurable}, \ \Big\| \mathrm{D}_{\frac{n}{2} + \varepsilon} \big( \sqrt{\psi_\varepsilon} m \big) \Big\|_2 < \infty \Big\},$$ and $L_2^\circ(\R^n, \mu_\varepsilon)$ is the mean-zero subspace of $L_2(\R^n, \mu_\varepsilon)$, then we claim that $$\Lambda: h \in L_2^\circ(\R^n, \mu_\varepsilon) \mapsto m_h \in \mathsf{W}_{(\frac{n}{2},\varepsilon)}^2(\R^n),$$ $$m_h(\xi) = \frac{1}{\sqrt{\psi_\varepsilon(\xi)}} \int_{\R^n} h(x) \big( e^{2\pi i \langle \xi,x \rangle}-1 \big) \, d\mu_\varepsilon(x),$$ extends to a surjective isometry. Indeed, let $h$ be a Schwartz function in $L_2^\circ(\R^n, \mu_\varepsilon)$ with compact support away from $0$ and write $\omega_\varepsilon(x) = 1/|x|^{n+2\varepsilon}$ for the density $d\mu_\varepsilon(x)/dx$. In that case, the function $h \omega_\varepsilon$ is a mean-zero Schwartz function in $L_2(\R^n)$ and $$\sqrt{\psi_\varepsilon} m_h \, = \, \int_{\R^n} h(x) e^{2\pi i \langle \cdot ,x \rangle} \, d\mu_\varepsilon(x) \, = \, \big( h \omega_\varepsilon \big)^\vee.$$ In particular, this implies that we have 
\begin{eqnarray*}
\big\| \Lambda(h) \big\|_{\mathsf{W}_{(\frac{n}{2},\varepsilon)}^2(\R^n)} & = & \Big\| \mathsf{D}_{\frac{n}{2} + \varepsilon} \big( \sqrt{\psi_\varepsilon} m_h \big) \Big\|_{L_2(\R^n)} \\ & = & \Big\| \frac{1}{\sqrt{\omega_\varepsilon}} \big( h \omega_\varepsilon \big)^{\vee \wedge} \Big\|_{L_2(\R^n)} \ = \ \|h\|_{L_2^\circ(\R^n,\mu_\varepsilon)}.
\end{eqnarray*}
Therefore, since smooth mean-zero compactly supported functions away from $0$ are clearly dense in $L_2^\circ(\R^n,\mu_\varepsilon)$, we may extend $\Lambda$ to an isometry. Now, to show subjectivity we observe from elementary facts on Sobolev spaces \cite{Ad} that the class of Schwartz functions is dense in our $\mathsf{W}$-space. Therefore, it suffices to show that $\Lambda^{-1}(m)$ exists for any such Schwartz function $m$. By Plancherel theorem, we find $$\Big\| \mathsf{D}_{\frac{n}{2} + \varepsilon} \big( \sqrt{\psi_{\varepsilon}} m \big) \Big\|_{L_2(\R^n)} = \Big\| \frac{1}{\sqrt{\omega_\varepsilon}} \widehat{\sqrt{\psi_\varepsilon} m} \Big\|_{L_2(\R^n)} = \Big\| \frac{1}{w_\varepsilon} \widehat{\sqrt{\psi_\varepsilon} m} \Big\|_{L_2(\R^n,\mu_\varepsilon)}.$$ Since $h = \frac{1}{\omega_\varepsilon} \widehat{\sqrt{\psi_\varepsilon} m}$ is also an Schwartz function, we may write 
\begin{eqnarray*}
\sqrt{\psi_\varepsilon(\xi)} m(\xi) & = & \int_{\R^n} h(x) \omega_\varepsilon(x) e^{2\pi i \langle \xi, x \rangle} dx \\ & = & \int_{\R^n} h(x) \big( e^{2\pi i \langle \xi, x \rangle} - 1 \big) d\mu_\varepsilon(x) + \int_{\R^n} h(x) \omega_\varepsilon(x) dx \ = \ \mathrm{A}(\xi) + \mathrm{B}. 
\end{eqnarray*}
The function $\mathrm{A}(\xi)$ is well-defined since both $h$ and $e^{2\pi i \langle \xi, x \rangle} - 1$ are elements of the Hilbert space $L_2(\R^n, \mu_\varepsilon)$, so its product is absolutely integrable. Moreover, since $h \omega_\varepsilon \in \mathcal{S}(\R^n)$ we see that $\mathrm{A}(0) = 0$ and conclude $\mathrm{B} = \sqrt{\psi_\varepsilon(0)} m(0) = 0$. This means that $h = \Lambda^{-1}(m)$ is a mean-zero function in $L_2(\R^n,\mu_\varepsilon)$, as desired. In other words we know that every $m \in \mathsf{W}_{(\frac{n}{2},\varepsilon)}^2(\R^n)$ satisfies $$m(\xi) \, = \, \frac{\big\langle h, e^{2\pi i \langle \xi, \cdot \rangle}-1 \big\rangle_{\mu_\varepsilon}}{\sqrt{\psi_\varepsilon(\xi)}} \quad \mbox{and} \quad \Big\| \mathsf{D}_{\frac{n}{2}+\varepsilon} \big( \sqrt{\psi_\varepsilon} m \big) \Big\|_2 = \|h\|_{L_2(\R^n,\mu_\varepsilon)}.$$ This shows that every $m$ in our Sobolev-type space is a $\psi_\varepsilon$-Riesz transform and its norm coincides with that of its symbol $h$. The only drawback is that we use a complex Hilbert space for our cocycle, so the inner product is not the right one. Consider the cocycle $(\H_\varepsilon, \alpha_\varepsilon, b_\varepsilon)$ with  
\begin{eqnarray*}
\H_\varepsilon & = & L_2(\R^n, \mu_\varepsilon; \R^2) \ \simeq \ L_2(\R^n, \mu_\varepsilon; \C), \\ [8pt] b_\varepsilon(\xi) & = & \Big( \cos(2\pi \langle \xi, \cdot \, \rangle) - 1, \, \sin(2 \pi \langle \xi, \cdot \, \rangle) \Big) \ \simeq \ e^{2 \pi i \langle \xi, \cdot \, \rangle} - 1, \\ [3pt] \alpha_{\varepsilon, \xi}(f) & = & \left( \begin{array}{lr} \cos(2\pi \langle \xi, \cdot \, \rangle) & - \sin(2\pi \langle \xi, \cdot \, \rangle) \\ \sin(2\pi \langle \xi, \cdot \, \rangle) & \cos(2\pi \langle \xi, \cdot \, \rangle) \end{array} \right) {{f_1}\choose{f_2}} \ \simeq \ e^{2 \pi i \langle \xi, \cdot \, \rangle} f.
\end{eqnarray*}
Then, it is easily checked that we have 
\begin{itemize}
\item If $h$ is $\R$-valued and odd $\displaystyle \big\langle h, e^{2\pi i \langle \xi, \cdot \rangle}-1 \big\rangle_{\mu_\varepsilon} = i \Big\langle {{0}\choose{h}}, b_\varepsilon(\xi) \Big\rangle_{\H_\varepsilon}$.

\vskip3pt

\item If $h$ is $\R$-valued and even $\displaystyle \big\langle h, e^{2\pi i \langle \xi, \cdot \rangle}-1 \big\rangle_{\mu_\varepsilon} = \Big\langle {{h}\choose{0}}, b_\varepsilon(\xi) \Big\rangle_{\H_\varepsilon}$.
\end{itemize}
Therefore, decomposing $$m = \Big( \mathrm{Re}(m_{\mathrm{odd}}) + \mathrm{Re}(m_{\mathrm{even}}) \Big) + i \Big( \mathrm{Im}(m_{\mathrm{odd}}) + \mathrm{Im}(m_{\mathrm{even}}) \Big)$$ and noticing the elementary inequalities $$\Big\| \mathrm{Re/Im}(m_{\mathrm{odd/even}}) \Big\|_{\mathsf{W}_{(\frac{n}{2}, \varepsilon)}^2(\R^n)} \, \le \, \|m \|_{\mathsf{W}_{(\frac{n}{2}, \varepsilon)}^2(\R^n)}$$ we see how every element in our Sobolev-type space decomposes as a sum of four Riesz transforms whose $\mathsf{M}_p$-norms are all dominated by $c(p) \|m\|_{\mathsf{W}_{(\frac{n}{2}, \varepsilon)}^2(\R^n)}$. \fin  

\begin{lemma} \label{HMLemma2}
Let $(\G,\psi)$ be a discrete group equipped with a conditionally negative length giving rise to an $n$-dimensional cocycle $(\H_\psi, \alpha_\psi, b_\psi)$. If $1 < p < \infty$ and $\varepsilon > 0$, then $$\|m\|_{\mathsf{M}_p(\widehat{\mathbf{G}})} \, \lesssim_{c(p)} \, \inf_{m = \widetilde{m} \circ b_\psi} \Big\| \mathsf{D}_{\frac{n}{2} + \varepsilon} \big( \sqrt{\psi_\varepsilon} \widetilde{m} \big) \Big\|_{L_2(\R^n)}.$$
\end{lemma}

\dem According to the proof of Lemma \ref{HMLemma}, $$m(g) = \sum_{j=1}^4 \frac{\big\langle h_j, b_\varepsilon \circ b_\psi(g) \big\rangle_{\H_\varepsilon}}{\sqrt{\psi_\varepsilon(b_\psi(g))}} = \frac{1}{\sqrt{\psi_\varepsilon(b_\psi(g))}} \Big\langle \sum_{j=1}^4 h_j, b_\varepsilon \circ b_\psi(g) \Big\rangle_{\H_\varepsilon} $$ for any $\widetilde{m} \in \mathsf{W}_{(n/2,\varepsilon)}^2(\R^n)$ satisfying $m = \widetilde{m} \circ b_\psi$ and certain $h_j \in L_2(\R^n, \mu_\varepsilon)$. Next we observe that each of the four summands above can still be regarded as Riesz transforms on $\V$ with respect to the following cocycle \vskip-10pt
\begin{eqnarray*}
\H_{\psi,\varepsilon} & = & L_2(\R^n, \mu_\varepsilon; \R^2), \\ [7pt] b_{\psi,\varepsilon}(g) & = & \Big( \cos(2\pi \langle b_\psi(g), \cdot \, \rangle) - 1, \, \sin(2 \pi \langle b_\psi(g), \cdot \, \rangle) \Big), \\ [3pt] \alpha_{\psi,\varepsilon,g}(f) & = & \left( \begin{array}{lr} \cos(2\pi \langle b_\psi(g), \cdot \, \rangle) & - \sin(2\pi \langle b_\psi(g), \cdot \, \rangle) \\ \sin(2\pi \langle b_\psi(g), \cdot \, \rangle) & \cos(2\pi \langle b_\psi(g), \cdot \, \rangle) \end{array} \right) {{f_1\circ \alpha_{\psi,g^{-1}}(\cdot)}\choose{f_2 \circ \alpha_{\psi,g^{-1}}(\cdot)}}.
\end{eqnarray*}

\vskip5pt

\noindent Then we conclude observing $\|m\|_{\mathsf{M}_p(\widehat{\mathbf{G}})} \le c(p) \|\sum_j h_j\|_{\H_\varepsilon} \lesssim c(p) \|\widetilde{m}\|_{\mathsf{W}_{(\frac{n}{2},\varepsilon)}^2(\R^n)}$. \fin

\demB Let $\G_{\psi,0} = \{g \in \G \, : \, \psi(g)=0 \}$ denote the subgroup of elements with vanishing $\psi$-length, which trivializes for injective cocycles. According to  $m = \widetilde{m} \circ b_\psi$, the multiplier $m$ is constant on $\G_{\psi,0}$ and takes the value $m(e)$. This means that the Fourier multiplier associated to $m_0 = m(e) \chi_{\G_{\psi,0}}$ is nothing but $m(e)$ times the conditional expectation onto $\mathcal{L}(\G_{\psi,0})$, so $$\|m\|_{\mathsf{M}_p(\widehat{\mathbf{G}})} \, \le \, |m(e)| +  \|m - m_0\|_{\mathsf{M}_p(\widehat{\mathbf{G}})}.$$ Since $m-m_0 = (\widetilde{m} - m(e) \delta_0) \circ b_\psi$ and $\widetilde{m} = \widetilde{m} - m(e) \delta_0$ almost everywhere, we may just proceed by assuming that $m(g)=0$ for all $g \in \G_{\psi,0}$. Let $\eta$ be a radially decreasing smooth function on $\R$ with $\chi_{(-1,1)} \le \eta \le \chi_{(-2,2)}$ and set $\phi(\xi) = \eta(\xi) - \eta(2\xi)$, so that me may construct a standard Littlewood-Paley partition of unity $\phi_j(\xi) = \phi(2^{-j} \xi)$ for $j \in \Z$. Note that we have $$\sum_{j \in \Z} \phi_j(\xi) = \begin{cases} 1 & \mbox{if} \ \xi \neq 0, \\ 0 & \mbox{if} \ \xi = 0. \end{cases}$$ Using three similar constructions with suitably chosen generating functions yields another partition of unity $(\rho_j)_{j \in \Z}$ with $\rho_j \equiv 1/3$ on the support of $\phi_j$. We shall be working with the radial Littlewood-Paley partition of unity in $\R^n$ given by $\varphi_j = \rho_j \circ | \ |^2$. On the other hand, if we set $\varphi_{1j} = \varphi_{2j} = \sqrt{\phi_j}$ we may write $m$ as follows (recall that $m$ vanishes where $\psi$ does)
\begin{eqnarray*}
m & = & \sum_{j \in \Z} (\phi_j \circ \psi) \, m \ = \ 3 \sum_{j \in \Z} (\phi_j \circ \psi) \, (\rho_j \circ \psi) \, m \\ & = & 3 \sum_{j \in \Z} (\varphi_{1j} \circ \psi) \, (\varphi_{2j} \circ \psi) (\rho_j \circ \psi) \, m \ = \ 3 \sum_{j \in \Z} (\varphi_{1j} \circ \psi) \, (\varphi_{2j} \circ \psi) \, m_j.
\end{eqnarray*}
Since we have $$m_j(g) = \rho_j(\psi(g)) m(g) = \rho_j(|b_\psi(g)|^2) \widetilde{m}(b_\psi(g)) = (\varphi_j \widetilde{m}) (b_\psi(g)),$$ we deduce that $m_j = \widetilde{m}_j \circ b_\psi$ with $\widetilde{m}_j = \varphi_j \widetilde{m}$. We know by assumption that $$\sup_{j \in \Z} \Big\| \mathsf{D}_{\frac{n}{2} + \varepsilon} \big( \sqrt{\psi_\varepsilon} \, \varphi_j \, \widetilde{m} \big) \Big\|_2 \, < \, \infty.$$ According to the proof Lemma \ref{HMLemma2}, this implies that we may write $m_j$ as a Riesz transform $m_{\gamma,h_j}$ with respect to the length function $\gamma = \psi_\varepsilon \circ b_\psi$, whose associated infinite-dimensional cocycle was also described in the proof of Lemma \ref{HMLemma2}. Since the families $\varphi_{1j}, \varphi_{2j}$ satisfy the assumptions of Lemma \ref{Littlewood-Paley2}, we may combine this result with the proof of Lemma \ref{HMLemma} to obtain
\begin{eqnarray*}
\|m\|_{\mathsf{M}_p(\widehat{\mathbf{G}})} & = & \Big\| \sum_{j \in \Z} \varphi_{1j}(\psi_1(\cdot)) \varphi_{2j}(\psi_2(\cdot)) m_{\gamma,h_j} \Big\|_{\mathsf{M}_p(\widehat{\mathbf{G}})} \\ & \le & c(p,n) \, \sup_{j \in \Z} \|h_j\|_{\H_{\gamma}} \ \le \ c(p,n) \, \sup_{j \in \Z} \Big\| \mathsf{D}_{\frac{n}{2} + \varepsilon} \big( \sqrt{\psi_\varepsilon} \, \varphi_j \, \widetilde{m} \big) \Big\|_2. 
\end{eqnarray*}
The dependence on $n$ of $c(p,n)$ comes from the Littlewood-Paley inequalities. \fin

\noindent In the following result, we use the standard notation $$\|f\|_{\stackrel{\hskip-5pt \circ}{\mathsf{H}_\alpha^2}(\Omega)} = \big\| (| \ |^{\alpha} \widehat{f} \, )^{\vee} \big\|_{L_2(\Omega)} \quad \mbox{and} \quad \|f\|_{\mathsf{H}_\alpha^2(\Omega)} = \big\| ((1 + | \ |^2)^{\frac{\alpha}{2}} \widehat{f} \, )^{\vee} \big\|_{L_2(\Omega)}.$$

\begin{corollary} \label{SobVsClassical}
If $1 < p < \infty$ and $0 < \varepsilon < [\frac{n}{2}] + 1 - \frac{n}{2}$, we also have $$\|m\|_{\mathsf{M}_p(\widehat{\mathbf{G}})} \, \lesssim_{c(p,n,\varepsilon)} |m(e)| + \inf_{m = {\widetilde{m}}^{\null} \circ b_\psi} \left\{ \sup_{j \in \Z} \big\| \varphi_0 \, \widetilde{m}(2^j \cdot ) \big\|_{\mathsf{H}_{\frac{n}{2}+\varepsilon}^2(\R^n)} \right\}$$ for any pair $(\G,\psi)$ which gives rise to an $n$-dimensional cocycle $(\H_\psi, \alpha_\psi, b_\psi)$. 
\end{corollary}

\dem We have 
\begin{eqnarray*}
\Big\| \mathsf{D}_{\frac{n}{2} + \varepsilon} \Big( \sqrt{\psi_\varepsilon} \varphi_j \widetilde{m} \Big) \Big\|_2 & = & \sqrt{\mathrm{k}_n(\varepsilon)} \, \Big\| \mathsf{D}_{\frac{n}{2} + \varepsilon} \big( | \ |^\varepsilon \varphi_j \widetilde{m} \big) \Big\|_2 \\ & = & \sqrt{\mathrm{k}_n(\varepsilon)} \, \Big\| \mathsf{D}_{\frac{n}{2} + \varepsilon} \big( | \ |^\varepsilon \varphi_0 \widetilde{m}(2^j \cdot) \big) \Big\|_2
\end{eqnarray*}
with $\mathrm{k}_n(\varepsilon) \sim \frac{\pi^{n/2}}{\Gamma(n/2)} \frac{1}{\varepsilon}$. Indeed, it is simple to check that $\mathsf{W}_{(n/2,\varepsilon)}^2(\R^n)$ has a dilation invariant norm. It therefore suffices to show that the inequality below holds up to a constant $c(n,\varepsilon)$ $$\Big\| \mathsf{D}_{\frac{n}{2} + \varepsilon} \big( | \ |^{\varepsilon} f \big) \Big\|_2 \, \lesssim \, \Big\| \big( 1 + | \ |^2 \big)^{\frac12(\frac{n}{2}+\varepsilon)} \widehat{f} \, \Big\|_2$$ for functions supported by (say) the corona $\Omega = \mathrm{B}_2(0) \setminus \mathrm{B}_1(0)$, which is the form of the support of $\varphi_0$. In other words, we need to show that $$\Lambda_\varepsilon: f \in \mathsf{H}_{\frac{n}{2} + \varepsilon}^2(\Omega) \mapsto | \ |^\varepsilon f \in \ \stackrel{\hskip-18pt \circ}{\mathsf{H}_{\frac{n}{2}+\varepsilon}^2}(\Omega)$$ defines a bounded operator. These two families of spaces satisfy the expected interpolation identities in the variable $\alpha = (1-\theta) \alpha_0 + \theta \alpha_1$ with respect to the complex method. To prove that $\Lambda_\varepsilon$ is bounded we use Stein's interpolation. Assume by homogeneity that $f$ is in the unit ball of $\mathsf{H}_{\frac{n}{2}+\varepsilon}^2(\Omega)$ and let $\theta = 2\varepsilon \delta_{n \, \mathrm{odd}} + \varepsilon \delta_{n \, \mathrm{even}}$. Thus, there exists $F$ analytic on the strip satisfying $F(\theta) = f$ and $$\max \Big\{ \sup_{t \in \R} \|F(it)\|_{\mathsf{H}_{\frac{n}{2}}^2(\Omega)}, \sup_{t \in \R} \|F(1+it)\|_{\mathsf{H}_{[\frac{n}{2}]+1}^2(\Omega)} \Big\} \, \le \, 1.$$ Given $\delta > 0$, define now the following analytic family of operators $$\mathcal{L}_z(F) = \exp \big( \delta (z - \theta)^2 \big) | \ |^{\varepsilon z/\theta} F(z),$$ so that $\mathcal{L}_\theta(F) = \Lambda_\varepsilon(f)$. We claim that the following estimates hold
\begin{eqnarray*}
\big\| \mathcal{L}_{it}(F) \big\|_{\stackrel{\hskip-7pt \circ}{\mathsf{H}_{\frac{n}{2}}^2}(\Omega)} & \le & c(n) \big( 1 + |t| \big)^{\frac{n}{2}} \, e^{-\delta t^2} e^{\delta \theta^2}, \\ \big\| \mathcal{L}_{1+it}(F) \big\|_{\stackrel{\hskip-19pt \circ}{\mathsf{H}_{[\frac{n}{2}]+1}^2}(\Omega)} & \le & c(n) \big(1 + |t| \big)^{[\frac{n}{2}]+1} e^{-\delta t^2} e^{\delta (1-\theta)^2}. 
\end{eqnarray*}
The three lines lemma yields the assertion with a constant of the form $$c(n) \, \sqrt{\mathrm{k}_n(\varepsilon)} \, \frac{e^{\delta ((1-\theta) \theta^2 + \theta (1-\theta)^2)}}{\delta^{\frac{n}{2}}}.$$ If $\alpha = \frac{\varepsilon}{\theta} (1+it)$ and $u_\alpha(\xi) = |\xi|^\alpha$, our second claim follows from the simple inequality
\begin{eqnarray*}
\lefteqn{\hskip-15pt \big\| \mathcal{L}_{1+it}(F) \big\|_{\stackrel{\hskip-19pt \circ}{\mathsf{H}_{[\frac{n}{2}]+1}^2}(\Omega)}} \\ & \lesssim & c(n) |e^{\delta(1+it - \theta)^2} | \sup_{0 \le |\beta| \le [\frac{n}{2}]+1} \| \partial_\beta u_\alpha \|_{L_\infty(\Omega)} \big\| F(1+it) \big\|_{\mathsf{H}_{[\frac{n}{2}]+1}^2(\Omega)}. 
\end{eqnarray*}
A similar argument shows that the map $f \mapsto | \ |^{it} f$ is contractive on $L_2(\Omega)$ and bounded on $\mathsf{H}_{[\frac{n}{2}]+1}^2(\Omega)$ up to a constant $c(n) (1 + |t|)^{[\frac{n}{2}]+1}$. By complex interpolation $$\big\| | \ |^{it} f \big\|_{\mathsf{H}_{\frac{n}{2}}^2(\Omega)} \, \le \, c(n) (1 + |t|)^{\frac{n}{2}} \|f\|_{\mathsf{H}_{\frac{n}{2}}^2}.$$ Since $\varepsilon/\theta \in \{\frac12, 1\}$, we easily obtain our first claim from the estimate above. \fin

\begin{remark} \label{ConstantsSobolev}
\emph{A similar argument works for arbitrary unimodular groups. On the other hand, a quick look at the constant we obtain in the proof of Corollary \ref{SobVsClassical} shows that our Sobolev type norm is more appropriate than the classical one in terms of dimensional behavior. Namely, it is easily checked that the constant $c(n)$ above grows linearly with $n$. In particular, we obtain a constant $$c(n) \sqrt{\mathrm{k}_n(\varepsilon)} \sim \frac{n \, \pi^{n/4}}{\sqrt{\Gamma(n/2)}} \, \frac{1}{\sqrt{\varepsilon}}$$ which decreases to $0$ very fast with $n$. We pay a price for small $\varepsilon$ though.}
\end{remark}

\begin{remark} \label{CRdFS}
\emph{The Coifman--Rubio de Francia--Semmes theorem \cite{CRdFS} shows that functions in $\R$ of bounded $2$-variation define $L_p$-bounded Fourier multipliers for $1 < p < \infty$. In this $1$-dimensional setting it can be proved that our abstract Sobolev condition implies bounded $2$-variation. In summary, if we set $\mathsf{HM}_\R$ for the class of H\"ormander-Mihlin multipliers in $\R$, $\mbox{$\psi$-}\mathsf{Riesz}_\R$ for the multipliers in $\R$ satisfying the hypotheses of Theorem B1 and $\mathsf{CRS}_\R$ for the Coifman--Rubio de Francia--Semmes class, Corollary \ref{SobVsClassical} and the comment above yield $$\mathsf{HM}_\R \subset \mbox{$\psi$-}\mathsf{Riesz}_\R \subset \mathsf{CRdFS}_\R.$$ In higher dimensions, Xu extended the notion of $q$-variation to generalize CRS theorem \cite{Xu}. Although we do not know how to compare our condition in Theorem B1 with Xu's one, ours seems easier to check in many cases.}
\end{remark}

\begin{remark}
\emph{Corollary \ref{SobVsClassical} implies $$\|m\|_{\mathsf{M}_p(\widehat{\mathbf{G}})} \, \lesssim_{c(p,n)} \, |m(e)| + \sup_{0 \le |\beta| \le [\frac{n}{2}]+1} \big\| | \ |^{|\beta|} \partial_\beta \widetilde{m} \big\|_\infty.$$ This already improves the main result in \cite{JMP0} for $1 < p < \infty$. In the case of unimodular groups whose Haar mease does not have an atom at $e$, the term $|m(e)|$ can be removed from Theorem B1, Corollary \ref{SobVsClassical} and the inequality above.}
\end{remark}

\subsection{A dimension free formulation}

A quick look at our argument for Theorem B1 shows that the only dependence of the constant we get on the dimension of the cocycle comes from the use of our Littlewood-Paley inequalities in Lemma \ref{Littlewood-Paley2}. The proof of Theorem B2 just requires to replace that result by Lemma \ref{Littlewood-Paley3}.

\demBB As in the statement, let $\psi$ be a conditionally negative length whose associated cocycle $(\H_\psi, \alpha_\psi, b_\psi)$ is finite-dimensional; let $\widetilde{m}$ be a lifting multiplier for that cocycle, so that $m = \widetilde{m} \circ b_\psi$; and let $\varphi: \R_+ \to \C$ be an analytic function in the class $\mathcal{J}$ considered in Paragraph \ref{LPEstimates}. Arguing as in Theorem B1 we may assume with no loss of generality that $m$ vanishes where $\psi$ does. In other words, certain noninjective cocycles may be used as long as $m$ is constant on a nontrivial subgroup of $\G$. Since $ds/s$ is dilation invariant $$m(g) \, = \, \mathrm{k}_\varphi \int_{\R_+} \varphi(s \psi(g))^3 m(g) \, \frac{ds}{s} \, = \, \mathrm{k}_\varphi \int_{\R_+} \varphi(s \psi(g))^2 m_s(g) \, \frac{ds}{s}$$ for some constant $\mathrm{k}_\varphi \neq 0$. Note that $$m_s(g) = \varphi(s \psi(g)) m(g) \, \Rightarrow \, m_s = \widetilde{m}_s \circ b_\psi \ \mbox{with} \ \widetilde{m}_s = \varphi(s | \cdot |^2) \widetilde{m}.$$ Then we proceed again as in Theorem B1. Indeed, we know by assumption that $$\esssup_{s > 0} \Big\| \mathsf{D}_{\frac{\dim \H_\psi}{2} + \varepsilon} \big( \sqrt{\psi_\varepsilon} \, \widetilde{m}_s \big) \Big\|_2 \, < \, \infty.$$ According to the proof Lemma \ref{HMLemma2}, this implies that we may write $m_s$ as a Riesz transform $m_{\gamma,h_s}$ with respect to the length function $\gamma = \psi_\varepsilon \circ b_\psi$ for almost every $s > 0$. Since the families $\varphi_{1}, \varphi_{2} = \varphi$ satisfy the assumptions of Lemma \ref{Littlewood-Paley2}, we may combine this result with the proof of Lemma \ref{HMLemma} to obtain
\begin{eqnarray*}
\|m\|_{\mathsf{M}_p(\widehat{\mathbf{G}})} \!\!\!\! & = & \!\!\! \Big\| \int_{\R_+} \varphi(s \psi(\cdot))^2 \, m_{\gamma,h_s} \, \frac{ds}{s} \Big\|_{\mathsf{M}_p(\widehat{\mathbf{G}})} \\ \!\!\!\! & \le & \!\!\! c(p) \esssup_{s > 0} \|h_s\|_{\H_{\gamma}} \le c(p) \esssup_{s > 0} \Big\| \mathsf{D}_{\frac{\dim \H_\psi}{2} + \varepsilon} \big( \sqrt{\psi_\varepsilon} \, \varphi(s |\cdot|^2) \, \widetilde{m} \big) \Big\|_2. \hskip8pt \square 
\end{eqnarray*}

\begin{remark} \label{Abstract1}
\emph{Lemma \ref{HMLemma} admits generalizations for any conditionally negative length $\gamma: \R^n \to \R_+$ whose associated measure $\nu$ is absolutely continuous with respect to the Lebesgue measure and such that $1 - \cos(2\pi \langle \xi, \cdot \rangle) \in L_1(\R^n, \nu)$ for all $\xi \in \R^n$. Indeed, if we set $d\nu(x) = u(x) dx$ we also have $$\|m\|_{\mathsf{M}_p(\R^n)} \, \lesssim_{c(p)} \, \Big\| \frac{1}{\sqrt{u}} \widehat{\sqrt{\gamma} \, m} \Big\|_{L_2(\R^n)}$$ as long as the space determined by the right hand side admits a dense subspace of functions $m$ for which $\sqrt{\gamma} \, m$ satisfies the Fourier inversion theorem. The Schwartz class was enough for our choice $(\gamma,\nu) = (\psi_\varepsilon,\mu_\varepsilon)$. Lemma \ref{HMLemma2} can also be extended when $\nu$ is invariant under the action $\alpha_\psi$ of the chosen finite-dimensional cocycle $(\H_\psi, \alpha_\psi, b_\psi)$. This invariance is necessary to make sure that the construction in the proof of Lemma \ref{HMLemma2} yields a well-defined cocycle out of $\gamma$ and $\psi$.}  
\end{remark}

\begin{remark} \label{Abstract2}
\emph{Although the above mentioned applications are of independent interest, it is perhaps more significant to read our approach as a way to relate certain kernel reproducing formulas with some differential operators/Sobolev norms in von Neumann algebras. Namely, consider any  length $$\psi(g) = \tau_\psi \Big( \big( 2 \lambda(e) - \lambda(g) - \lambda(g^{-1}) \big) \Big) \quad \mbox{with} \quad \tau_\psi(f) = \tau(f \omega_\psi)$$ for some positive invertible density $\omega_\psi$, and construct the spaces 
\begin{eqnarray*}
L_2(\widehat{\mathbf{G}},\tau_\psi) & = & \Big\{ \hskip2pt h \hskip1pt \in L_0(\widehat{\mathbf{G}}) \ \big| \ \|h\|_{2,\psi} = \tau_\psi (|h|^2)^{\frac12} < \infty \Big\}, \\ \mathsf{W}_\psi^2(\widehat{\mathbf{G}},\tau) & = & \Big\{ m \in \ell_\infty(\G) \ \big| \ \|m\|_{\mathsf{W},\psi} = \Big\| \lambda \big( \sqrt{\psi} \, m \big) \, \frac{1}{\sqrt{\omega_\psi}} \Big\|_{L_2(\widehat{\mathbf{G}})} < \infty \Big\}.
\end{eqnarray*}
Then, there exists a linear isometry $$\Lambda_\psi: h \in L_2^\circ(\widehat{\mathbf{G}},\tau_\psi) \mapsto \frac{\tau_\psi (b_\psi(\cdot) h )}{\sqrt{\psi(\cdot)}} \in \mathsf{W}_\psi^2(\widehat{\mathbf{G}},\tau),$$ where $L_2^\circ(\widehat{\mathbf{G}},\tau_\psi)$ denotes the subspace of mean-zero elements. The surjectivity of $\Lambda_\psi$ depends as above on the existence of a dense subspace of our Sobolev type space admitting Fourier inversion in $\V$ (in fact, this can be used in the opposite direction to identify nice Schwartz-type classes in group algebras). Whenever the map $\Lambda_\psi$ is surjective, it relates the \lq\lq Sobolev norm\rq\rq${}$ of a Riesz transform with the $L_2$-norm of its symbol, which in turn dominates its multiplier norm up to $c(p)$. We have not explored applications in this general setting.}  
\end{remark}

\begin{remark} \label{Inf-dimTh}
\emph{Using Laplace transforms in the spirit of Stein \cite{St}, we can prove Littlewood-Paley estimates in discrete time and also $L_p$ bounds for smooth Fourier multipliers even for infinite-dimensional cocycles. These results will appear in a forthcoming paper.}
\end{remark}

\subsection{Limiting Besov type conditions} \label{Besov}

Let $\eta$ be a radially decreasing smooth function on $\R^n$ with $\chi_{\mathrm{B}_1(0)} \le \eta \le \chi_{\mathrm{B}_2(0)}$ and set $\varphi(\xi) = \eta(\xi) - \eta(2\xi)$, so that me may construct a standard Littlewood-Paley partition of unity $\varphi_k(\xi) = \varphi(2^{-k} \xi)$ for $k \in \Z$. Consider the function $\phi = 1 - \sum_{j \ge 1} \varphi_j$ and let $\alpha \in \R$, $1 \le p,q \le \infty$. Then the Besov space $B_{\alpha q}^p(\R^n)$ and its homogeneous analogue are defined as subspaces of tempered distributions $f \in \mathcal{S}'(\R^n)$ in the following way 
\begin{eqnarray*}
\stackrel{\hskip-5pt \circ}{B_{\alpha q}^p} \hskip-3pt (\R^n) & = & \Big\{ f \in \mathcal{S}'(\R^n) \ \big| \ ||| f |||_{\alpha q}^p \, = \, \Big( \sum_{k \in \Z} 2^{kq \alpha} \| \widehat{\varphi}_k * f\|_p^q \Big)^{\frac{1}{q}} < \infty \Big\}, \\ B_{\alpha q}^p (\R^n) & = & \Big\{ f \in \mathcal{S}'(\R^n) \ \big| \ \|f\|_{\alpha q}^p = \|\widehat{\phi} * f \|_p + \Big( \sum_{k \ge 1} 2^{kq\alpha} \| \widehat{\varphi}_k * f \|_p^q \Big)^{\frac{1}{q}} < \infty \Big\}.
\end{eqnarray*}
Note that $\| \ \|_{\alpha q}^p$ is a norm while $||| \ |||_{\alpha q}^p$ is a seminorm. Besov spaces refine Sobolev spaces in an obvious way. For instance, it is straightforward to show that we have $$B_{\alpha 2}^2(\R^n) \simeq \mathsf{H}_\alpha^2(\R^n) \quad \mbox{and} \quad \stackrel{\hskip-5pt \circ}{B_{\alpha 2}^2} \hskip-3pt (\R^n) \simeq \hskip3pt \stackrel{\hskip-6pt \circ}{\mathsf{H}_\alpha^2} \hskip-3pt (\R^n)$$ with constants depending on the dimension $n$. It is a very natural question to study how can we modify the Sobolev $(\frac{n}{2}+\varepsilon)$-condition in the H\"ormander-Mihlin theorem when $\varepsilon$ approaches $0$. This problem has been studied notably by Seeger \cite{Se3,Se2}, see also \cite{C2,Lz,Se1} and the references therein. In terms of $L_p$-bounded Fourier multipliers for $1 < p < \infty$ the best known result is $$\|m\|_{\mathsf{M}_p(\R^n)} \, \lesssim_{c(p,n)} \, \sup_{j \in \Z} \big\| \varphi_0 \, m(2^j \cdot) \big\|_{B_{\frac{n}{2},1}^2(\R^n)}^{\null},$$ where $\varphi_0(\xi) = \varphi(\xi) = \eta(\xi) - \eta(2\xi)$ as defined above. Of course, we may not expect to replace the Besov space in the right hand side by the bigger one $B_{n/2,2}^2(\R^n)$ or its homogeneous analogue $$\sup_{j \in \Z} \Big( \sum_{k \in \Z} 2^{nk} \big\| \widehat{\varphi}_k * \big( \varphi_0 m(2^j \cdot) \big) \big\|_2^2\Big)^{\frac12} \, = \, \sup_{j \in \Z} \Big( \sum_{k \in \Z} 2^{nk} \big\| \widehat{\varphi}_k * \big( \varphi_j m \big) \big\|_2^2\Big)^{\frac12}.$$ The identity above holds from the dilation invariance (homogeneity) of the norm we use. In the following theorem we show that a log-weighted form of this space is included in $\mathsf{M}_p(\R^n)$. Similar results for Euclidean multipliers were already studied by Baernstein/Sawyer in \cite{BS}, the main difference being that they impose a more demanding $\ell_1$-Besov condition. Our argument is also very different. We formulate our result in the context of discrete group von Neumann algebras, although it also holds for arbitrary unimodular groups. The main idea is to replace the former measures $d\mu_\varepsilon(x) = \omega_\varepsilon (x) dx$ with $\omega_\varepsilon(x) = |x|^{-(n+2\varepsilon)}$ used to prove Theorem B1, by the limiting measure $d\nu(x) = u(x) dx$ with $$u(x) = \frac{1}{|x|^n} \Big( \chi_{\mathrm{B}_1(0)}(x) + \frac{1}{1 + \log^2 |x|} \chi_{\R^n \setminus \mathrm{B}_1(0)}(x) \Big).$$ Let us also consider the associated length $$\gamma(\xi) \, = \, 2 \int_{\R^n} \big( 1 - \cos(2 \pi \langle \xi, x \rangle) \big) \, u(x) \, dx.$$ Remarks below provide more examples and a comparison with Seeger's results. 

\begin{lemma} \label{Calculation}
The length above satisfies $$\gamma(\xi) \, \sim \, \frac{1}{1+|\log |\xi||} \chi_{\mathrm{B}_1(0)}(\xi) + \big( 1 + |\log|\xi|| \big) \chi_{\R^n \setminus \mathrm{B}_1(0)}(\xi).$$
\end{lemma}

\dem This follows from elementary calculations. \fin

\begin{theorem} \label{BesovThm}
Let $(\G,\psi)$ be a discrete group with a conditionally negative length giving rise to an $n$-dimensional cocycle $(\H_\psi, \alpha_\psi, b_\psi)$. Let $(\varphi_j)_{j \in \Z}$ denote a standard radial Littlewood-Paley partition of unity in $\R^n$. Then, if $1 < p < \infty$ the following estimate holds 
\begin{eqnarray*}
\|m\|_{\mathsf{M}_p(\widehat{\mathbf{G}})} & \lesssim_{c(p,n)} & |m(e)| + \inf_{m = {\widetilde{m}}^{\null} \circ b_\psi} \left\{ \sup_{j \in \Z} \Big\| \frac{1}{\sqrt{u}} \big( \sqrt{\gamma} \, \varphi_j \, \widetilde{m} \big)^{\wedge} \Big\|_{L_2(\R^n)} \right\} \\ & \sim_{c(p,n)} & |m(e)| + \inf_{m = {\widetilde{m}}^{\null} \circ b_\psi} \left\{ \sup_{j \in \Z} \Big( \sum_{k \in \Z} 2^{nk} \mathrm{w}_k \big\| \widehat{\varphi}_k * ( \sqrt{\gamma} \, \varphi_j \, \widetilde{m} ) \big\|_2^2 \Big)^{\frac12} \right\},
\end{eqnarray*}
where $u, \gamma$ are as above and the weights $\mathrm{w}_k$ are of the form $\delta_{k \le 0} + k^2 \delta_{k > 0}$ for $k \in \Z$.
\end{theorem}

\dem Arguing as in Theorem B1 we may assume with no loss of generality that $m$ vanishes where $\psi$ does, so that $m(e)=0$. 
According to Remark \ref{Abstract1} and Lemma \ref{Calculation}, Lemmas \ref{HMLemma} and \ref{HMLemma2} apply in our setting $(\gamma,\nu,u)$. Moreover, according to our argument for Theorem B1, we find $$m(g) \, = \, \sum_{j \in \Z} (\phi_j(\psi(g)))^2 \, (\varphi_j \widetilde{m})(b_\psi(g))$$ for certain smooth partition of unity $\phi_j$. If $$\sup_{j \in \Z} \Big\| \frac{1}{\sqrt{u}} \big( \sqrt{\gamma} \, \varphi_j \, \widetilde{m} \big)^{\wedge} \Big\|_2 \, < \, \infty,$$ the (generalized) proofs of Lemmas \ref{HMLemma} and \ref{HMLemma2} imply that $m_j = (\varphi_j \widetilde{m}) \circ b_\psi$ is a Riesz transform $m_{\zeta,h_j}$ with respect to the length function $\zeta = \gamma \circ b_\psi$, which comes from a composed cocycle as in Lemma \ref{HMLemma2} thanks to the orthogonal invariance of $\nu$ (radial density). Since the families $\varphi_{1j} = \varphi_{2j} = \phi_j$ satisfy the assumptions of Lemma \ref{Littlewood-Paley2}, we combine this with the proof of Lemma \ref{HMLemma} to obtain
\begin{eqnarray*}
\|m\|_{\mathsf{M}_p(\widehat{\mathbf{G}})} & = & \Big\| \sum_{j \in \Z} \varphi_{1j}(\psi(\cdot)) \varphi_{2j}(\psi(\cdot)) m_{\zeta,h_j} \Big\|_{\mathsf{M}_p(\widehat{\mathbf{G}})} \\ & \le & c(p,n) \, \sup_{j \in \Z} \|h_j\|_{\H_{\zeta}} \ \le \ c(p,n) \, \sup_{j \in \Z} \Big\| \frac{1}{\sqrt{u}} \big( \sqrt{\gamma} \, \varphi_j \, \widetilde{m} \big)^{\wedge} \Big\|_2. 
\end{eqnarray*}
This proves the first estimate. The second one is clear from the definition of $u$. \fin

\begin{remark} \label{BesovGral}
\emph{According to Remark \ref{Abstract1}, other limiting measures $\nu$ apply as long as we know that $1 - \cos (2\pi \langle \xi, \cdot \rangle)$ belongs to $L_1(\R^n, \nu)$ and the measure $\nu$ is invariant under the cocycle action $\alpha_\psi$. In particular, any radial measure $d\nu(x) = u(x) dx$ such that $u(s) (s^2 \chi_{(0,1)}(s) + \chi_{(1,\infty)}(s)) \in L_1(\R_+,ds)$ satisfies these conditions. Note that any such measure will provide an associated length of polynomial growth, so that the associated Sobolev type spaces admits the Schwartz class as a dense subspace satisfying the Fourier inversion formula, as demanded by Remark \ref{Abstract1}. If fact, it is conceivable that for slow-increasing lengths $\gamma$ $$\big\| \widehat{\varphi}_k * \big( \sqrt{\gamma} \varphi_j \widetilde{m} \big) \big\|_2^2 \, \sim \, \gamma(2^j) \big\| \widehat{\varphi}_k * \big( \varphi_j \widetilde{m} \big) \big\|_2^2.$$ Therefore, under this assumption we would finally get $$\|m\|_{\mathsf{M}_p(\widehat{\mathbf{G}})} \, \lesssim_{c(p,n)} \, |m(e)| + \inf_{m = {\widetilde{m}}^{\null} \circ b_\psi} \left\{ \sup_{j \in \Z} \Big( \sum_{k \in \Z} 2^{nk} \frac{\gamma(2^j)}{u(2^k)} \big\| \widehat{\varphi}_k * ( \varphi_j \, \widetilde{m} ) \big\|_2^2 \Big)^{\frac12} \right\}.$$}
\end{remark}

\begin{remark} \label{BesovAnalysis}
\emph{The Besov space $B_{n/2,1}^2(\R^n)$ used in Seeger's result is not dilation invariant. Thus, in order to compare his estimates with ours, we first need to dilate from $\varphi_0 m(2^j \cdot)$ to $\varphi_j m$. Easy calculations yield 
\begin{eqnarray*}
\big\| \varphi_0 m(2^j \cdot) \big\|_{B_{\frac{n}{2},1}^2(\R^n)} & = & \big\| \widehat{\phi} * \big( \varphi_0 m(2^j \cdot) \big) \big\|_2 + \sum_{k \ge 1} 2^{\frac{nk}{2}} \big\| \widehat{\varphi}_k * \big( \varphi_0 m(2^j \cdot) \big) \big\|_2 \\ & \sim & \Big( \sum_{k+j \le 0} \big\| \widehat{\varphi}_k * \big( \varphi_j m \big) \big\|_2^2 \Big)^{\frac12} + \sum_{k+j \ge 1} 2^{\frac{nk}{2}} \big\| \widehat{\varphi}_k * \big( \varphi_j m \big) \big\|_2. 
\end{eqnarray*}
On the other hand, our estimate in Theorem \ref{BesovThm} gives $$\Big( \sum_{k \le 0} 2^{nk} \big\| \widehat{\varphi}_k * \big( \sqrt{\gamma} \varphi_j m \big) \big\|_2^2 \Big)^{\frac12} + \Big( \sum_{k \ge 1} 2^{nk} k^2 \big\| \widehat{\varphi}_k * \big( \sqrt{\gamma} \varphi_j m \big) \big\|_2^2 \Big)^{\frac12}.$$ It seems we get better estimates for $k \le \min(0,-j)$ and worse for $k > \max(0,-j)$.}
\end{remark}

\section{\textbf{Analysis in free group branches}}

In this section we prove Theorem C. We shall use the same terminology as in Section \ref{ExSect} for the natural cocycle of $\F_\infty$ associated to the word length $| \ |$. In particular, recall that the Hilbert space is $\H_{| \ |} = \R[\F_\infty] / \R \delta_e$, with $\alpha_{| \ |} = \lambda$ and $b_{| \ |}(g) = \delta_g + \R \delta_e$. The system $\xi_g = \delta_g - \delta_{g^-}$ with $g \neq e$ forms an orthonormal basis of $\H_{| \ |}$, where $g^-$ is the word which results after deleting the last generator on the right of $g$. By a branch of $\F_\infty$, we mean a subset $\mathrm{B} = (g_k)_{k \ge 1}$ with $g_k=g_{k+1}^-$. We will say that $g_1$ is the root of $\mathrm{B}$.

\demC Given $t > 0$, let $$\widetilde{m}_t(j) = tje^{-tj} \widetilde{m}(j).$$ Our hypotheses on $\widetilde{m}$ imply that $$|\widetilde{m}_t(j)-\widetilde{m}_t(j-1)| \lesssim te^{-tj} + te^{-\frac{tj}2} \lesssim te^{-\frac{tj}2}.$$ In particular, we find the following estimate 
\begin{eqnarray*}
\lefteqn{\hskip-10pt \sup_{t > 0} \sum_{g \in \mathrm{B}} \big| \widetilde{m}_t(|g|) \sqrt{|g|} - \widetilde{m}_t(|g^-|) \sqrt{|g^-|} \big|^2} \\ & \lesssim & \sup_{t > 0} \sum_{j \ge 1} \big| \widetilde{m}_t(j)-\widetilde{m}_t(j-1) \big|^2 j + |\widetilde{m}_t(j-1)|^2 \frac{1}{j} \\ & \lesssim & \sup_{t > 0} \int_{\R_+} t^2e^{-ts}sds + \int_{\R_+} (ts)^2e^{-2ts} \frac{ds}{s} = \, \int_{\R_+} x ( e^{-x} + e^{-2x} ) \, dx \, < \, \infty.
\end{eqnarray*}
Therefore, if we define $h_t = \sum_{g \neq e} \langle h_t, \xi_g \rangle_{\H_{| \ |}} \xi_g$ with $$\langle h_t,\xi_g \rangle_{\H_{| \ |}} = \delta_{g \in \mathrm{B}} \Big( \widetilde{m}_t(|g|)\sqrt{|g|}-\widetilde{m}_t(|g^-|)\sqrt{|g^-|} \Big),$$ it turns out that $(h_t)_{t > 0}$ is uniformly bounded in $\H_{| \ |}$. Moreover $$\delta_g = \sum_{h \le g} \xi_h = b_{| \ |}(g) \ \Rightarrow \ \widetilde{m}_t(|g|) = \frac{\langle b_{| \ |}(g), h_t \rangle_{\H_{| \ |}}}{\sqrt{|g|}} \quad \mbox{for all} \quad g \in \mathrm{B}.$$ In other words, $\widetilde{m}_t \circ | \ |$ coincides on $\mathrm{B}$ with the Fourier symbol of the Riesz transform $R_{| \ |, h_t}$ associated to the word length $| \ |$ in the direction of $h_t$. Assume (wlog) $p \ge 2$. Let us now fix $f \in L_p(\mathcal{L}(\F_\infty))$ with vanishing Fourier coefficients outside $\mathrm{B}$. By Littlewood-Paley theory \cite{JLX} of the Poisson semigroup $e^{-tA_{| \, |}}: \lambda(g) \mapsto e^{-t|g|} \lambda(g)$
\begin{eqnarray*}
\|T_{\widetilde{m}}f\|_p & \lesssim_{c(p)} & \Big( \int_{\R_+} \big| (tA_{| \, |})^2 e^{-2tA_{| \, |}} T_{\widetilde{m}}f \big|^2 + \big| (tA_{| \, |})^2 e^{-2tA_{| \, |}} T_{\widetilde{m}}f^* \big|^2 \, \frac {dt}t \Big)^\frac12 \\ & = & \Big( \int_{\R_+} \big| R_{| \ |,h_t}(tA_{| \, |}e^{-tA_{| \, |}}) f \big|^2 + \big| R_{| \ |,h_t} (tA_{| \, |}e^{-tA_{| \, |}})f^* \big|^2 \, \frac {dt}t \Big)^\frac12
\end{eqnarray*}
since $T_{{\widetilde{m}}_t} \!\! = \! R_{| \ |, h_t}$ on $\mathrm{B}$. Apply then the integral version of Lemma \ref{Littlewood-Paley} and \cite{JLX}. \fin 

We will say that a family of branches $\mathcal{T} = \{\mathrm{B}_k \, : \, k \ge 1\}$ forms a partition of the free group when $\F_\infty = \{e\}\cup (\cup_{k} \mathrm{B}_k)$ and the $\mathrm{B}_k$'s are pairwise disjoint. Given a branch $\mathrm{B} \in \mathcal{T}$ let us write $g_{\mathrm{B},1}$ for its root. We will say that $\mathrm{B}$ is a principal branch when its root $|g_{\mathrm{B},1}|=1$. If $\mathrm{B}$ is not a principal branch, there is a unique branch $\mathrm{B}_-$ in $\mathcal{T}$ which contains $g_{\mathrm{B},1}^-$. Given $g \in \F_\infty$, define $\Pi_{\mathrm{B}}g$ to be the biggest element in $\mathrm{B}$ which is smaller than or equal to $g$. If there is no such element, set $\Pi_{\mathrm{B}}g = e$. Now let us fix a standard Littlewood-Paley partition of unity in $\R_+$. That is, given $\eta: \R_+ \to \R_+$ a smooth decreasing function with $\chi_{(0,1)} \le \eta \le \chi_{(0,2)}$, set $\phi(\xi) = \eta(\xi) - \eta(2\xi)$ and $\phi_k(\xi) = \phi(2^{-k} \xi)$ for $k \in \Z$. Assume in addition that $\sqrt{\phi}$ is Lipschitz (as we assume in Theorem C). Then construct $$\varphi_1 = \sum_{j \le 1} \phi_j \quad \mbox{and} \quad \varphi_j = \phi_j$$ for $j \ge 2$, so that $\sum_{j \ge 1} \varphi_j = 1$.
Define $\Lambda_j: \lambda(g) \mapsto \sqrt{\varphi_j(|g|)} \lambda(g)$ and $$\Lambda_{j,k}: \, \lambda(g) \ \mapsto \ \frac{\sqrt{\varphi_j(|\Pi_{\mathrm{B}_k}g|)|\Pi_{\mathrm{B}_k}g|} - \sqrt{\varphi_j(|\Pi_{\mathrm{B}_k^-}g|)|\Pi_{\mathrm{B}_k^-}g|}}{\sqrt{|g|}} \lambda(g)$$ for any $g \in \F_\infty$ with the convention that $\Pi_{\mathrm{B}_k^-}g=e$ if $\mathrm{B}_k^-$ does not exist. 

\begin{lemma} \label{BranchLPLemma}
If $1 < p < 2$ and $f \in L_p(\mathcal{L}(\F_\infty))$ $$\inf_{\Lambda_{j,k}f = a_{j,k}+b_{j,k}} \Big\| \Big( \sum_{j,k} a_{j,k}^* a_{j,k}^{\null} + \widetilde{b}_{j,k}^{\null} \widetilde{b}_{j,k}^* \Big)^{\frac12} \Big\|_p \, \lesssim_{c(p)} \, \|f\|_p.$$
\end{lemma}

\dem Define $h_{j,k} = \sum_{g \neq e} \langle h_{j,k}, \xi_g \rangle_{\H_{| \ |}} \xi_g$ with $$\langle h_{j,k},\xi_g\rangle = \delta_{g \in \mathrm{B}_k} \Big( \sqrt{\varphi_j(|g|){|g|}} - \sqrt{ \varphi_j(|g^-|)|g^-|} \Big).$$ To show that $h_{j,k} \in \H_{| \ |}$, we note that $\varphi_1$ and $\varphi_j$ ($j \ge 2$) are supported by $[0,4]$ and $[2^{j-1},2^{j+1}]$ respectively. Therefore, arguing as we did in the proof of Theorem C i), we obtain (using that $\sqrt{\phi}$ is Lipstchitz)
\begin{eqnarray*}
\|h_{j,k}\|_{\H_{| \ |}}^2 & \lesssim & \sum_{g \in \mathrm{B}_k} \big| \sqrt{\varphi_j(|g|)} - \sqrt{\varphi_j(|g^-|)} \big|^2 |g| + |\varphi_j(|g^-|)| \big( \sqrt{|g|} - \sqrt{|g^-|} \big)^2 \\ & \lesssim & \sum_{2^{j-1} \le i \le 2^{j+1}+1} \big| \sqrt{\varphi_j(i)} - \sqrt{\varphi_j(i-1)} \big|^2 i + |\varphi_j(i-1)| \frac1i \\ & \lesssim & \sum_{2^{j-1} \le i \le 2^{j+1}+1} \frac{i}{4^j} + \frac1i \ \lesssim \ 1.
\end{eqnarray*}
In particular, the family $(h_{j,k})$ is uniformly bounded in $\H_{| \ |}$ and $$\delta_g = \sum_{h \le g} \xi_h = b_{| \ |}(g) \ \Rightarrow \ \Lambda_{j,k}f \, = \, \sum_{g \in \F_\infty} \frac{\langle b_{| \ |}(g), h_{j,k} \rangle_{\H_{| \ |}}}{\sqrt{|g|}} \widehat{f}(g) \lambda(g) \, = \, R_{| \ |, h_{j,k}}f.$$ Now, according to the definition of $h_{j,k}$ it is easily checked that $\langle h_{j,k}, h_{j',k'} \rangle_{\H_{| \ |}}$ vanishes unless $k=k'$ and $|j - j'| \le 1$. This implies that the subsystems $(h_{2j,k})$ and $(h_{2j+1,k})$ are orthogonal and uniformly bounded above and below. Once this is known and splitting into two systems, the assertion follows immediately from Theorem A1 by standard considerations. The proof is complete. \fin

\demCC It clearly suffices to prove the result for $\mathrm{B}$ a principal branch. Let $\mathcal{T} = \{ \mathrm{B}_k : k \ge 1\}$ form a partition of $\F_\infty$ which contains $\mathrm{B}_1 = \mathrm{B}$ as a principal branch. Given $f \in L_p(\mathcal{L}(\F_\infty))$ with vanishing Fourier coefficients outside $\mathrm{B}$, it is then easily checked that $\Lambda_jf = \Lambda_{j,1}f$ and $\Lambda_{j,k}f=0$ for other values of $k$. The lower estimate then follows from Lemma \ref{BranchLPLemma}. For the upper estimate, we use the fact that $\varphi_j$ is a partition of unity together with the inequality in Remark \ref{MoreLPRem} iv). Namely, we obtain 
\begin{eqnarray*}
\|f\|_{L_p(\widehat{\mathbf{B}})} & = & \Big\| \sum_{j \ge 1} \Lambda_j^2 f \Big\|_{L_p(\widehat{\mathbf{B}})} \\
& = & \Big\| \sum_{j \ge 1} R_{| \ |,h_{j,1}} (a_j + b_j) \Big\|_{L_p(\widehat{\mathbf{B}})} \\ & \le & \Big\| \Big( \langle h_{j,1}, h_{k,1} \rangle \Big) \Big\|_{\mathcal{B}(\ell_2)}^{\frac12} \Big\| \Big( \sum_{j \ge 1} a_j^*a_j^{\null} + b_j^{\null} b_j^* \Big)^\frac12 \Big\|_{L_p(\mathcal{L}(\F_\infty))}.
\end{eqnarray*}
Now observe that $\langle h_{j,1}, h_{k,1} \rangle$ vanishes when $|j-k|>1$ and is bounded above and below otherwise. This  shows that the matrix above is bounded. We are done. \fin

\begin{corollary} \label{UpperBranch}
If $\mathrm{B}$ is any branch of $\F_\infty$ and $2 < p < \infty$ $$\Big\| \Big( \sum_{j \ge 1} |\Lambda_jf|^2 \Big)^\frac12 \Big\|_{L_p(\widehat{\mathbf{B}})} \, \lesssim_{c(p)} \, \|f\|_{L_p(\widehat{\mathbf{B}})},$$ with the multipliers $\Lambda_j: \lambda(g) \mapsto \sqrt{\varphi_j(|g|)} \lambda(g)$ defined as in Theorem \emph{C ii)}.
\end{corollary}

\dem It follows from the obvious extension of Lemma \ref{BranchLPLemma} for $p>2$. \fin

\begin{remark} \label{Bozejko}
\emph{Bo\.zejko-Fendler theorem \cite{BF} shows that Fourier summability fails in $L_p(\mathcal{L}(\F_n))$ when $|\frac12 - \frac{1}{p}| > \frac16$ and the partial sums are chosen to lie in a sequence of increasing balls with respect to the word length. This may be regarded as some sort of Fefferman's disc multiplier theorem \cite{Fef} for the free group algebra although discreteness might allow some room for Fourier summability near $L_2$ in the spirit of Bochner-Riesz multipliers. This result indicates that we might not expect Littlewood-Paley estimates for nontrivial branches arising from sharp (not smooth) truncations in our partitions of unity, as it holds for $\R$ or $\Z$.}
\end{remark}

\section*{{\bf Appendix A. Operator algebraic tools}}
\label{Back}

Along this paper we have used some concepts from noncommutative integration which include noncommutative $L_p$ spaces and sums of independent noncommuting random variables. In the context of group von Neumann algebras, we have also used crossed products, length functions and cocycles. In this section we briefly review these notions for the readers who are not familiar with them.

\vskip3pt

\noindent \textbf{Noncommutative integration.} Part of von Neumann algebra theory has evolved as the noncommutative form of measure theory and integration. A von Neumann algebra \cite{KR,Ta} is a unital weak-operator closed $\mathrm{C}^*$-algebra and, according to the Gelfand-Naimark-Segal theorem, any such algebra $\M$ embeds in the algebra $\mathcal{B}(\mathcal{H})$ of bounded linear operators on a Hilbert space $\mathcal{H}$. We write $\1_\M$ for the unit. The positive cone $\M_+$ is the set of positive operators in $\M$ and a trace $\tau: \M_+ \to [0,\infty]$ is a linear map satisfying $$\tau(f^*f) = \tau(ff^*).$$ It is normal if $\sup_\alpha \tau(f_\alpha) = \tau(\sup_\alpha f_\alpha)$ for bounded increasing nets $(f_\alpha)$ in $\M_+$; it is semifinite if for any non-zero $f \in \M_+$ there exists $0 < f' \le f$ such that $\tau(f') < \infty$; and it is faithful if $\tau(f) = 0$ implies that $f = 0$. The trace $\tau$ plays the r\^ole of the integral in the classical case. A von Neumann algebra is semifinite when it admits a normal semifinite faithful (n.s.f. in short) trace $\tau$. Any operator $f$ is a linear combination $f_1 - f_2 + if_3 - if_4$ of four positive operators. Thus, we can extend $\tau$ to the whole algebra $\M$ and the tracial property can be restated in the familiar form $\tau(fg) = \tau(gf)$. Unless explicitly stated, the pair $(\M,\tau)$ will denote a semifinite von Neumann algebra equipped with a n.s.f. trace. We will refer to it as a \emph{noncommutative measure space}. Note that commutative von Neumann algebras correspond to classical $L_\infty$ spaces.

According to the GNS construction, the noncommutative analog of measurable sets (characteristic functions) are orthogonal projections. Given $f \in \M_+$, the support of $f$ is the least projection $q$ in $\M$ such that $qf = f = fq$ and is denoted by $\mbox{supp} \hskip1pt f$. Let $\mathcal{S}_\M^+$ be the set of all $f \in \M_+$ such that $\tau(\mbox{supp} \hskip1pt f) < \infty$ and set $\mathcal{S}_\M$ to be the linear span of $\mathcal{S}_\M^+$. If we write $|f|=\sqrt{f^*f}$, we can use the spectral measure $d\gamma: \R_+ \to \mathcal{B}(\mathcal{H})$ of $|f|$ to define $$|f|^p = \int_{\R_+} s^p \, d \gamma(s) \quad \mbox{for} \quad 0 < p < \infty.$$ We have $f \in \mathcal{S}_\M \Rightarrow |f|^p \in \mathcal{S}_\M^+ \Rightarrow \tau(|f|^p) < \infty$. If we set $\|f\|_p = \tau( |f|^p )^{\frac1p}$, we obtain a norm in $\mathcal{S}_\M$ for $1 \le p < \infty$ and a $p$-norm for $0 < p < 1$. Using that $\mathcal{S}_\M$ is an involutive strongly dense ideal of $\M$, we define the \emph{noncommutative $L_p$ space} $L_p(\M)$ associated to the pair $(\M, \tau)$ as the completion of $(\mathcal{S}_\M, \| \ \|_p)$. On the other hand, we set $L_\infty(\M) = \M$ equipped with the operator norm. Many fundamental properties of classical $L_p$ spaces like duality, real and complex interpolation, H\"older inequalities, etc hold in this setting. Elements of $L_p(\M)$ can also be described as measurable operators affiliated to $(\M,\tau)$, we refer to Pisier/Xu's survey \cite{PX2} for more information and historical references. Note that classical $L_p$ spaces are denoted in the noncommutative terminology as $L_p(\Omega,\mu) = L_p(\M)$ where $\M$ is the commutative von Neumann algebra $L_\infty(\Omega,\mu)$.

\renewcommand{\theequation}{A.1}
\addtocounter{equation}{-1}

A unital, weakly closed $*$-subalgebra is called a von Neumann subalgebra. A conditional expectation $\mathsf{E}: \mathcal{M}
\to \mathcal{N}$ from a von Neumann algebra $\mathcal{M}$ onto a von Neumann subalgebra $\mathcal{N}$ is a positive contractive projection. It is called normal if the adjoint map $\mathsf{E}^*$ sends $L_1(\mathcal{M})$ to $L_1(\mathcal{N})$. In this case, the restriction map $\mathsf{E}_1 = \mathsf{E}^* \! \! \mid_{L_1(\mathcal{M})}$ satisfies $\mathsf{E}_1^* = \mathsf{E}$. Note that such normal conditional expectation exists if and only if the restriction of $\tau$ to the von Neumann subalgebra $\mathcal{N}$ remains semifinite, see \cite{Ta} for further details. Any such conditional expectation is trace preserving $\tau \circ \mathsf{E} = \tau$ and satisfies the bimodule property $$\mathsf{E}(a_1 b \hskip1pt a_2) = a_1 \mathsf{E}(b) \hskip1pt a_2 \quad \mbox{for all} \quad a_1, a_2 \in \mathcal{N} \ \mbox{and} \ b \in \mathcal{M}.$$ Given von Neumann algebras $\mathcal{N} \subset \mathcal{A}, \mathcal{B} \subset \M$, we will say that $\mathcal{A}, \mathcal{B}$ are independent over $\mathsf{E}$ whenever $\mathsf{E}(ab) = \mathsf{E}(a) \mathsf{E}(b)$ for all $a \in \mathcal{A}$ and $b \in \mathcal{B}$. Similarly, we will say that a family of random variables $(f_j)_{j \in \mathcal{J}}$ in $\M$ is fully independent over $\mathsf{E}$ if the von Neumann algebras generated by any two disjoint subsets of $(f_j)_{j \in \mathcal{J}}$ are independent over $\mathsf{E}$. The noncommutative analog of Rosenthal inequality \cite{Ro} was obtained in \cite{JX3} and read as follows for $p \ge 2$. If the random variables $(f_j)_{j \in \mathcal{J}} \subset L_p(\M)$ satisfy $\mathsf{E}(f_j) = 0$ and are fully independent over $\mathsf{E}$, then we find 
\begin{equation} \label{Rosenthal}
\hskip3pt \frac{1}{p} \Big\| \sum_{j \in \mathcal{J}} f_j \Big\|_p \, \sim \, \Big( \sum_{j \in \mathcal{J}} \|f_j\|_p^p \Big)^{\frac{1}{p}} + \Big\| \Big( \sum_{j \in \mathcal{J}} \mathsf{E}(f_j^* f_j) \Big)^{\frac12} \Big\|_p + \Big\| \Big( \sum_{j \in \mathcal{J}} \mathsf{E}(f_j f_j^*) \Big)^{\frac12} \Big\|_p.
\end{equation} 

\vskip3pt

\noindent \textbf{Group von Neumann algebras.} Let $\G$ be a discrete group with left regular representation $\lambda: \G \to \mathcal{B}(\ell_2(\G))$ given by $\lambda(g) \delta_h = \delta_{gh}$, where the $\delta_g$'s form the unit vector basis of $\ell_2(\G)$. Write $\mathcal{L}(\G)$ for its \emph{group von Neumann algebra}, the weak operator closure of the linear span of $\lambda(\G)$ in $\mathcal{B}(\ell_2(\G))$. Consider the standard trace $\tau(\lambda(g)) = \delta_{g=e}$ where $e$ denotes the identity of $\G$. Any $f \in \mathcal{L}(\G)$ has a Fourier series expansion of the form $$\sum_{g \in \G} \widehat{f}(g) \lambda(g) \quad \mbox{with} \quad \tau(f) = \widehat{f}(e).$$ Define $$L_p(\widehat{\mathbf{G}}) = L_p(\V,\tau) \equiv \mbox{Closure of $\mathcal{L}(\G)$ wrt $\|f\|_{L_p(\widehat{\mathbf{G}})} = \big( \tau |f|^p \big)^\frac1p$},$$ the natural $L_p$ space over the noncommutative measure space $(\mathcal{L}(\mathrm{G}), \tau)$. Note that when $\G$ is abelian we get the $L_p$ space on the dual group equipped with its normalized Haar measure, after identifying $\lambda(g)$ with the character $\chi_g$.

Now, given another noncommutative measure space $(\M, \varphi)$ with $\M \subset \mathcal{B(\H)}$ assume there exists a trace preserving action $\alpha: \G \to \mathrm{Aut}(\M)$. Define the \emph{crossed product algebra} $\M \rtimes_\alpha \G$ as the weak operator closure of the $*$-algebra generated by $\1_\M \otimes \lambda(\G)$ and $\rho(\M)$ in $\mathcal{B}(\ell_2(\G; \H))$. The $*$-representation $\rho: \M \to \mathcal{B}(\ell_2(\G; \H))$ is given by $\rho(f) = \sum_{h \in \G} \alpha_{h^{-1}}(f) \otimes e_{h,h}$, with $e_{g,h}$ the matrix units for $\ell_2(\G)$. A generic element of $\M \rtimes_\alpha \G$ has the form $\sum_{g \in \G} f_g \rtimes_\alpha \lambda(g)$, with $f_g \in \M$. Playing with $\lambda$ and $\rho$, it is clear that $\M \rtimes_\alpha \G$ sits in $\M \bar\otimes \mathcal{B}(\ell_2(\G))$
\begin{eqnarray*}
\summ_g f_g \rtimes_\alpha \lambda(g) & = & \summ_g \rho(f_g) \big( \1_\M \otimes \lambda(g) \big) \\
& = & \summ_{g,h,h'} \big( \alpha_{h^{-1}}(f_g) \otimes e_{h,h} \big) \big( \mathbf{1}_\M \otimes  e_{gh',h'} \big)
\\ & = & \summ_{g,h} \alpha_{h^{-1}}(f_g) \otimes e_{h,g^{-1}h}  \
= \  \summ_{g,h} \alpha_{g^{-1}}(f_{gh^{-1}}) \otimes e_{g,h}.
\end{eqnarray*}
Similar computations lead to
\begin{itemize}
\item $(f \rtimes_\alpha \lambda(g))^* = \alpha_{g^{-1}}(f^*) \rtimes_\alpha \lambda(g^{-1})$,

\item $(f \rtimes_\alpha \lambda(g)) (f' \rtimes_\alpha \lambda(g')) = f \alpha_g(f') \rtimes_\alpha \lambda(gg')$,

\item $\varphi \rtimes_\alpha \tau (f \rtimes_\alpha \lambda(g)) = \varphi \otimes \tau (f \otimes \lambda(g)) = \delta_{g=e} \varphi(f)$.
\end{itemize}
Since $\alpha$ will be fixed, we relax the terminology and write $\sum_g f_g \rtimes \lambda(g) \in \M \rtimes \G$.

\renewcommand{\theequation}{A.2}
\addtocounter{equation}{-1}

Let us now consider semigroups of operators $\S_\psi = (\S_{\psi,t})_{t \ge 0}$ on $\mathcal{L}(\G)$ which act diagonally on the trigonometric system. In other words, $\S_{\psi,t} : \lambda(g) \mapsto e^{-t \psi(g)} \lambda(g)$ for some function $\psi: \G \to \R_+$. $\S_\psi = (\S_{\psi,t})_{t \ge 0}$ defines a \emph{noncommutative Markov semigroup} when:
\begin{itemize}
\item[i)] $\S_{\psi,t}(\1_\V) = \1_\V$ for all $t \ge 0$,

\vskip2pt

\item[ii)] Each $\S_{\psi,t}$ is normal and completely positive on $\V$,

\vskip2pt

\item[iii)] Each $\S_{\psi,t}$ is self-adjoint, i.e. $\tau ((\S_{\psi,t}f)^*g) = \tau (f^*(\S_{\psi,t}g))$ for $f,g \in \V$,

\vskip2pt

\item[iv)] $\S_{\psi,t}f \rightarrow f$ as $t \rightarrow 0^+$ in the weak-$*$ topology of $\V$.
\end{itemize}
These conditions are reminiscent of Stein's notion of diffusion semigroup \cite{St}. They imply that $\S_{\psi,t}$ is completely contractive, trace preserving and also extends to a semigroup of contractions on $L_p(\V)$ for any $1 \le p \le \infty$. As in the classical case $\S_\psi$ always admits an infinitesimal generator $$-A_\psi = \lim_{t \rightarrow 0} \frac{\S_{\psi,t}-id_\V}{t} \ \quad \mbox{with} \quad \ \S_{\psi,t} = \exp(-tA_\psi).$$ In the $L_2$ setting, $A_\psi$ is an unbounded operator defined on $$\mathrm{dom}_2(A_\psi) = \Big\{ f \in L_2(\widehat{\mathbf{G}}) \, \big| \, \lim_{t \rightarrow 0} \frac{\S_{\psi,t}f - f}{t} \in L_2(\widehat{\mathbf{G}}) \Big\}.$$ As an operator in $L_2(\widehat{\mathbf{G}})$, $A_\psi$ is positive and so we may define the subordinated Poisson semigroup $\P_\psi = (\P_{\psi,t})_{t \ge 0}$ by $\P_{\psi,t} = \exp(-t\sqrt{A_\psi})$. This is again a Markov semigroup. Note that $P_t$ is chosen so that $(\partial_t^2-A_\psi)\P_{\psi,t}=0$. 
In general, we let $-A_{\psi,p}$ the generator of the realization of $\S_\psi = (\S_{\psi,t})_{t \ge 0}$ on $L_p(\V)$. It should be noticed that $\mathrm{ker}A_{\psi,p}$ is a complemented subspace of $L_p(\V)$. Let $E_p$ denote the corresponding projection and $J_p = id_{L_p(\widehat{\mathbf{G}})} - E_p$. Consider the complemented subspaces $$L_p^\circ(\widehat{\mathbf{G}}) = J_p(L_p(\widehat{\mathbf{G}})) = \Big\{ f \in L_p(\widehat{\mathbf{G}}) \, \big| \ \lim_{t \to \infty} \S_{\psi,t} f = 0 \Big\}.$$ The associated gradient form or \lq\lq carr\'e du champs\rq\rq${}$ is defined as $$\Gamma_\psi(f_1,f_2) \, = \, \frac12 \Big( A_\psi(f_1^*)f_2 + f_1^*A_\psi(f_2) - A_\psi(f_1^*f_2) \Big).$$ Since $\S_\psi = (\S_{\psi,t})_{t \ge 0}$ is a Fourier multiplier, we get $A_\psi(\lambda(g)) = \psi(g) \lambda(g)$ and 
\begin{equation} \label{Gradient}
\Gamma_\psi(f_1,f_2) \, = \, \sum_{g,h \in \G} \overline{\widehat{f}_1(g)} \widehat{f}_2(h) \, \frac{\psi(g^{-1}) + \psi(h) - \psi(g^{-1}h)}{2} \, \lambda(g^{-1}h).
\end{equation}
The crucial condition $\Gamma_\psi(f,f) \ge 0$ is characterized in the following paragraph.

\vskip3pt

\noindent \textbf{Length functions and cocycles.} A \emph{left cocycle} $(\H, \alpha, b)$ associated to a discrete group $\G$ is a triple given by a Hilbert space $\H$, an isometric action $\alpha: \G \to \mathrm{Aut}(\H)$ and a map $b: \G \to \H$ so that $$\alpha_g(b(h)) = b(gh) - b(g).$$ A \emph{right cocycle} satisfies the relation $\alpha_g(b(h)) = b(hg^{-1}) - b(g^{-1})$ instead. In this paper, we say that $\psi: \G \to \R_+$ is a \emph{length function} if it vanishes at the identity $e$, $\psi(g) = \psi(g^{-1})$ and $$\summ_g \beta_g = 0 \Rightarrow \summ_{g,h} \overline{\beta}_g \beta_h \psi(g^{-1}h) \le 0.$$ Those functions satisfying the last condition are called conditionally negative. It is straightforward to show that length functions take values in $\R_+$. In what follows, we only consider cocycles with values in real Hilbert spaces. Any cocycle $(\H, \alpha,b)$ gives rise to an associated length function $\psi_b(g) = \langle b(g), b(g) \rangle_\H$, as it can be checked by the reader. Reciprocally, any length function $\psi$ gives rise to a left and a right cocycle. This is a standard application of the ideas around Schoenberg's theorem \cite{Sc}, which claims that $\psi: \G \to \R_+$ is a length function if and only if the associated semigroup $\S_\psi= (\S_{\psi,t})_{t \ge 0}$ given by $S_{\psi,t}: \lambda(g) \mapsto \exp(-t\psi(g)) \lambda(g)$ is Markovian on $\V$. Let us collect these well-known results.

\addtocounter{Atheorem}{+2}

\begin{Alemma} \label{Lemmapsi}
If $\psi: \G \to \R_+$ is a length function$\, :$
\begin{itemize}
\item[i)] The Gromov forms
\begin{eqnarray*}
K_\psi^1(g,h) & = & \frac{\psi(g) + \psi(h) - \psi(g^{-1}h)}{2}, \\
K_\psi^2(g,h) & = & \frac{\psi(g) + \psi(h) - \psi(gh^{-1})}{2},
\end{eqnarray*}
define positive matrices on $\G \times \G$ and lead to $$\Big\langle \summ_g a_g \delta_g, \summ_h b_h
\delta_h \Big\rangle_{\psi,j} = \summ_{g,h} a_g K_\psi^j(g,h) b_h$$ on the group algebra $\R[\G]$ of finitely supported real functions on $\G$.

\vskip5pt

\item[ii)] Let $\H_\psi^j$ be the Hilbert space completion of $$(\R[\G]/N_\psi^j, \langle \cdot \hskip1pt, \cdot \rangle_{\psi,j}) \quad \mbox{with} \quad N_\psi^j = \mbox{null space of} \ \langle \cdot \hskip1pt , \cdot \rangle_{\psi,j}.$$ If we consider the mapping $b_\psi^j: g \in \G \mapsto \delta_g + N_\psi^j \in \H_\psi^j$
\begin{eqnarray*}
\alpha_{\psi,g}^1 \Big( \sum_{h \in \G} a_h b_\psi^1(h) \Big) & = & \sum_{h \in \G} a_h \big( b_\psi^1(gh) - b_\psi^1(g) \big), \\ \alpha_{\psi,g}^2 \Big( \sum_{h \in \G} a_h b_\psi^2(h) \Big) & = & \sum_{h \in \G} a_h \big( b_\psi^2(hg^{-1}) - b_\psi^2(g^{-1}) \big),
\end{eqnarray*}
determine isometric actions $\alpha_\psi^j: \G \to \mathrm{Aut}(\H_\psi^j)$ of $\G$ on $\H_\psi^j$.

\vskip5pt

\item[iii)] Imposing the discrete topology on $\H_\psi^j$, the semidirect product $\G_\psi^j = \H_\psi^j \rtimes \G$ becomes a discrete group and we find the following group homomorphisms
\begin{eqnarray*}
\pi_\psi^1: g \in \G & \mapsto & b_\psi^1(g) \rtimes g \in \G_\psi^1, \\ \pi_\psi^2: g \in \G & \mapsto & b_\psi^2(g^{-1}) \rtimes g \in \G_\psi^2.
\end{eqnarray*}
\end{itemize}
\end{Alemma}

The previous lemma allows us to introduce two pseudo-metrics on our discrete group $\G$ in terms of the length function $\psi$. Indeed, a short calculation leads to the crucial identities $$\psi(g^{-1}h) \, = \, \Big\langle b_\psi^1(g) - b_\psi^1(h), b_\psi^1(g) - b_\psi^1(h) \Big\rangle_{\psi,1} \, = \, \big\| b_\psi^1(g) - b_\psi^1(h) \big\|_{\H_\psi^1}^2,$$ $$\psi(gh^{-1}) \, = \, \Big\langle b_\psi^2(g) - b_\psi^2(h), b_\psi^2(g) - b_\psi^2(h) \Big\rangle_{\psi,2} \, = \, \big\| b_\psi^2(g) - b_\psi^2(h) \big\|_{\H_\psi^2}^2.$$ In particular, $$\mathrm{dist}_1(g,h) = \sqrt{\psi(g^{-1}h)} = \|b_\psi^1(g) - b_\psi^1(h) \|_{\H_\psi^1}$$ defines a pseudo-metric on $\G$, which becomes a metric when the cocycle map is injective. Similarly, we may work with $\mathrm{dist}_2(g,h) = \sqrt{\psi(gh^{-1})}$. When the cocycle map is not injective the inverse image of $0$ $$\G_0 = \Big\{ g \in \G \, \big| \, \psi(g) = 0 \Big\}$$ defines a subgroup. The following elementary observation is relevant.

\begin{Alemma} \label{LemLeftRight}
Let $(\H_1, \alpha_1, b_1)$ and $(\H_2, \alpha_2, b_2)$ be a left and a right cocycle on $\G$. Assume that the associated length functions $\psi_{b_1}$ and $\psi_{b_2}$ coincide, then we find an isometric isomorphism $$\Lambda_{12}: b_1(g) \in \H_1 \mapsto b_2(g^{-1}) \in \H_2.$$ In particular, given a length function $\psi$ we see that $\H_{\psi}^1 \simeq \H_\psi^2$ via $b_\psi^1(g) \mapsto b_\psi^2(g^{-1})$.
\end{Alemma}

\dem By polarization, we see that $$\big\langle b_1(g), b_1(h) \big\rangle_{\H_1} = \frac{1}{2} \Big( \|b_1(g)\|_{\H_1}^2 + \|b_1(h)\|_{\H_1}^2- \|b_1(g) - b_1(h)\|_{\H_1}^2 \Big).$$ Since $b_1(g) - b_1(h) = \alpha_{1,h}(b_1(h^{-1}g))$, we obtain
\begin{eqnarray*}
\big\langle b_1(g), b_1(h) \big\rangle_{\H_1} & = & \frac{\psi_{b_1}(g) + \psi_{b_1}(h) - \psi_{b_1}(g^{-1}h)}{2} \\ & = & \frac{\psi_{b_2}(g) + \psi_{b_2}(h) - \psi_{b_2}(g^{-1}h)}{2} \ = \ \big\langle b_2(g^{-1}), b_2(h^{-1}) \big\rangle_{\H_2}.
\end{eqnarray*}
The last identity uses polarization and $b_2(g^{-1}) - b_2(h^{-1}) = \alpha_{2,h}(b_2(g^{-1}h))$.
\fin

\begin{Aremark}
\emph{According to Schoenberg's theorem, Markov semigroups of Fourier multipliers in $\V$ are in one-to-one correspondence with conditionally negative length functions $\psi: \G \to \R_+$. Lemma \ref{Lemmapsi} automatically gives 
\begin{eqnarray*}
\Gamma_\psi(f,f) & = & \frac12 \Big( A_\psi(f^*) f + f^* A_\psi(f) - A_\psi(f^*f) \Big) \\ & = & \sum_{g,h \in \G} \overline{\widehat{f}_1(g)} \widehat{f}_2(h) \, K_\psi(g,h) \, \lambda(g^{-1}h).
\ \ge \ 0.
\end{eqnarray*}}
\end{Aremark}

\begin{Atheorem} \label{CharacterizationCN}
Let $\Pi_0$ denote the space of trigonometric polynomials in $\V$ whose Fourier coefficients have vanishing sum. A given function $\psi: \G \to \R_+$ defines a conditionally negative length if and only if there exists a positive linear functional $\tau_\psi: \Pi_0 \to \C$ satisfying the identity $$\psi(g) \, = \, \tau_\psi \Big(2 \lambda(e) - \lambda(g) - \lambda(g^{-1}) \Big).$$
\end{Atheorem}

\dem Assume first that $\psi: \G \to \R_+$ satisfies the given identity for some positive linear functional $\tau_\psi: \Pi_0 \to \C$. To show that $\psi$ is a conditionally negative length it suffices to construct a cocycle $(\H_\psi, \alpha_\psi, b_\psi)$ so that $\psi(g) = \langle b_\psi(g), b_\psi(g) \rangle_{\H_\psi}$. Since $\Pi_0$ is a $*$-subalgebra, $\langle f_1, f_2 \rangle_{\H_\psi} = \tau_\psi(f_1^* f_2)$ is well-defined on $\Pi_0 \times \Pi_0$. If we quotient out the null space of this bracket, we may define $\H_\psi$ as the completion of such a quotient. If $N_\psi$ denotes the null space, let $$\alpha_{\psi,g}(f + N_\psi) = \lambda(g) f + N_\psi \quad \mbox{and} \quad b_\psi(g) = \lambda(g) - \lambda(e) + N_\psi.$$ It is easily checked that $(\H_\psi, \alpha_\psi, b_\psi)$ defines a left cocycle on $\G$. Moreover, since $2\lambda(e) - \lambda(g) - \lambda(g^{-1}) = |\lambda(g) - \lambda(e)|^2$, our assumption can be rewritten in the form $\psi(g) = \langle b_\psi(g), b_\psi(g) \rangle_{\H_\psi}$ as expected. Let us now prove the converse. Assume that $\psi: \G \to \R_+$ defines a conditionally negative length and define $$\tau_\psi (\lambda(g) - \lambda(e)) \, = \, - \frac12 \psi(g).$$ Since the polynomials $\lambda(g) - \lambda(e)$ span $\Pi_0$, $\tau_\psi$ extends to a linear functional on $\Pi_0$ which satisfies $\tau_\psi(2 \lambda(e) - \lambda(g) - \lambda(g^{-1}) = \frac12 (\psi(g) + \psi(g^{-1})) = \psi(g)$. Therefore, it just remains to show that $\tau_\psi: \Pi_0 \to \C$ is positive. Let $f = \sum_g a_g \lambda(g) \in \Pi_0$ so that $\sum_g a_g = 0$. In particular, we also have $f = \sum_g a_g (\lambda(g) - \lambda(e))$. By the conditional negativity of $\psi$ we find
\begin{eqnarray*}
\tau_\psi (|f|^2) & \!\!\! = \!\!\! & \sum_{g,h \in \G} \overline{a_g} a_h \tau_\psi \Big( \big( \lambda(g) - \lambda(e) \big)^* \big( \lambda(h) - \lambda(e) \big) \Big) \\ & \!\!\! = \!\!\! & \sum_{g,h \in \G} \overline{a_g} a_h \tau_\psi \Big( \lambda(g^{-1}h) - \lambda(g^{-1}) - \lambda(h) + \lambda(e) \Big) \\ & \!\!\! = \!\!\! & \sum_{g,h \in \G} \overline{a_g} a_h \frac{\psi(g^{-1}) + \psi(h) - \psi(g^{-1}h)}{2} \, = \, - \frac12 \sum_{g,h \in \G} \overline{a_g} a_h \psi(g^{-1}h) \, \ge \, 0.
\end{eqnarray*}
This shows that our functional $\tau_\psi: \Pi_0 \to \C$ is positive. The proof is complete. \fin

\begin{Aremark}
\emph{Given $\G$ discrete, consider a cocycle $(\H,\alpha,b)$. As we have seen in Lemma \ref{Lemmapsi}, it is convenient to impose the discrete topology on $\H$. If $\dim \H = n$ we sometimes make that clear writing $\R_{\mathrm{disc}}^n$ instead, even for $n=\infty$. As a discrete abelian group we find $\mathcal{L}(\R_{\mathrm{disc}}^n) \simeq L_\infty(\R_{\mathrm{bohr}}^n)$, the algebra of essentially bounded functions on the Bohr compactification of $\R^n$. Elements of this algebra have the form $$f \sim \summ_\xi \widehat{f}(\xi) \mbox{b-exp}_\xi$$ where $\mbox{b-exp}_\xi$ the $\xi$-th character on $\R^n_{\mathrm{bohr}}$ for each $\xi \in \R^n_{\mathrm{disc}}$. Its restriction to $\R^n$ is the standard character $\exp_\xi(x) = \exp(2 \pi i \langle \xi, x \rangle)$. Note that this construction is still meaningful for infinite dimensions, where as the Lebesgue measure does not extend to this case. According to de Leeuw's compactification theorem \cite{dL}, a Fourier multiplier is bounded on $L_p(\R^n)$ iff it is also bounded in $L_p(\R_{\mathrm{bohr}}^n)$.}
\end{Aremark}

\section*{\textbf{Appendix B. A geometric perspective}}

In this appendix we will describe the tangent module associated with a given length function, and how it can be combined with Riesz transform estimates. Let us recall that a Hilbert module over an algebra $\A$ is a vector space $\mathrm{X}$ with a bilinear map $m: (\xi,a) \in \mathrm{X} \times \A \to \xi a \in \mathrm{X}$ and a sesquilinear form $\langle \ ,\ \rangle: \mathrm{X} \times \mathrm{X} \to \A$ such that $\langle \xi,\eta a \rangle = \langle \xi, \eta \rangle a$, $\langle \xi a,\eta \rangle = a^*\langle \xi,\eta \rangle$ and $\langle \xi,\xi \rangle \ge 0$. We refer to Lance's book \cite{La} for more information. Define $\mathcal{L}(\mathrm{X})$ as the $\mathrm{C}^*$-algebra of right-module maps $T$ which admit an adjoint. That is, there exits a linear map $S: \mathrm{X} \to \mathrm{X}$ such that $\langle S \xi,\eta \rangle=\langle \xi,T\eta\rangle$. A Hilbert bimodule is additionally equipped with a $*$-homomorphism $\pi: \A \to \mathcal{L}(\mathrm{X})$ and a \emph{derivation} $\delta: \A \to \mathrm{X}$ is a linear map which satisfies the Leibniz rule $$\delta(ab) = \pi(a)\delta(b)+\delta(a)b.$$ A typical example for such a derivation is given by an inclusion $\A \subset \mathcal{M}$ into some von Neumann algebra $\mathcal{M}$, a conditional expectation $E: \mathcal{M} \to \A''$, and a  vector $\xi$ such that $\delta(a) = a\xi-\xi a$. We have seen above that for a conditionally negative length function $\psi$, we can construct an associated (left) cocycle $(\H_{\psi}, \alpha_\psi, b_\psi)$. In the following, we will assume that the $\mathbb{R}$-linear span of $b_\psi(\G)$ is $\H_{\psi}$. Let us recall the Brownian functor $B: \H_{\psi}\to L_2(\Omega)$ which comes with an extended action $\alpha$ of $\G$ on $L_{\infty}(\Omega)$. This construction is usually called the gaussian measure space construction. Here the derivation is given by $$\delta_\psi: \lambda(g) \mapsto B(b_\psi(g)) \rtimes \lambda(g).$$ We have already encountered the corresponding bimodule in context with the Khintchine inequality. The proof of the following lemma is obvious. In fact, here left and right action coincide.

\begin{Blemma} 
The Hilbert bimodule $$\Omega_{\psi}(\G) = \delta_\psi(\mathbb{C}[\G]) (\mathbf{1} \rtimes \lambda(\G))$$ is exactly given by the vector space $\mathrm{X}_{\psi} = \big\{ \sum_g B(\xi_{g}) \rtimes \lambda(g): \xi_{g} \in \H_{\psi} \big\}$.
\end{Blemma}

\dem For $\xi=b_\psi(g)$ we consider $\delta_\psi(\lambda(g)) (\mathbf{1} \rtimes \lambda(g^{-1})) = B(\xi)$. Since $\H_{\psi}$ is the real linear span, we deduce that $\mathrm{X}_{\psi}$ is contained in $\delta_	\psi(\mathbb{C}[\G])(\mathbf{1} \rtimes \lambda(\G))$. The converse is obvious. Moreover, since $\delta_\psi$ is a derivation, it is easy to see that $\delta(\mathbb{C}[\G]) (\mathbf{1} \rtimes \lambda(\G))$ is invariant under the left action. \fin

Of course, the gaussian functor $B$ is not really necessary to describe the bimodule $\Omega_{\psi}(\G) \simeq \H_{\psi} \rtimes \G$. The following proposition shows that our previous results extend to the tangent module and not only to differential forms with `constant' coefficients given by elements in $\H_{\psi} \subset \Omega_{\psi}(\G)$. Given $\xi=\sum_{h} B(\xi_h) \rtimes \lambda(h) \in \Omega_{\psi}(\G)$, define the extended Riesz transform in the direction of $\xi$ as follows 
\begin{eqnarray*}
\mathbf{R}_{\psi, \xi} f & = & \sum_{h \in \G} \lambda(h^{-1}) R_{\psi, \xi_h} f \\ [3pt] & = & 2 \pi i \sum_{g,h \in \G} \frac{\langle b_\psi(g), \xi_h \rangle_{\H_\psi}}{\sqrt{\psi(g)}} \widehat{f}(g) \lambda(h^{-1}g) \\ & = & \mathsf{E}_\V \Big[ \Big( \sum_{h \in \G} B(\xi_h) \rtimes \lambda(h) \Big)^* \Big( \sum_{g \in \G} \frac{\widehat{f}(g)}{\sqrt{\psi(g)}} B(b_\psi(g)) \rtimes \lambda(g) \Big) \Big] \\ [5pt] & = & \mathsf{E}_\V \big( \xi^* \delta_\psi A_\psi^{-\frac12} f \big).
\end{eqnarray*}
Note that we recover the Riesz transforms $R_{\psi,h}$ for $\xi = B(h) \rtimes \lambda(e)$.

\begin{Bproposition} If $2\le p<\infty$
\begin{itemize}
\item[i)] $\big\| \mathbf{R}_{\psi,\xi}: L_p(\widehat{\mathbf{G}}) \to L_p(\widehat{\mathbf{G}}) \big\| \, \lesssim_{c(p)} \, \|\xi\|_{\Omega_{\psi}(\G)}$,

\vskip9pt

\item[ii)] $\displaystyle \Big\| \Big( \sum_{j \ge 1} |\mathbf{R}_{\psi,\xi_j}(f)|^2 \Big)^{\frac12} \Big\|_p \, \lesssim_{c(p)} \, \Big\| \Big( \langle \xi_j,\xi_k \rangle \Big) \Big\|_{\mathcal{B}(\ell_2) \bar{\otimes} \V}^{\frac12} \|f\|_p$,

\item[iii)] $\displaystyle \Big\| \Big( \sum_{j \ge 1} |\mathbf{R}_{\psi,\xi_j}(f_j)|^2 \Big)^{\frac12} \Big\|_p \, \lesssim_{c(p)} \, \sup_{j \ge 1} \big\| \mathsf{E}_\V(\xi_j^* \xi_j) \big\|_{\V}^{\frac12} \Big\| \Big( \sum_{j \ge 1} |f_j|^2 \Big)^{\frac12} \Big\|_p$.
\end{itemize}
\end{Bproposition}

\dem Assertion i) follows from ii) or iii). The second and third assertions follow from the well-known Cauchy-Schwarz inequality for conditional expectations, whose noncommutative form follows from Hilbert module theory \cite{JDoob} $$\big\| \mathsf{E}_\V \otimes id_{\mathcal{B}(\ell_2)}(xy) \big\|_p \le \big\| \mathsf{E}_\V \otimes id_{\mathcal{B}(\ell_2)} (xx^*) \big\|_{\infty}^{\frac12} \big\| \big( \mathsf{E}_\V \otimes id_{\mathcal{B}(\ell_2)} (y^*y) \big)^{\frac12} \big\|_p.$$ Indeed, for ii) we pick $$x = \sum_{j \ge 1} \xi_j^* \otimes e_{j1} \quad \mbox{and} \quad y = \delta_\psi A_\psi^{-\frac12} f \otimes e_{11}.$$ On the other hand, for iii) we take $$x = \sum_{j \ge 1} \xi_j^* \otimes e_{jj} \quad \mbox{and} \quad y = \sum_{j \ge 1} \delta_\psi A_\psi^{-\frac12} f_j \otimes e_{j1}.$$ After it, the result follows from the cb-extension of Theorem A2 in Remark \ref{OS}. \fin

\begin{Bremark}
{\rm The `adjoint' of $\mathbf{R}_{\psi,\xi}$ given by $$\mathbf{R}_{\psi,\xi}^{\dagger}(f) = \mathbf{R}_{\psi,\xi}(f^*)^*$$ can be written as $\mathbf{R}_{\psi,\xi}^{\dagger}(f) = \mathsf{E}_\V(\xi^* \delta_\psi A_\psi^{-\frac12}(f^*))^* = - \mathsf{E}_\V(\delta_\psi A_\psi^{-\frac12} (f) \xi)$. Hence $$\Big\| \Big( \sum_{j \ge 1} | \mathbf{R}_{\psi,\xi_j}^{\dagger}(f_j)^*|^2 \Big) \Big\|_p
\, \lesssim_{c(p)} \, \sup_{j \ge 1} \big\| \mathsf{E}_\V (\xi_j^* \xi_j) \big\|_\V^{\frac12} \Big\| \Big( \sum_{j \ge 1} |f_j^*|^2 \Big)^{\frac12} \Big\|_p.$$ It is however, in general difficult to find an element $\eta$ such that $\mathbf{R}_{\psi,\eta}^{\dagger}(f) = \mathbf{R}_{\psi,\xi}(f)$ unless $\G$ is commutative and the action is trivial. This is a particular challenge if we want to extend the results from above literally to $p<2$, because then we need both a row and a column bound to accommodate the decomposition $R(f)=a+b$ in the tangent module $\Omega_{\psi}(\G)$.}
\end{Bremark}

Let us now indicate how to construct the corresponding real spectral triple. We first recall that $\Omega_{\psi}(\G)$ is a quotient of the universal object $\Omega_{\bullet}(\G) \subset \mathbb{C}[\G] \otimes \mathbb{C}[\G]$ spanned by $\delta_{\bullet}(a)b$, where the universal derivation is $\delta_{\bullet}(a)=a \otimes \mathbf{1} - \mathbf{1} \otimes a$. Using $J(a \otimes b) = b^* \otimes a^*$, the universal object becomes a real bimodule. In other words the left and right representations $\pi_\ell(a)(\xi)=(a \otimes \mathbf{1}) \xi$ and $\pi_r(a)(\xi)=\xi(\mathbf{1} \otimes a)$ are related via $J\pi_r(a^*)J=\pi_\ell(a)$. In our concrete situation, we have a natural isometry $J(x)=x^*$ on $\mathrm{M} = L_{\infty}(\Omega) \rtimes \G$, which leaves the subspace $L_2(\Omega_{\psi}(\G)) \subset L_2(\mathrm{M})$ invariant. Hence $\Omega_{\psi}(\G)$ is a quotient of $\Omega_{\bullet}(\G)$. The Dirac operator for this spectral triple is easy to construct. The underlying Hilbert space is $\H = L_2(\widehat{\mathbf{G}}) \oplus L_2(\Omega_{\psi}(\G))$ and $$D = \Big( \begin{array}{lc} 0 & \delta_\psi \\ \delta_\psi^*& 0 \end{array} \Big).$$ Note that $\delta_\psi^*(B(\xi) \rtimes \lambda(g)) = \langle \xi, b_\psi(g) \rangle \lambda(g)$ is densely defined. Now, in order to extend $D$ to a larger Hilbert space, we must replace the gaussian construction by the semicircular construction. As in the gaussian category, we have a function $$s: \H \to \Gamma_0(\H)$$ into the von Neumann algebra generated by semicircular random variables, and a representation $\alpha: O(\H) \to {\rm Aut}(\Gamma_0(\H))$ with $s(o(h)) = \alpha_o(s(h))$. This allows us to define $\delta_{free}^\psi(\lambda(g)) = s(b_\psi(g)) \rtimes \lambda(g)$. The proof for the boundedness of the corresponding Riesz transforms $$f \mapsto \delta_{free}^\psi A_\psi^{-\frac12} f$$ follows from the corresponding Khintchine inequality for $x = \sum_{\xi,g} s(\xi) \rtimes \lambda(g)$. Namely, we find $$\|x\|_{p} \, \sim_{\sqrt{p}} \, \max \Big\{ \big\| \mathsf{E}_\V (x^*x)^{\frac12} \big\|_p, \big\|\mathsf{E}_\V(xx^*)^{\frac12} \big\|_p \Big\}$$ for $2 \le p < \infty$ and $$\|x\|_{p'} \, \sim_{\sqrt{p}} \, \inf_{x = x_1+x_2} \big\|\mathsf{E}_\V(x_1^*x_1)^{\frac12} \big\|_p + \big\| \mathsf{E}_\V(x_2 x_2^*)^{\frac12} \big\|_p.$$ In fact, it turns out that $$\mathsf{E}_\V(|x_{free}|^2) = \mathsf{E}_\V(|x_{gauss}|^2)$$ for $x_{free} = \sum_{\xi,g} s(\xi) \rtimes \lambda(g)$ and $x_{gauss} = \sum_{\xi,g} B(\xi) \rtimes \lambda(g)$. This means the bimodule $\mathrm{X}_\psi$ can be realized either with gaussian or semicircular variables.

\begin{Bproposition}
The tuple $$(\mathbb{C}[\G], L_2(\Gamma_0(H) \rtimes \G),D_{free},J)$$ forms a real spectral triple. Given $f \in \mathbb{C}[\G]$, we have $$\big\| [D,f] \big\| \sim \max \Big\{ \big\| \Gamma(f,f) \big\|_{\V}^{\frac12}, \big\| \Gamma(f^*,f^*) \big\|_{\V}^{\frac12} \Big\}.$$ 
\end{Bproposition}

\dem Using the inclusion $\V \subset L_2(\widehat{\mathbf{G}})$ via $f \mapsto f 1_{\tau}$ in the GNS construction, the definition of the crossed product and recalling that $\delta_{free}^\psi(\lambda(g))=s(b_\psi(g)) \rtimes \lambda(g)$, we get
\begin{eqnarray*}
\delta_{free}^\psi(e_g) \!\!\! & = & \!\!\! \delta_{free}^\psi(\lambda(g) 1_{\tau}) = \delta_{free}^\psi (\lambda(g)) 1_{\tau} \\ \!\!\! & = & \!\!\! s(b_\psi(g))(e_g) = \alpha_{g^{-1}}(s(b_\psi(g))) \otimes e_g = -s(b_\psi(g^{-1})) \otimes e_g.
\end{eqnarray*}
Let $h \in \G$, then we deduce
\begin{eqnarray*}
\lefteqn{\hskip-10pt (\lambda(h) \delta_{free}^\psi - \delta_{free}^\psi \lambda(h))e_g} \\ & = & -s(b_\psi(g^{-1}))e_{hg} + s(b_\psi(g^{-1}h^{-1}))e_{hg} = -\lambda(h)[\alpha_{g^{-1}}s(b_\psi(h^{-1})](e_g).
\end{eqnarray*}
Since $\delta_{free}^\psi(\lambda(h))=s(b_\psi(h)) \rtimes \lambda(h) = -(0 \rtimes \lambda(h)) (s(b_\psi(h^{-1})) \rtimes \mathbf{1})$, we get $$[\lambda(h),\delta_{free}^\psi] = \delta_{free}^\psi(\lambda(h)).$$ By linearity, we deduce that $[f,\delta_{free}^\psi] = \delta_{free}^\psi(f)$. Now, we shall use a central limit procedure. Consider the crossed product $\Gamma_0(\ell_2^m(H))\rtimes \G$. Then, the copies $\pi_j(\Gamma_0(H)\rtimes \G)$ given by the $j$-th coordinate are freely independent over $\V$. Thus Voiculescu's inequality from \cite{JAW} applies and yields that for $x=\sum_{\xi,g} a_{\xi,g}s(\xi)\rtimes \lambda(g)$ and the sum of independent copies $$\Big\| \sum_{j=1}^m a_{\xi,g} s(\xi \otimes e_j) \lambda(g) \Big\|_{\Gamma_0 \rtimes \G} \le \|x\| + \sqrt{m} \|\mathsf{E}_\V(x^*x)\|_\V^{\frac12} + \sqrt{m} \|\mathsf{E}_\V(xx^*)\|_\V^{\frac12}.$$ Dividing by $\sqrt{m}$ and observing that $x$ and $\frac{1}{\sqrt{m}} \sum_{j=1}^m a_{\xi,g} s(\xi \otimes e_j) \lambda(g)$ are equal in distribution we find indeed, $$\|x\|_{\Gamma_0 \rtimes \G} \, \sim \, \max \Big\{ \big\| \mathsf{E}_\V(x^*x) \big\|_\V^{\frac12}, \big\| \mathsf{E}_\V(xx^*) \big\|_\V^{\frac12} \Big\}.$$ Thus for a differential form $\omega \in \Omega_{free}^\psi(\G)$, we get $$\|\omega\| \sim \max\{\|\omega\|_{\mathrm{X}_\psi}, \|\omega^*\|_{\mathrm{X}_\psi}\}.$$ In particular, we conclude that $$\|\delta_{free}^\psi(f)\| \, \sim \, \max \Big\{ \big\| \Gamma(f,f) \big\|_\V^{\frac12}, \big\| \Gamma(f^*,f^*) \big\|_\V^{\frac12} \Big\}.$$ This expression is certainly finite for $f \in \mathbb{C}[\G]$ and the proof is complete. \fin

\begin{Bremark}
{\rm  If we are only interested in the norm of the commutator $$\big\| [D,f] \big\|_{\mathcal{B}(\H)} \, = \, \max \Big\{ \big\| \Gamma(f,f) \big\|_\V^{\frac12}, \big\| \Gamma(f^*,f^*) \big\|_\V^{\frac12} \Big\}$$ for the Hilbert space $\mathcal{H}$ from above, there is no need to resort to free probability. Indeed, we have $[\delta_\psi,\lambda(g)] = L_{b_\psi(g)} \lambda(g)$, where $L_{b_\psi(g)} =b_\psi(g) \otimes x$ is the usual rank one operator from Hilbert $\mathrm{C}^*$-modules. It then follows that $$\big\| [D,f] \big\|_{	\mathcal{B}(\H)} \, = \, \max \Big\{ \big\| \Gamma(f,f) \big\|_\V^{\frac12}, \big\| \Gamma(f^*,f^*) \big\|_\V^{\frac12} \Big\}$$
gives precisely the Lip-norm considered in \cite{JM} and found in the free tangent module.}
\end{Bremark}

As in \cite{CGIS}, we have to deal with the fact that this spectral triple might have some degenerated parts, but in many calculations of the $\zeta$-function of $|D|$ the kernel is usually ignored. In fact, if $D(e_g)=\frac{1}{\psi(g)}e_g$ defines a compact operator, then $[F,\pi(a)]$ is compact for all $a\in \mathbb{C}[\G]$. This follows from the boundedness of $[D,a]$ \cite[Proposition 2.5]{CGIS} and then \cite[Proposition 2.8]{CGIS} applies. In our situation, we know that $\delta_\psi = R_\psi A_\psi^{1/2}$ and hence $\delta_\psi$ vanishes on $\H_0={\rm  span} \{e_g: \psi(g)=0\}$. Clearly, we see that $\delta_\psi^*\delta_\psi = A_\psi$ is the generator of our semigroup which also vanishes on $\H_0$. On the hand, we have $$\delta_\psi \delta_\psi^* = R_\psi A_\psi R_\psi^*$$ and hence the range of $\delta_\psi \delta_\psi^*$ is given by the first Hodge projection $\Pi_{hdg} = R_\psi R_\psi^*$. This can be described explicitly. Indeed, for $g$ with $\psi(g) \neq 0$ we denote by $Q_g$ the projection onto the span of $B(b_\psi(g))\in L_2(\Omega)$ and get $$R_\psi R_\psi^* = \sum_g Q_g \otimes e_{gg}.$$ With the help of the Hodge projection, we may now consider the signature of the Dirac operator $$F = D(D^*D)^{-1/2} \Pi_{hdg}.$$ 

\begin{Bproblem} 
\emph{Show that $F$ admits dimension free estimates in $L_p$.}
\end{Bproblem} 

\section*{\textbf{Appendix C. Meyer's problem for Poisson}}

Let  $\Delta=\partial_x^2$ be the Laplacian operator on ${\Bbb R}^n$. Classical theory of semigroups of operators yields that the fractional laplacians $(-\Delta)^\beta$ with $0 < \beta < 1$ are closed densely defined operators on $L_p({\Bbb R}^n)$ \cite[Chapter 9, Section 11]{Yosida} and regarding them as Fourier multipliers we see that the Schwartz class lies in the domain of any of them. Moreover, they generate Markov semigroups on $L_\infty({\Bbb R}^n)$. When $\beta = \frac12$ we get the Poisson semigroups $P_t = \exp(-t\sqrt{-\Delta})$. In this appendix we shall show that Meyer's problem \eqref{MeyerProblem} fails for this generator when $p \le \frac{2n}{n+1}$. Let us first give a formula for the corresponding carr\'e du champs $\Gamma_{\frac12}$. 

\begin{Clemma} For any  Schwartz function $f$, we have
\[ \Gamma_{\frac12}(f,f) = \int_0^{\infty}  P_{t} |\nabla P_{t}f|^2 dt \]
where $\nabla g(x,t)=(\partial_{x_1}f, \partial_{x_2}f,\ldots, \partial_{x_n}f,\partial_tf)$ includes spatial and time variables.
\end{Clemma}

\dem Let
\[
F_t={(\partial_t P_t)}(|P_tf|^2)-P_t((\partial_t P_tf^{*})P_tf)-P_t(P_tf^{*}(\partial_t P_tf)).
\]
Then, since $\partial_t^2 P_t + \Delta P_t = 0$ we deduce
\begin{eqnarray*}
\partial_t F_t &=& \partial_t^2 P_t (|P_tf|^2)-P_t(\partial_t^2 P_tf^{*} P_tf)-P_t(P_tf^{*}%
\partial_t^2 P_tf)-2P_t(|\partial_t P_tf|^2) \\
&=&-\Delta P_t(|P_tf|^2)-P_t((-\Delta P_tf^{*})P_tf)-P_t(P_tf^{*}(-\Delta P_tf))-2P_t(|\partial_t P_tf|^2) \\
&=&-2P_t(|\partial_x P_tf|^2)-2P_t(|P_t^{\prime }f|^2)=-2P_t|\nabla P_tf|^2.
\end{eqnarray*}
Note that $F_0=\lim_{t\rightarrow 0}F_t=2\Gamma _{\frac 12}(f,f)$ by the definition of Carre du Champ and $F_t\rightarrow 0$  as $
t\rightarrow \infty$. We get
\begin{align*}
2\Gamma _{\frac 12}(f,f) &=\int_0^\infty -\partial_t
F_tdt  =2\int_0^{\infty}
 P_{t} |\nabla P_{t}f|^2 dt . \qedhere
\end{align*}
\fin

\begin{Cproposition}  The equivalence \eqref{MeyerProblem} fails for Poisson semigroups $P_t=e^{-t\sqrt{-\Delta}}$ on ${\Bbb R}^n$ for any $1 < p \le \frac {2n}{n+1}$ with $n\geq 2$. More precisely, for any non-zero Schwartz function $f$  
we have $$\Gamma_{\frac12}(f,f)^\frac12 \notin L_p(\R^n) \quad \mbox{for any} \quad p \le \frac {2n}{n+1}.$$ 
\end{Cproposition}

\dem We follow an argument from \cite{Fef0}. Fix a non-zero $f \in L_p(\R^n)$ and $|x|>4$. Then
\begin{eqnarray*}
\Gamma_{\frac12}(f,f)& =& \int_0^{\infty}
 P_{t} |\nabla P_{t}f|^2 dt\\
 & =& \int_0^{\infty}\int_{{\Bbb R}^n} c_n\frac {t}{(|x-y|^2+t^2)^{\frac {n+1}2}}
  |\nabla P_{t}f(y)|^2 dydt  \\
  &\geq&c_n\int_1^{2}\int_{|y|<1} \frac {1}{|x|^{n+1}}
  |\nabla P_{t}f(y)|^2 dydt =c_n c_f \frac {1}{|x|^{n+1}}
\end{eqnarray*}
for $c_f = \int_1^2\int_{|y|<1} |\nabla P_{t}f(y)|^2dy dt > 0$ (since $f \neq 0$) and any $|x| > 4$. Then $$\Gamma_{\frac12}(f,f)^\frac12 \in L_p(\R^n) \ \Rightarrow \ \frac p2 (n+1) > n \ \Rightarrow \ p >\frac {2n}{n+1}.$$
We then conclude that  $$\big\| \Gamma_{\frac12}(f,f)^\frac12 \big\|_p = \infty \quad \mbox{while} \quad \big\| (-\Delta )^\frac14f \big\|_p<\infty$$ for any non-zero Schwartz function $f$ with $p \le \frac {2n}{n+1}$. Therefore, \eqref{MeyerProblem} fails. \fin

\noindent Thus, our revision of \eqref{MeyerProblem} in this paper is needed even for commutative semigroups.

\vskip3pt

\noindent \textbf{Acknowledgement.} We are indebted to Anthony Carbery, Andreas Seeger and Jim Wright for interesting comments and references. Junge is partially supported by the NSF DMS-1201886, Mei by the NSF DMS-1266042 and Parcet by the ERC StG-256997-CZOSQP and the Spanish grant MTM2010-16518. Junge and Parcet are also supported in part by ICMAT Severo Ochoa Grant SEV-2011-0087 (Spain).

\bibliographystyle{amsplain}

\newpage

\null

\vskip30pt

\hfill \noindent \textbf{Marius Junge} \\
\null \hfill Department of Mathematics
\\ \null \hfill University of Illinois at Urbana-Champaign \\
\null \hfill 1409 W. Green St. Urbana, IL 61891. USA \\
\null \hfill\texttt{junge@math.uiuc.edu}

\vskip2pt

\hfill \noindent \textbf{Tao Mei} \\
\null \hfill Department of Mathematics
\\ \null \hfill Wayne State University \\
\null \hfill 656 W. Kirby Detroit, MI 48202. USA \\
\null \hfill\texttt{mei@wayne.edu}

\vskip2pt

\enlargethispage{2cm}

\hfill \noindent \textbf{Javier Parcet} \\
\null \hfill Instituto de Ciencias Matem{\'a}ticas \\ \null \hfill
CSIC-UAM-UC3M-UCM \\ \null \hfill Consejo Superior de
Investigaciones Cient{\'\i}ficas \\ \null \hfill C/ Nicol\'as Cabrera 13-15.
28049, Madrid. Spain \\ \null \hfill\texttt{javier.parcet@icmat.es}
\end{document}